\documentclass{article}
\usepackage[T1]{fontenc}
\usepackage{lmodern} 
\usepackage[a4paper]{geometry} 
\usepackage{fullpage} 

\usepackage[french]{babel}
\usepackage{amsmath,amsfonts,amssymb,stmaryrd}
\usepackage{enumerate}
\usepackage{amscd}
\usepackage{amsthm}
\usepackage[all]{xy}
\usepackage{dsfont}

\newtheorem{déf}{Définition}[subsection] 

\newtheorem{propo}[déf]{Proposition}
\newtheorem{theo}[déf]{Théorème}
\newtheorem{cor}[déf]{Corollaire}
\newtheorem{lem}[déf]{Lemme}
\newtheorem*{theoNonNum}{Théorème} 
\newtheorem*{lemNonNum}{Lemme} 

\theoremstyle{definition}
\newtheorem{rmq}[déf]{Remarque}

\theoremstyle{definition}

\theoremstyle{definition}
\newtheorem{nota}[déf]{Notation}

\newenvironment{démo}
{\noindent{\textit{Démonstration}.}}{\hfill$\square$}

\begin{document}

\title{Correspondance de Jacquet-Langlands et distinction: cas de certaines séries 
discrètes non cuspidales de niveau $0$}
\author{Charlène Coniglio-Guilloton}

\maketitle

\begin{center}
\textbf{Résumé :}
\end{center}
On considère $\mathbb{K} / \mathbb{F}$ une extension quadratique séparable modérément ramifiée de 
corps locaux non archimédiens. 
Soit $\mathcal{D}$ une $\mathbb{F}$-algèbre à division centrale d'indice $n$ pair. 
Soit $G = (\mathcal{D} \otimes_{\mathbb{F}} \mathbb{K})^{\times} \simeq {\rm GL}_2 (\Delta)$ où $\Delta$ 
est une $\mathbb{K}$-algèbre à division centrale d'indice $n/2$. 
Alors $G/ \mathcal{D}^{\times}$ est un espace symétrique. 
On considère $(\pi, V)$ une représentation lisse et irréductible de niveau $0$, non cuspidale, 
membre de la série discrète de $G$. On obtient des conditions de $\mathcal{D}^{\times}$-distinction pour 
$\pi$. Pour cela, on utilise la paramétrisation de Silberger et Zink de $\pi$ via une paire admissible modérée ainsi 
que le système de coefficients sur l'immeuble de $G$ associé à $\pi$ donné par Schneider et 
Stuhler. Puis, dans les cas étudiés, on vérifie que 
la correspondance de Jacquet-Langlands préserve la distinction.

\begin{center}
\textbf{Abstract :}
\end{center}
Let $\mathbb{K} / \mathbb{F}$ be a quadratic separable extension of non archimedean local fields. 
Let $\mathcal{D}$ be a division algebra of center $\mathbb{F}$ and dimension $n^2$ over $\mathbb{F}$. 
One assumes that $n$ is even. Let $G = (\mathcal{D} \otimes_{\mathbb{F}} \mathbb{K})^{\times} \simeq {\rm GL}_2 (\Delta)$ 
where $\Delta$ is a division algebra of center $\mathbb{K}$ and of dimension $(n/2)^2$ over $\mathbb{K}$. 
Then $G/ \mathcal{D}^{\times}$ is a symmetric space. 
One considers $(\pi, V)$ a level zero discrete series representation of $G$. One assumes that $(\pi, V)$ is not 
supercuspidal. Using the parametrization of Silberger and Zink of $\pi$ in terms of tame admissible pairs and the coefficient 
system of the Bruhat-Tits building associated to $\pi$ by Schneider and Stuhler, we obtain some conditions of 
$\mathcal{D}^{\times}$-distinguishness of $\pi$. Then, in certain cases, we prove that the Jacquet-Langlands 
correspondence preserves distinction.

\section*{Introduction}

\noindent
Nous noterons $\mathbb{K} / \mathbb{F}$ une extension quadratique séparable modérément ramifiée 
de corps locaux non archimédiens. 
Soit $n$ un entier naturel non nul. On note $G_{\mathbb{F}} = {\rm GL}_n (\mathbb{F})$ 
et $G_{\mathbb{K}} = {\rm GL}_n (\mathbb{K})$. On fixe un diviseur $d$ de $n$ 
($dm = n$) et $\mathcal{D}$ une $\mathbb{F}$-algèbre à division centrale d'indice $d$ 
(i.e de dimension $d^2$ sur son centre $\mathbb{F}$). Enfin, on note 
$H_{\mathbb{F}} = {\rm GL}_m (\mathcal{D})$ et 
$H_{\mathbb{K}} = ({\rm M}_m (\mathcal{D}) \otimes_{\mathbb{F}} \mathbb{K})^{\times}$. 
Il existe un diviseur $\mu$ de $n$ et une $\mathbb{K}$-algèbre à division centrale $\Delta$ 
d'indice $\delta = n / \mu$ tels que $H_{\mathbb{K}} = {\rm GL}_{\mu} (\Delta)$. 
On a le résultat suivant :

\begin{theoNonNum}
(\cite{JacquetLanglands}, \cite{Rogawski}, \cite{DeligneKazhdanVigneras}, \cite{Badulescu})
Il existe une unique bijection, appelée correspondance de Jacquet-Langlands :
$$
JL : \mathcal{R}^2 (G_{\mathbb{K}}) \rightarrow \mathcal{R}^2 (H_{\mathbb{K}})
$$
telle que pour tous $(g, \widetilde{g}) \in H_{\mathbb{K}} \times G_{\mathbb{K}}$ 
elliptiques réguliers de même polynôme minimal, et toute représentation 
$\pi \in \mathcal{R}^2 (G_{\mathbb{K}})$, on a:
$$
\Theta_{\pi} (\widetilde{g}) 
= (-1)^{\mu \times (\delta-1)} \Theta_{JL (\pi)} (g)
$$
où $\mathcal{R}^2 (G_{\mathbb{K}})$ (resp. $\mathcal{R}^2 (H_{\mathbb{K}})$) sont les  
 classes d'isomorphisme des représentations lisses irréductibles membres de la série discrète 
de $G_{\mathbb{K}}$ (resp. $H_{\mathbb{K}}$) et $\Theta$ désigne le caractère d'Harish Chandra.
\end{theoNonNum}

\noindent
Lorsque $n = 2$ et $d=2$, la correspondance de Jacquet-Langlands est l'application identité. 
Dans ce cas, les travaux de J. Hakim et de D. Prasad nous montrent qu'une série discrète (en niveau quelconque) de 
${\rm GL}_2 (\mathbb{K})$ 
est ${\rm GL}_2 (\mathbb{F})$-distinguée si et seulement si elle est 
$\mathcal{D}^{\times}$-distinguée (on pourra se référer dans \cite{Hakim} au 
Théorème 9.1 page 21 pour les représentations cuspidales de caractère 
central trivial ainsi qu'au Théorème 7.1 page 16 pour la représentation de Steinberg et ses tordues 
ou on pourra retrouver ce résultat dans \cite{Prasad2} Théorème C). 
Plus généralement, les travaux de J. Hakim et F. Murnaghan dans \cite{HakimMurnaghan1} 
(Théorème 11.1 page 1887) nous donnent des critères de 
${\rm GL}_n (\mathbb{F})$-distinction pour les cuspidales modérées de ${\rm GL}_n (\mathbb{K})$ (là aussi en niveau quelconque). 
Enfin, dans \cite{Coniglio}, nous avons montré qu'une représentation cuspidale de niveau $0$ de 
${\rm GL}_n (\mathbb{K})$ est ${\rm GL}_n (\mathbb{F})$-distinguée si et seulement si son image par la 
correspondance de Jacquet-Langlands est ${\rm GL}_m (\mathcal{D})$-distinguée.\\

\noindent
Dans le travail qui suit, nous généralisons les résultats évoqués précédemment 
dans le cas de certaines représentations lisses 
irréductibles de $H_{\mathbb{K}} = {\rm GL}_{\mu} (\Delta)$, de niveau $0$, 
membres de la série discrète, non cuspidales, 
dans le cas 
particulier où $\mu = 2$. 
Avec les articles \cite{SilbergerZink1} et \cite{SilbergerZink2}, 
A. Silberger et E. W. Zink montrent que la correspondance 
de Jacquet-Langlands se restreint en une bijection :
$$
JL: \mathcal{R}_0^2 ({\rm GL}_n (\mathbb{K})) \rightarrow \mathcal{R}_0^2 ({\rm GL}_2 (\Delta))
$$
où $\mathcal{R}_0^2 ({\rm GL}_n (\mathbb{K}))$ (resp. $\mathcal{R}_0^2 ({\rm GL}_2 (\Delta))$) 
sont les classes d'isomorphisme des 
représentations lisses irréductibles membres de la série discrète de niveau $0$.
De plus, les articles \cite{SilbergerZink1} et \cite{SilbergerZink2} nous donnent une paramétrisation de ces séries discrètes 
de niveau $0$ via des paires admissibles modérées. 
Notre objectif est de déterminer des conditions de $\mathcal{D}^{\times}$-distinction de 
telles représentations (où $\mathcal{D}$ est une $\mathbb{F}$-algèbre à division centrale d'indice $n = d = 2 \delta$). 
Nous n'obtiendrons pas de résultats généraux mais, dans les cas que nous étudions,  
nous montrerons que  la correspondance de Jacquet-Langlands préserve la distinction. 
Notre méthode de travail est la suivante : soit $(\pi, V)$ une série discrète non cuspidale de niveau $0$ de 
${\rm GL}_2 (\Delta)$. 
On note 
$(\chi_f, \mathbb{K}_f)$ 
(où $\chi_f$ est en particulier un caractère modéré $\mathbb{K}$-régulier de $\mathbb{K}_f^{\times}$) 
la paire admissible modérée associée à $\pi$ par \cite{SilbergerZink2}. 
Notons $X_{\mathbb{K}}$ l'immeuble de Bruhat-Tits de ${\rm GL}_2 (\Delta)$,  
$X_{\mathbb{F}}$ celui de $\mathcal{D}^{\times}$ (qui est réduit à un point) 
et $j : X_{\mathbb{F}} \mapsto X_{\mathbb{K}}$ l'injection naturelle entre ces immeubles. 
On notera $X_0$ l'ensemble des sommets de $X_{\mathbb{K}}$ et $X_{(1)}$ l'ensemble des arêtes orientées. 
Si $s$ est un sommet de $X_{\mathbb{K}}$, on note $\mathcal{A}_s$ l'ordre héréditaire associé à $s$, 
$\mathcal{P}_s$ le radical de Jacobson de $\mathcal{A}_s$, $\mathcal{U}_s^1 = 1+ \mathcal{P}_s$ 
et $\overline{G}_s = \overline{\mathcal{A}}_s^{\times} = \mathcal{A}_s^{\times} / \mathcal{U}_s^1$. 
Si $a$ est une arête, on note de même $\mathcal{A}_a$ l'ordre héréditaire associé, 
$\mathcal{P}_a$ le radical de Jacobson de $\mathcal{A}_a$ et $\mathcal{U}_a^1 = 1+ \mathcal{P}_a$.
On utilise les travaux \cite{SchneiderStuhler2} de P. Schneider et U. Stuhler 
en \ref{DefinitionSystCoeffShneiderEtStuhler} 
pour définir sur $X_{\mathbb{K}}$ un système de coefficients :

\begin{theoNonNum} (\cite{SchneiderStuhler2})
Pour tout sommet $s$ de $X_{\mathbb{K}}$ et toute arête $a$, 
on notera $V_s = V^{\mathcal{U}_s^1}$ l'espace des vecteurs fixes de $V$ sous l'action de 
$\mathcal{U}_s^1$ et, de même, on notera $V_a = V^{\mathcal{U}_a^1}$.\\
Soit $C_0$ l'espace des $0$-chaînes c'est-à-dire l'ensemble des fonctions 
$f : X_0 \rightarrow V$ à support fini telles que pour tout sommet 
$s$ dans $X_0$, on a $f(s) \in V^{\mathcal{U}_s^1}$. 
On note $C_1$ l'espace des $1$-chaînes c'est-à-dire l'ensemble des fonctions 
$f : X_{(1)} \rightarrow V$ à support fini telles que pour toute arête 
orientée $\langle s, t \rangle$, on a 
$f (\langle s, t \rangle) \in V^{\mathcal{U}_{\{ s, t \}}^1}$ et 
$f(\langle s, t \rangle) = - f(\langle t, s \rangle)$. 
On définit enfin l'opérateur bord :
$$
\partial_1 : C_1 \rightarrow C_0, f \mapsto \partial_1 (f) \, \, \text{où} \, \, 
\partial_1 (f) (s_0) = \sum_{s_1 \sim s_0} f(\langle s_1, s_0 \rangle)
$$
On a un complexe augmenté exact de $\mathbb{C}$-espaces vectoriels, qui est aussi une suite 
exacte de ${\rm GL}_2 (\Delta)$-modules :
$$
0 \rightarrow C_1 \xrightarrow{\partial_1} C_0  
\xrightarrow{\varepsilon} V \rightarrow 0
$$
où $\varepsilon : C_0 \rightarrow V$ est l'augmentation définie par 
$\varepsilon : w \mapsto \sum_{s \in X_0} w (s)$.
\end{theoNonNum}

En particulier, $V$ s'identifie au module d'homologie du complexe, on a :  
$V \simeq \frac{C_0}{\partial_1 (C_1)}$. 
Puis, par exactitude à droite du foncteur ${\rm Hom}$, on a la suite duale :
$$
0 \rightarrow {\rm Hom}_{\mathcal{D}^{\times}} (\pi, \mathds{1}) 
\xrightarrow{\varepsilon^{\ast}} {\rm Hom}_{\mathcal{D}^{\times}} (C_0, \mathds{1}) 
\xrightarrow{\partial_1^{\ast}} {\rm Hom}_{\mathcal{D}^{\times}} (C_1, \mathds{1})
$$
On en déduit que:
${\rm ker} (\partial_1^{\ast}) = {\rm Im} (\varepsilon^{\ast}) 
\simeq {\rm Hom}_{\mathcal{D}^{\times}} (\pi, \mathds{1})$. 
Ainsi, la représentation $\pi$ est $\mathcal{D}^{\times}$-distinguée si et seulement si l'application 
$\partial_1^{\ast}$ n'est pas injective.
Dans la suite, nous étudions donc le noyau de l'application $\partial_1^{\ast}$. 
Puisque $C_0 \simeq \bigoplus_{s \in X_0} V_s$, si 
$\varphi \in {\rm ker} (\partial_1^{\ast})$, on peut voir $\varphi$ comme une famille $(\varphi_s)_{s \in X_0}$ où 
$\varphi_s$ est la restriction de $\varphi$ à $V_s$. 
Les espaces $V_s$, comme $\overline{G}_s$-modules, 
sont entièrement connus grâce aux travaux de \cite{SilbergerZink2} et dépendent uniquement de 
la paire admissible modérée $(\chi_f, \mathbb{K}_f)$ associée à $\pi$. Nous explicitons ces espaces en 
\ref{ParamSilbergerZinkNonCuspidales} :

\begin{theoNonNum}
 La représentation 
$V_s$ de $\overline{G}_s$ se décompose simplement de la 
façon suivante :
$$
V_s \simeq ( \bigoplus_{\nu = 0}^{f-1} \overline{\chi}_0^{\Phi^{\nu}} \otimes {\rm St}_{\overline{G}_s} ) 
\oplus (\bigoplus_{0 \leq \nu_1 < \nu_2 < f} {\rm Ind}_{\overline{B}_s}^{\overline{G}_s} 
(\overline{\chi}_0^{\Phi^{\nu_1}} \otimes \overline{\chi}_0^{\Phi^{\nu_2}} ))  
$$ 
où $\Phi$ est un générateur du groupe de Galois 
${\rm Gal} (k_{\Delta, 2} / k_{\mathbb{K}})$, 
${\rm St}_{\overline{G}_s}$ est la représentation de Steinberg de $\overline{G}_s$, 
$\overline{B}_s$ est le sous-groupe de Borel standard de $\overline{G}_s$, 
$\chi = \chi_f \circ {\rm N}_{\mathbb{K}_d / \mathbb{K}_f}$ et 
$\overline{\chi}_0$ est l'unique caractère de $k_{\Delta}^{\times}$ qui vérifie :
$$
\overline{\chi} = \overline{\chi}_0 \circ {\rm N}_{k_{\Delta,2} / k_{\Delta}}
$$
\end{theoNonNum}
\noindent
(Ici, $k_{\Delta, 2}$ est une extension quadratique du corps résiduel de $\Delta$, et, 
si $r \in \mathbb{N}^{\ast}$, on note $\mathbb{K}_r$ l'unique extension non ramifiée 
de $\mathbb{K}$ de degré $r$ contenue dans $\overline{\mathbb{K}}$, une clôture algébrique de 
$\mathbb{K}$).\\
Les conditions pour que $\varphi \in {\rm Hom}_{\mathcal{D}^{\times}} (C_0, \mathds{1})$ appartienne au 
noyau de $\partial_1^{\ast}$ sont données en \ref{ConditionsDistinctionFormeLineaire} :

\begin{lemNonNum}
La forme linéaire $\varphi \in {\rm Hom}_{\mathcal{D}^{\times}} (C_0, \mathds{1})$ appartient au 
noyau de $\partial_1^{\ast}$ si et seulement si 
pour tous sommets voisins $s$ et $t$, les applications  
$\varphi_{s}$ et $\varphi_{t}$ coïncident sur $V_{\{s, t\}}$,
et pour tout $g$ dans $\mathcal{D}^{\times}$, pour tout sommet $s$, on a 
$\varphi_{s} = \varphi_{g.s} \circ \pi (g)$.
\end{lemNonNum}

\noindent
La première condition est analogue à la condition de recollement dans la définition des faisceaux. 
Pour pouvoir expliciter la deuxième condition, qui est la condition de $\mathcal{D}^{\times}$-équivariance, 
il faut en particulier connaître les 
$\mathcal{D}^{\times}$-orbites d'un sommet de $X_{\mathbb{K}}$. 
Nous déterminons ces $\mathcal{D}^{\times}$-orbites dans la partie 
\ref{PartieOrbitesSommets}. 
Enfin, si $s$ est un sommet et $\mathcal{N}_s$ le normalisateur de $\mathcal{A}_s^{\times}$ 
dans ${\rm GL}_2 (\Delta)$, alors $\mathcal{N}_s$ agit sur $V_s$ et, pour tout $g$ dans $\mathcal{N}_s$, 
$\varphi_{s} = \varphi_{s} \circ \pi (g)$. 
Pour comprendre l'action de $\mathcal{N}_s$, nous avons besoin de connaître le 
type étendu maximal de niveau $0$ 
associé à $\pi$ (au sens de Silberger et Zink) dont nous rappelons la définition dans la partie \ref{TypesEtendusMax}.
Malheureusement, \cite{SilbergerZink2} ne décrivent pas l'action de tout le normalisateur. 
Nous ne pouvons pas faire les calculs dans le cas général : on doit donc se restreindre à 
quelques cas particuliers.\\
Dans notre étude, nous nous sommes surtout intéressés au cas où l'extension $\mathbb{K} / \mathbb{F}$ est totalement 
ramifiée, modérément ramifiée. Dans ce cas, on obtient tout d'abord un résultat de multiplicité $1$ en 
\ref{TheoremeMultiplicite1NonCusp} : 

\begin{theoNonNum}
On a :
$$
{\rm dim}_{\mathbb{C}} ({\rm Hom}_{\mathcal{D}^{\times}} (\pi, \mathds{1})) \leq 1
$$
\end{theoNonNum}

\noindent 
Notons que de tels résultats de multiplicité $1$ sont connus pour ${\rm GL}_n (\mathbb{K})$, 
on pourra se référer à \cite{Flicker} pour vérifier que 
$({\rm GL}_n (\mathbb{K}), {\rm GL}_n (\mathbb{F}))$ est une paire de Gelfand. 
La preuve de ce résultat s'étend facilement à  
$({\rm GL}_{\mu} (\Delta), {\rm GL}_m (\mathcal{D}))$ lorsque 
$\mathbb{F}$ est de caractéristique nulle. 
En effet, on montre en annexe le résultat suivant :

\begin{theoNonNum}
Supposons que $\mathbb{F}$ est de caractéristique nulle. Alors, la paire  
$({\rm GL}_{\mu} (\Delta), {\rm GL}_m (\mathcal{D}))$ est une paire de Gelfand-Kazhdan. 
\end{theoNonNum}

Si $\varphi$ appartient à ${\rm Hom}_{\mathcal{D}^{\times}} (C_0, \mathds{1})$ et $s$ 
est un sommet, on note 
$\varphi_s = \sum_{\nu = 0}^{f-1} \varphi_s^{\nu} + \sum_{0 \leq \nu_1 < \nu_2 < f} \varphi_s^{\nu_1, \nu_2}$ 
la décomposition de $\varphi_s$ sur :
$$
V_s \simeq ( \bigoplus_{\nu = 0}^{f-1} \overline{\chi}_0^{\Phi^{\nu}} \otimes {\rm St}_{\overline{G}_s} ) 
\oplus (\bigoplus_{0 \leq \nu_1 < \nu_2 < f} {\rm Ind}_{\overline{B}_s}^{\overline{G}_s} 
(\overline{\chi}_0^{\Phi^{\nu_1}} \otimes \overline{\chi}_0^{\Phi^{\nu_2}} ))  
$$
Supposons que l'extension $\mathbb{K} / \mathbb{F}$ est totalement 
ramifiée. On notera $k \simeq k_{\mathbb{K}} \simeq k_{\mathbb{F}}$ le corps 
résiduel de $\mathbb{K}$. Alors, $j(X_{\mathbb{F}})$ est un sommet $s_0$ 
et on obtient les conditions suivantes de $\mathcal{D}^{\times}$-distinction en 
\ref{ConclusionCNDistinction} :

\begin{theoNonNum}
Si $f$ est pair, alors 
$\pi$ n'est pas $\mathcal{D}^{\times}$-distinguée.\\
Supposons que $\pi$ soit $\mathcal{D}^{\times}$-distinguée. Alors $f$ est impair, 
$\overline{\chi}_0$ est non trivial sur $k^{\times}$ mais trivial sur les carrés de 
$k^{\times}$, et pour tout sommet $s$ et tous entiers $0 \leq \nu_1 < \nu_2 \leq f-1$, on a :
$$
\varphi_s^{\nu_1, \nu_2} = 0
$$
et pour tout $\nu \in \{ 0, \cdots, f-1 \}$, pour tout sommet $s$ distinct de $s_0$, 
$\varphi_{s}^{\nu}$ est entièrement déterminée \mbox{par $\varphi_{s_0}^{\nu}$}. En particulier, 
si $\varphi_{s_0}^{\nu} = 0$, alors $\varphi_{s}^{\nu} = 0$.
\end{theoNonNum}

\noindent
Dans la partie \ref{DistinctionEtJLNonCusp}, on montre le résultat suivant :

\begin{theoNonNum}
Si $f$ est pair, alors $\pi$ et son image par la 
correspondance de Jacquet-Langlands ne sont 
pas distinguées, de sorte que Jacquet-Langlands préserve bien la distinction. 
\end{theoNonNum}

Enfin, dans une dernière partie, on étudie, avec la même méthode, la représentation de Steinberg 
de ${\rm GL}_2 (\Delta)$, dans le cas où l'extension $\mathbb{K}/ \mathbb{F}$ 
est non ramifiée, puis totalement ramifiée. On montre alors en \ref{NonDistinctionSteinbergEtJacquetLanglands}
le résultat suivant :

\begin{theoNonNum}
Soit $(\pi, V)$ la représentation de Steinberg de ${\rm GL}_2 (\Delta)$. 
Alors $\pi$ n'est pas $\mathcal{D}^{\times}$-distinguée et la correspondance de Jacquet-Langlands préserve 
la distinction dans ce cas.
\end{theoNonNum}


\section{Notations principales.}

Dans toute la suite, on fixe $\mathbb{K} / \mathbb{F}$ une extension quadratique 
séparable modérément ramifiée de corps locaux non archimédiens. On fixe $\mathcal{D}$ une 
$\mathbb{F}$-algèbre à division centrale d'indice $d$. 
On suppose que $d$ est pair de sorte que l'on peut considérer $\mathbb{K}$ plongé dans $\mathcal{D}$, 
et on note $\Delta$ le commutant de $\mathbb{K}$ dans 
$\mathcal{D}$. Alors $\Delta$ est une $\mathbb{K}$-algèbre à division centrale d'indice 
$\delta = d/2$ et on a un isomorphisme de $\mathbb{K}$-algèbres :
$$
\mathcal{D} \otimes_{\mathbb{F}} \mathbb{K} \rightarrow 
{\rm End}_{\Delta} (\mathcal{D}) \simeq {\rm M}_2 (\Delta) , 
x \otimes k \mapsto [f_{x \otimes k} : y \mapsto xyk]
$$
Ainsi $H = \mathcal{D}^{\times}$ est un sous-groupe de 
$G = ({\rm End}_{\Delta} (\mathcal{D}))^{\times} \simeq {\rm GL}_2 (\Delta)$. 
On note $\varpi_{\mathbb{F}}$ (resp. $\varpi_{\mathbb{K}}$, $\varpi_{\mathcal{D}}$, 
$\varpi_{\Delta}$) des uniformisantes 
de $\mathbb{F}$ (resp. $\mathbb{K}$, $\mathcal{D}$, $\Delta$), $\mathcal{O}_{\mathbb{F}}$ (resp. 
$\mathcal{O}_{\mathbb{K}}$, $\mathcal{O}_{\mathcal{D}}$, $\mathcal{O}_{\Delta}$) 
les anneaux de valuation,
$\mathcal{P}_{\mathbb{F}}$ (resp. 
$\mathcal{P}_{\mathbb{K}}$, $\mathcal{P}_{\mathcal{D}}$, $\mathcal{P}_{\Delta}$) 
les idéaux de valuation, 
$k_{\mathbb{F}}$ (resp. 
$k_{\mathbb{K}}$, $k_{\mathcal{D}}$, $k_{\Delta}$) les corps résiduels et 
$v_{\mathbb{F}}$ (resp. 
$v_{\mathbb{K}}$, $v_{\mathcal{D}}$, $v_{\Delta}$) les valuations normalisées.\\
Si l'extension $\mathbb{K} / \mathbb{F}$ est totalement ramifiée, 
on identifiera les corps résiduels de 
$\mathbb{K}$ et de $\mathbb{F}$ et on notera $k = k_{\mathbb{F}} \simeq k_{\mathbb{K}}$. 
On notera $q$ le cardinal de $k_{\mathbb{F}}$ et $Q$ celui de $k_{\Delta}$. 
On fixe $\overline{\mathbb{K}}$ une clôture algébrique de $\mathbb{K}$ 
et pour tout entier naturel $l$, 
on note $\mathbb{K}_l$ l'extension non ramifiée de degré $l$ de $\mathbb{K}$ 
contenue dans 
$\overline{\mathbb{K}}$ et $k_{\mathbb{K}, l}$ son corps résiduel. 
On fixe $\Phi$ un générateur du groupe de Galois ${\rm Gal} (k_{\mathcal{D}} / k_{\mathbb{F}})$. 
Soient ${\rm Nrd}_{\mathcal{D}}$ et ${\rm Nrd}_{\Delta}$ les normes réduites.\\
On notera $X_{\mathbb{F}}$ (resp. $X_{\mathbb{K}}$) 
l'immeuble de Bruhat-Tits de $\mathcal{D}^{\times}$ 
(resp. $G$) et $j : X_{\mathbb{F}} \hookrightarrow X_{\mathbb{K}}$ 
l'injection naturelle entre ces immeubles 
dont nous rappelons la définition en \ref{RappelsInjectionsImmeubleTits}. 
On notera $A_{\mathbb{K}}$ l'appartement standard de $X_{\mathbb{K}}$, $X_0$ l'ensemble des sommets de 
$X_{\mathbb{K}}$. Si $s$ et $t$ sont deux sommets dans $X_0$, $d(s, t)$ désigne la 
distance de $s$ à $t$. Si $s \in X_0$, on note $\mathcal{A}_s$ l'ordre héréditaire associé à $s$, 
$\mathcal{U}_s = \mathcal{A}_s^{\times}$, $\mathcal{P}_s$ le radical de Jacobson de~$\mathcal{A}_s$, 
$\mathcal{U}_s^1 = 1+ \mathcal{P}_s$. On note également 
$\overline{G}_s = \overline{\mathcal{A}}_s^{\times} = \mathcal{A}_s^{\times} / \mathcal{U}_s^1$, 
${\rm St}_{\overline{G}_s}$ la représentations de Steinberg de $\overline{G}_s$ et 
$\overline{B}_s$ le sous-groupe de Borel standard de $\overline{G}_s$. 
On note $\mathcal{R}_0^2 (G)$ l'ensemble des classes d'isomorphisme des membres de la série discrète de $G$ 
de niveau $0$. Le caractère trivial d'un groupe est noté $\mathds{1}$. 
Si $\mathbb{E} / \mathbb{M}$ est une extension finie de corps, 
${\rm N}_{\mathbb{E} / \mathbb{M}} : \mathbb{E}^{\times} \rightarrow \mathbb{M}^{\times}$ 
désigne l'application norme.\\
Toutes les représentations considérées seront supposées complexes et lisses.\\
On identifiera les sommets de $X_{\mathbb{K}}$ aux $\mathcal{O}_{\Delta}$-chaînes de 
réseaux de période $1$ d'un $\Delta$-espace vectoriel (à droite) de dimension $2$ (cf. 
\cite{BushnellFröhlich} pages 212 à 224). Si $L$ est un $\mathcal{O}_{\Delta}$-réseau d'un 
$\Delta$-espace vectoriel de dimension $2$, on notera 
$[L] = \{ L \varpi_{\Delta}^k : k \in \mathbb{Z} \}$ le sommet de $X_{\mathbb{K}}$ associé.


\section{Orbites des sommets de $X_{\mathbb{K}}$ sous l'action de $\mathcal{D}^{\times}$.}\label{PartieOrbitesSommets}

\begin{nota}
On fixe $(e_1, e_2)$ une $\Delta$-base de $\mathcal{D}$.  
Soit $\mathcal{S}_0$ l'ensemble des sommets de l'appartement standard $A_{\mathbb{K}}$ :
$$
\mathcal{S}_0 
= \{ s_k = [ e_1 \mathcal{O}_{\Delta} \oplus e_2 \mathcal{P}_{\Delta}^k] : k \in \mathbb{Z} \}
$$
\end{nota}

\subsection{Rappels sur les injections d'immeubles.}\label{RappelsInjectionsImmeubleTits}

On pourra se référer à \cite{Tits} (2.6 page 47) pour les résultats suivants:

\begin{theo} (\cite{Tits})
Il existe une injection naturelle $j : X_{\mathbb{F}} \rightarrow X_{\mathbb{K}}$
vérifiant les trois propriétés suivantes:
\begin{itemize}
\item[a)] L'application $j$ est $\mathcal{D}^{\times}$-équivariante, c'est-à-dire que 
pour tout $g$ dans $\mathcal{D}^{\times}$ et tout $x$ dans $X_{\mathbb{F}}$, on a 
$j(g.x) = g. j(x)$.
\item[b)] L'image de $j$ est incluse dans $ X_{\mathbb{K}}^{{\rm Gal} (\mathbb{K} / \mathbb{F})}$, 
où $ X_{\mathbb{K}}^{{\rm Gal} (\mathbb{K} / \mathbb{F})}$ désigne les 
éléments de $X_{\mathbb{K}}$ qui sont fixes sous l'action du groupe de Galois 
${\rm Gal} (\mathbb{K} / \mathbb{F})$.
\item[c)] L'application $j$ est affine, c'est-à-dire que pour tout appartement $\mathcal{A}$ de 
$X_{\mathbb{F}}$ (vu comme espace affine), il existe un appartement $\mathcal{B}$ de $X_{\mathbb{K}}$ tel que 
$j(\mathcal{A}) \subseteq \mathcal{B}$ et 
$j_{\vert \mathcal{A}} : \mathcal{A} \rightarrow \mathcal{B}$ est une application affine.
\end{itemize}
\end{theo}

Nous allons utiliser le théorème suivant démontré dans \cite{PrasadYu} (Théorème 1.9 page 555):

\begin{theo}\label{TheoremeYuPrasad} (\cite{PrasadYu})
Soit $G$ un groupe réductif défini sur un corps local non archimédien $\mathbb{L}$ de 
carcatéristique résiduelle $p$. Soit $F \subseteq {\rm Aut}_{\mathbb{L}} (G)$ un 
groupe fini d'ordre non divisible par $p$. On note $G^{F}$ les points fixes de $G$ sous l'action de $F$ 
et $H = (G^{F})^{\circ}$ la composante de l'unité de $G^{F}$. 
Soient $X_G$ et $X_H$ les immeubles de Bruhat-Tits de $G$ et de $H$ respectivement. 
Alors, via l'injection naturelle des immeubles $j : X_H \hookrightarrow X_G$, on peut identifier 
$X_H$ à $X_G^{F}$, l'ensemble des points fixes de $X_G$ sous l'action de $F$.

\end{theo}

On en déduit le résultat suivant:

\begin{propo}
L'application $j$ est unique et 
${\rm Im} (j) = X_{\mathbb{K}}^{{\rm Gal} (\mathbb{K} / \mathbb{F})}$.
\end{propo}

On explicite alors $j$ dans le cas où l'extension $\mathbb{K} / \mathbb{F}$ est non ramifiée 
(on pourra se référer à \cite{Coniglio}, propriété 3.2.12) :

\begin{propo}
Supposons que l'extension $\mathbb{K} / \mathbb{F}$ est non ramifiée. 
L'immeuble $X_{\mathbb{F}}$ de $\mathcal{D}^{\times}$ possède un seul sommet 
$[\mathcal{O}_{\mathcal{D}}]$ et on a 
$j([\mathcal{O}_{\mathcal{D}}]) = m_0$ 
où $m_0$ est le milieu de l' arête $[s_0, s_1 ]$ dans l'immeuble $X_{\mathbb{K}}$.
\end{propo}

On explicite également $j$ dans le cas où $\mathbb{K} / \mathbb{F}$ est 
totalement ramifiée, modérément ramifiée (cf. \cite{Coniglio}, propriété 3.2.20) :

\begin{propo}
On suppose que l'extension $\mathbb{K} / \mathbb{F}$ est 
totalement ramifiée, modérément ramifiée. 
On a 
$j ([\mathcal{O}_{\mathcal{D}}]) = s_0$.
\end{propo}

\subsection{Cas où l'extension $\mathbb{K} / \mathbb{F}$ est non ramifiée.}

On suppose dans cette partie que l'extension $\mathbb{K} / \mathbb{F}$ est non ramifiée. 
On détermine les $\mathcal{D}^{\times}$-orbites des sommets de $X_{\mathbb{K}}$ :

\begin{propo}\label{OrbitesCasNonRam}
Les $\mathcal{D}^{\times}$-orbites des sommets de $X_{\mathbb{K}}$ 
sont exactement les sphères de centre $m_0$ et de rayon $k+1/2$ pour $k \in \mathbb{N}$.
\end{propo}

\begin{démo}
La base $(e_1, e_2)$ étant fixée, on identifie $G$ à ${\rm GL}_2 (\Delta)$.\\ 
On remarque que 
\mbox{$v_{\Delta} (\varpi_{\mathcal{D}}^2) = v_{\Delta} (\varpi_{\Delta}) = 1 
= v_{\Delta} (t_0^2)$} où :
$$
t_0 
= \left( \begin{array}{cc}
   0 & 1\\
   \varpi_{\Delta} & 0
  \end{array} \right)
$$
Il existe donc $u$ dans ${\rm GL}_2 (\mathcal{O}_{\Delta})$ (donc $u$ fixe $s_0$) tel que 
$\varpi_{\mathcal{D}} = t_0 u$. Comme $t_0$ échange $s_0$ et $s_1$, l'uniformisante  
$\varpi_{\mathcal{D}}$ fixe $m_0$ et  
échange les sommets $s_0$ et $s_1$. De plus,  
$\mathcal{O}_{\mathcal{D}}^{\times} \subseteq {\rm GL}_2 (\mathcal{O}_{\Delta})$ fixe $s_0$ et $m_0$, 
donc fixe point par point l'arête 
$[s_0, s_1]$. 
Enfin, puisque $v_{\Delta} (\varpi_{\mathcal{D}}^2) = v_{\Delta} (\varpi_{\Delta})$, 
$\varpi_{\mathcal{D}}^2$ fixe $s_0$, et donc fixe point par point l'arête 
$[s_0, s_1]$.\\
Puisque $\mathcal{D}^{\times}$ agit via des isométries sur l'immeuble 
$X_{\mathbb{K}}$, il est clair que la $\mathcal{D}^{\times}$-orbite d'un sommet $s$ de $X_{\mathbb{K}}$ 
est contenue dans la sphère de centre $m_0$ et de rayon 
$d (s, m_0)$.\\
Soit $s$ un sommet de $X_{\mathbb{K}}$. On distingue deux cas :
 \begin{itemize}
\item[$\ast$] Si $d(s_0,s) = d(s_1, s) +1$. 
Soit $s_k$ dans l'appartement $A_{\mathbb{K}}$ tel que $s_1 \in [s_0, s_k]$ et :
$$
d (m_0, s_k) = d (s, m_0)
$$
Alors $k>0$, $d(m_0, s) = k-1/2$ et $s_k = [ e_1 \mathcal{O}_{\Delta} \oplus e_2 \mathcal{P}_{\Delta}^k]$. 
Les sommets $s_0$ et $s$ sont dans un même appartement $A$. Il existe donc $ (f_1, f_2)$, une 
$\Delta$-base de $\mathcal{D}$ telle que :
$$
s_0 = [ f_1 \mathcal{O}_{\Delta} \oplus f_2 \mathcal{O}_{\Delta}], \, 
s = [ f_1 \mathcal{O}_{\Delta} \oplus f_2 \mathcal{P}_{\Delta}^k]
$$
Ainsi, 
$s_0 = [ e_1 \mathcal{O}_{\Delta} \oplus e_2 \mathcal{O}_{\Delta}] 
= [ f_1 \mathcal{O}_{\Delta} \oplus f_2 \mathcal{O}_{\Delta}]$. 
Quitte à remplacer $f_i$ par $f_i \varpi_{\Delta}^l$, on peut supposer~que :
$$
e_1 \mathcal{O}_{\Delta} \oplus e_2 \mathcal{O}_{\Delta} = 
f_1 \mathcal{O}_{\Delta} \oplus f_2 \mathcal{O}_{\Delta}
$$
Soit $d$ dans $\mathcal{D}^{\times}$ tel que $d.e_1 = f_1$. On sait que 
$d$ s'écrit sous la forme 
$d = \varpi_{\mathcal{D}}^m x$ où $x \in \mathcal{O}_{\mathcal{D}}^{\times}$. 
Supposons tout d'abord que $m$ est pair. 
Alors, puisque $v_{\Delta} (\varpi_{\mathcal{D}}^2) = 1$, on a :
$$
d = \varpi_{\mathcal{D}}^m x = \varpi_{\Delta}^{m/2} x^{'}, 
x^{'} \in {\rm GL}_2 (\mathcal{O}_{\Delta})
$$
On en déduit que $d.s_0 = s_0$ donc :
$$
d. (e_1 \mathcal{O}_{\Delta} + e_2 \mathcal{O}_{\Delta}) 
= f_1 \mathcal{O}_{\Delta} + d.e_2 \mathcal{O}_{\Delta} 
\in [f_1 \mathcal{O}_{\Delta} + f_2 \mathcal{O}_{\Delta}]
$$
et donc $d. (e_1 \mathcal{O}_{\Delta} + e_2 \mathcal{O}_{\Delta}) 
= f_1 \mathcal{O}_{\Delta} + f_2 \mathcal{O}_{\Delta}$. De plus, puisque 
$\mathcal{P}_{\Delta}^k \subseteq \mathcal{O}_{\Delta}$ :
\begin{eqnarray*}
 d.s_k 
& = & [ d.(e_1 \mathcal{O}_{\Delta} + e_1 \mathcal{P}_{\Delta}^k + e_2 \mathcal{P}_{\Delta}^k)] 
 =  [ d.(e_1 \mathcal{O}_{\Delta}) + d. (e_1 \mathcal{O}_{\Delta} + e_2 \mathcal{O}_{\Delta}) \varpi_{\Delta}^k]\\
& = & [ f_1 \mathcal{O}_{\Delta} + (f_1 \mathcal{O}_{\Delta} + f_2 \mathcal{O}_{\Delta}) \varpi_{\Delta}^k]
 =  [ f_1 \mathcal{O}_{\Delta} + f_2 \mathcal{P}_{\Delta}^k ] = s
\end{eqnarray*}
Montrons qu'il est impossible que $m$ soit impair. 
Raisonnons par l'absurde. Supposons que $m$ est impair. 
Alors $d$ échange les sommets $s_0$ et $s_1$, ainsi $d.s_1 = s_0$ et :
$$
d.(e_1 \mathcal{O}_{\Delta} + e_2 \mathcal{P}_{\Delta}) 
= f_1 \mathcal{O}_{\Delta} + d.e_2 \mathcal{P}_{\Delta} 
\in s_0 = [f_1 \mathcal{O}_{\Delta} + f_2 \mathcal{O}_{\Delta}]
$$
On en déduit que  
$d.(e_1 \mathcal{O}_{\Delta} + e_2 \mathcal{P}_{\Delta}) 
= f_1 \mathcal{O}_{\Delta} + f_2 \mathcal{O}_{\Delta}$. 
Puisque $\mathcal{P}_{\Delta}^{k-1} \subseteq \mathcal{O}_{\Delta}$, on a :
$$
d.s_k 
= [d.(e_1 \mathcal{O}_{\Delta} + e_1 \mathcal{P}_{\Delta}^{k-1} 
     + e_2 \mathcal{P}_{\Delta}^{k})] 
= [f_1 \mathcal{O}_{\Delta} + 
     (f_1 \mathcal{O}_{\Delta} + f_2 \mathcal{O}_{\Delta}) \varpi_{\Delta}^{k-1}] 
=  [f_1 \mathcal{O}_{\Delta} + f_2 \mathcal{P}_{\Delta}^{k-1}]
$$
Notons $\widetilde{s} = [f_1 \mathcal{O}_{\Delta} + f_2 \mathcal{P}_{\Delta}^{k-1}]$. 
Alors :
$$
d (\widetilde{s}, s_0) = d (d.s_k, s_0) = d (s_k, d^{-1}.s_0) 
= d (s_k, s_1) = k-1
$$
et $d( \widetilde{s}, m_0) = d (s_k, m_0) = k-1/2$. 
De plus, par hypothèse, $d(s_1, s) = k-1$. On a donc deux chemins géodésiques reliant $s_0$ à $s$ :
$$
[s_0, s_1, \cdots, s] \, \, \text{et} \, \, 
[s_0, \cdots, \widetilde{s}, s]
$$
Ces deux chemins sont égaux. Par suite, $d(s_1, \widetilde{s}) = k-2$ et 
$d (\widetilde{s}, m_0) = k-2+ \frac{1}{2} = k - 3/2$. 
D'où une contradiction. 
\item[$\ast$] Si $d (s_0, s) = d(s_1, s) -1$. 
Comme précédemment, on choisit $s_k$ dans l'appartement $A_{\mathbb{K}}$ tel que  
$s_0 \in [s_k, s_1]$ et 
$d(m_0, s) = d (m_0, s_k)$. Ainsi $k \leq 0$, $d(m_0, s) = -k+ 1/2$ et :
$$
s_k = [ e_1 \mathcal{O}_{\Delta} + e_2 \mathcal{P}_{\Delta}^k ]
$$
Il existe $A$ un appartement contenant $s$ et $s_0$. Soit $(f_1, f_2)$ une 
$\Delta$-base de $\mathcal{D}$ telle que :
$$
s_0 = [ f_1 \mathcal{O}_{\Delta} + f_2 \mathcal{O}_{\Delta} ] , \, 
s = [ f_1 \mathcal{O}_{\Delta} + f_2 \mathcal{P}_{\Delta}^k ]
$$
On suppose à nouveau que 
$f_1 \mathcal{O}_{\Delta} + f_2 \mathcal{O}_{\Delta} = e_1 \mathcal{O}_{\Delta} + e_2 \mathcal{O}_{\Delta}$. 
Soit $d$ dans $\mathcal{D}^{\times}$ tel que $d.e_2 = f_2$. On écrit $d$ sous la forme 
$d = x \varpi_{\mathcal{D}}^m$ avec $x \in \mathcal{O}_{\mathcal{D}}^{\times}$. \\
Supposons que $m$ est pair.
Alors, comme précédemment, on vérifie que :
$$
d.(e_1 \mathcal{O}_{\Delta} + e_2 \mathcal{O}_{\Delta}) 
= f_1 \mathcal{O}_{\Delta} + f_2 \mathcal{O}_{\Delta}
$$
Ainsi, en utilisant le fait que $\mathcal{O}_{\Delta} \subseteq \mathcal{P}_{\Delta}^k$, on montre que 
$d.s_k = s$.\\
Si $m$ est impair, 
alors $d.s_1 = s_0$ d'où :
$$
d.(e_1 \mathcal{O}_{\Delta} + e_2 \mathcal{P}_{\Delta}) 
= d.e_1 \mathcal{O}_{\Delta} + f_2 \mathcal{P}_{\Delta} 
\in s_0 = [f_1 \mathcal{O}_{\Delta} + f_2 \mathcal{O}_{\Delta}]
$$
On en déduit que 
$d.(e_1 \mathcal{O}_{\Delta} + e_2 \mathcal{P}_{\Delta}) 
= f_1 \mathcal{O}_{\Delta} + f_2 \mathcal{O}_{\Delta}$. 
Puisque $k \leq 0$, on a $\mathcal{P}_{\Delta} \subseteq \mathcal{O}_{\Delta} 
\subseteq \mathcal{P}_{\Delta}^k$, d'où :
$$
 d.s_k 
= [ d.(e_1 \mathcal{O}_{\Delta} + e_2 \mathcal{P}_{\Delta} + e_2 \mathcal{P}_{\Delta}^k)] 
= [ f_1 \mathcal{O}_{\Delta} + f_2 \mathcal{O}_{\Delta} + f_2 \mathcal{P}_{\Delta}^k] = s
$$

\end{itemize}

\end{démo}

\subsection{Cas où l'extension $\mathbb{K} / \mathbb{F}$ est totalement ramifiée, 
modérément ramifiée.}

On suppose ici que l'extension $\mathbb{K} / \mathbb{F}$ est totalement ramifiée, 
modérément ramifiée.\\
Avec un raisonnement analogue à celui de la démonstration de la proposition 
\ref{OrbitesCasNonRam}, on détermine les $\mathcal{D}^{\times}$-orbites des sommets de $X_{\mathbb{K}}$ :

\begin{propo}
Les $\mathcal{D}^{\times}$-orbites des sommets de $X_{\mathbb{K}}$ 
sont exactement les sphères de centre $s_0$ et de rayon $k \in \mathbb{N}$.
\end{propo}

Nous montrons par la suite que tous les sommets d'une même $\mathcal{D}^{\times}$-orbite 
sont en fait contenus dans la même $\mathcal{O}_{\mathcal{D}}^{\times}$-orbite :

\begin{propo}\label{OrbitesCasTotRamifie}
Les $\mathcal{O}_{\mathcal{D}}^{\times}$-orbites des sommets de $X_{\mathbb{K}}$ sont 
aussi les sphères de centre $s_0$ et de rayon $k$, avec $k$ dans $\mathbb{N}$.
\end{propo}

\begin{démo}
On a 
$\mathcal{O}_{\mathcal{D}} = e_1 \mathcal{O}_{\Delta} 
+ e_2 \mathcal{O}_{\Delta}$. 
Pour $1 \leq i \leq 2$, comme $e_i \in \mathcal{O}_{\mathcal{D}}$, on a 
$v_{\mathcal{D}} (e_i) \geq 0$. 
Si on suppose que $v_{\mathcal{D}} (e_1) > 0$ et $v_{\mathcal{D}} (e_2) >0$, alors, 
pour tout $(\lambda_1, \lambda_2) \in \mathcal{O}_{\Delta} \times \mathcal{O}_{\Delta} \setminus \{ (0, 0) \}$, 
on a $e_1 \lambda_1 + e_2 \lambda_2 \notin \mathcal{O}_{\mathcal{D}}^{\times}$. 
On peut donc supposer que $v_{\mathcal{D}} (e_1) = 0$.  
Puisque l'extension $\mathbb{K} / \mathbb{F}$ est totalement ramifiée, modérément ramifiée, 
on peut supposer que $\varpi_{\mathbb{K}}^2 = \varpi_{\mathbb{F}}$ et donc 
$\varpi_{\mathcal{D}}^{d} = \varpi_{\mathbb{F}} = \varpi_{\mathbb{K}}^2 = \varpi_{\Delta}^{2\delta}$ 
(avec $2 \delta = d$). 
On en déduit que $\varpi_{\Delta} \in \varpi_{\mathcal{D}} \mathcal{O}_{\mathcal{D}}^{\times}$. 
Par suite, $\mathcal{O}_{\mathcal{D}} / \mathcal{O}_{\mathcal{D}} \varpi_{\Delta}$ est isomorphe 
au corps résiduel $k_{\mathcal{D}}$ de $\mathcal{D}$ et est donc un $k_{\Delta}$-espace vectoriel 
de dimension $2$. Or, si $v_{\mathcal{D}} (e_2) >0$, on a :
$$
\mathcal{O}_{\mathcal{D}} / \mathcal{O}_{\mathcal{D}} \varpi_{\Delta} 
= (e_1 \mathcal{O}_{\Delta} + e_2 \mathcal{O}_{\Delta}) / 
(\mathcal{O}_{\mathcal{D}} \varpi_{\Delta}) = \overline{e}_1 k_{\Delta}
$$
est un $k_{\Delta}$-espace vectoriel de dimension $1$. D'où une contradiction. 
On en déduit que $v_{\mathcal{D}} (e_1) = v_{\mathcal{D}} (e_2) = 0$. 
Ensuite, on reprend la démonstration précédente : on fixe $s$ un sommet de 
$X_{\mathbb{K}}$ et l'on note $k$ la distance de $s$ à $s_0$. On fixe 
$(f_1, f_2)$ une $\Delta$-base de $\mathcal{D}$ telle que :
$$
s = [f_1 \mathcal{O}_{\Delta} + f_2 \mathcal{P}_{\Delta}^k ]
$$
et on peut supposer que 
$\mathcal{O}_{\mathcal{D}} = 
e_1 \mathcal{O}_{\Delta} + e_2 \mathcal{O}_{\Delta} 
= f_1 \mathcal{O}_{\Delta} + f_2 \mathcal{O}_{\Delta}$. 
Alors, il existe $d$ dans $\mathcal{D}^{\times}$ tel que $d.s = s_k$. 
D'après la démonstration précédente, un tel élément $d$ s'écrit comme 
produit d'éléments $d_r$ tels qu'il existe $l_1$ dans 
$\{ e_1, e_2 \}$ et $l_2$ dans $\{ e_1, e_2, f_1, f_2 \}$ 
tel que $d_r.l_1 = l_2$. Comme 
$e_1, e_2, f_1, f_2 \in \mathcal{O}_{\mathcal{D}}^{\times}$, on a :
$$
d_r.l_1 = l_2 \Rightarrow 
v_{\mathcal{D}} (l_2) = 0 = v_{\mathcal{D}} (d_r.l_1) 
= v_{\mathcal{D}} (d_r) + v_{\mathcal{D}} (l_1) = 
v_{\mathcal{D}} (d_r)
$$
On en déduit que $d_r \in \mathcal{O}_{\mathcal{D}}^{\times}$ et donc 
$d \in \mathcal{O}_{\mathcal{D}}^{\times}$.

\end{démo}

\section{Systèmes de coefficients de Schneider et Stuhler.}

Nous allons rappeler ici les résultats de \cite{SchneiderStuhler1} et \cite{SchneiderStuhler2}.

\begin{nota}
On fixe $(\pi, V)$ une représentation lisse irréductible de $G$, membre de la série discrète, 
de niveau $0$, non cuspidale. 
Pour $l$ dans $\{ 0, 1 \}$, on note $X_l$ l'ensemble des $l$-simplexes de $X_{\mathbb{K}}$ 
et $X_{(l)}$ l'ensemble des $l$-simplexes orientés de $X_{\mathbb{K}}$.\\
Si $s$ et $t$ sont deux sommets de $X_{\mathbb{K}}$ reliés par une arête 
(on écrira alors $s \sim t$), on notera $a = \{ s, t \}$ 
l'arête non orientée reliant $s$ et $t$, et $\langle s, t \rangle$ l'arête orientée. On notera $\mathcal{A}_a$ 
l'ordre héréditaire associé à $a$, $\mathcal{U}_a = \mathcal{A}_a^{\times}$, $\mathcal{P}_a$ le 
radical de Jacobson de $\mathcal{A}_a$ et $\mathcal{U}_a^1 = 1+ \mathcal{P}_a$.
\end{nota}

\subsection{Définition d'un système de coefficients.}

\begin{déf}
(\cite{SchneiderStuhler2}).\\
Un système de coefficients d'espaces vectoriels $W$ sur l'arbre $X_{\mathbb{K}}$ est la donnée :
\begin{itemize}
 \item d'une famille d'espaces vectoriels $W_{\sigma}$ où $\sigma$ parcourt l'ensemble des simplexes de 
     $X_{\mathbb{K}}$.
 \item d'applications linéaires $r_{\sigma}^{\tau} : W_{\tau} \rightarrow W_{\sigma}$ pour chaque couple de 
    simplexe $\sigma \subseteq \tau$ telles que $r_{\sigma}^{\sigma} = Id_{W_{\sigma}}$ et 
    si $\sigma_1 \subseteq \sigma_2 \subseteq \sigma_3$, alors 
    $r_{\sigma_1}^{\sigma_3} = r_{\sigma_1}^{\sigma_2} \circ r_{\sigma_2}^{\sigma_3}$.
\end{itemize}
Soit $C_0$ l'espace des $0$-chaînes c'est-à-dire l'ensemble des applications  
$f : X_0 \rightarrow W$ à support fini telles que pour tout sommet 
$s$ dans $X_0$, on a $f(s) \in W_s$. 
On note $C_1$ l'espace des $1$-chaînes c'est-à-dire l'ensemble des applications  
$f : X_{(1)} \rightarrow W$ à support fini telles que pour toute arête 
orientée $\langle s, t \rangle$, on a 
$f (\langle s, t \rangle) \in W_{\{ s, t \}}$ et 
$f(\langle s, t \rangle) = - f(\langle t, s \rangle)$. 
On définit enfin l'opérateur bord :
$$
\partial_1 : C_1 \rightarrow C_0, f \mapsto \partial_1 (f) \, \, \text{où} \, \, 
\partial_1 (f) (s_0) = \sum_{s_1 \sim s_0} f(\langle s_1, s_0 \rangle)
$$
On a alors un complexe de chaînes :
$$
C_1 \xrightarrow{\partial_1} C_0
$$
\end{déf}

Dans leur article \cite{SchneiderStuhler2}, 
P.Schneider et U.Stuhler définissent un système de coefficients particulier sur~$X_{\mathbb{K}}$ 
d'espace vectoriel $V$ (l'espace de la représentation $\pi$) de la façon suivante :

\begin{nota}\label{DefinitionSystCoeffShneiderEtStuhler}
\begin{itemize}
\item Pour tout sommet $s$ de $X_{\mathbb{K}}$ et toute arête $a \in X_{(1)}$, 
on notera $V_s = V^{\mathcal{U}_s^1}$ l'espace vectoriel des vecteurs fixes de $V$ 
sous l'action de $\mathcal{U}_s^1$, et de même $V_a = V^{\mathcal{U}_a^1}$. 
\item Pour chaque couple de 
    simplexe $\sigma \subseteq \tau$, on note $r_{\sigma}^{\tau} : V_{\tau} \rightarrow V_{\sigma}$ 
l'inclusion $V_{\tau} \subseteq V_{\sigma}$. 
\end{itemize}
On note alors $C_0$ l'espace des $0$-chaînes (c'est-à-dire l'ensemble des applications  
$f : X_0 \rightarrow V$ à support fini telles que pour tout sommet 
$s$ dans $X_0$, on a $f(s) \in V^{\mathcal{U}_s^1}$)  
et $C_1$ l'espace des $1$-chaînes (c'est-à-dire l'ensemble des applications  
$f : X_{(1)} \rightarrow V$ à support fini telles que pour toute arête 
orientée $\langle s, t \rangle$, on a 
$f (\langle s, t \rangle) \in V^{\mathcal{U}_{\{ s, t \}}^1}$ et 
$f(\langle s, t \rangle) = - f(\langle t, s \rangle)$). 
\end{nota}

D'après \cite{SchneiderStuhler2} (Théorème II.3.1), on a alors une résolution de $V$ :

\begin{propo}\label{ResolutionPi}
On a un complexe augmenté exact de $\mathbb{C}$-espaces vectoriels :
$$
0 \rightarrow C_1 \xrightarrow{\partial_1} C_0  
\xrightarrow{\varepsilon} V \rightarrow 0
$$
où $\varepsilon : C_0 \rightarrow V$ est l'augmentation définie par :
$$
\varepsilon : w \mapsto \sum_{s \in X_0} w (s)
$$
\end{propo}

Dans la suite, il sera plus commode de voir l'espace des $0$-chaînes (resp. des $1$-chaînes) 
non plus comme un espace de fonctions mais comme une somme directe, indexée par les sommets $s$  
(resp. les arêtes $a$), des espaces $V_s$ (resp. $V_a$) définis précédemment :

\begin{lem}
On a deux isomorphismes de $\mathbb{C}$-espaces vectoriels :
$$
C_0 \rightarrow \bigoplus_{s \in X_0} V^{\mathcal{U}_s^1}, \, 
f \mapsto (f(s))_{s \in X_0}
$$
et :
$$
C_1 \rightarrow \bigoplus_{\langle s, t \rangle \in X_{(1)}} V^{\mathcal{U}_{\{ s, t \}}^1}, \, 
f \mapsto (f(\langle s, t \rangle))_{\langle s, t \rangle \in X_{(1)}}
$$
\end{lem}

\begin{rmq}
Via les identifications précédentes, l'opérateur bord $\partial_1$ est définie de la façon suivante :
$$
\partial_1 : C_1 = \bigoplus_{\langle s, t \rangle \in X_{(1)}} V^{\mathcal{U}_{\{ s, t \}}^1} 
\rightarrow C_0 = \bigoplus_{ s \in X_{(0)}} V^{\mathcal{U}_{s}^1} , \, 
(w_{\langle s, t \rangle})_{\langle s, t \rangle \in X_{(1)}} 
\mapsto (v_s)_{s \in X_{(0)}}
$$
telle que 
$v_s = \sum_{t \sim s} w_{\langle t, s \rangle}$.
\end{rmq}

Le complexe augmenté est en fait une suite exacte de $G$-modules pour l'action de $G$ 
sur les espaces $C_0$ et $C_1$, définie dans \cite{SchneiderStuhler2}, que nous rappelons ici :

\begin{propo}
(\cite{SchneiderStuhler2}).\\
Le groupe $G$ agit sur $C_0$ de la façon suivante : 
$$
g.(v_s)_{s \in X_{(0)}} = (v^{'}_s)_{s \in X_{(0)}}, \, \, \text{pour} \, \, 
g \in G \, \, \text{et} \, \, (v_s)_{s \in X_{(0)}} \in C_0
$$
où $v^{'}_s = \pi (g). v_{g^{-1}.s}$. 
De même, $G$ agit sur $C_1$ : 
$$
g. (w_{\langle s, t \rangle})_{\langle s, t \rangle \in X_{(1)}} 
= (w^{'}_{\langle s, t \rangle})_{\langle s, t \rangle \in X_{(1)}}, \, \, \text{pour} \, \, 
g \in G \, \, \text{et} \, \, 
(w_{\langle s, t \rangle})_{\langle s, t \rangle \in X_{(1)}} \in C_1
$$
où $w^{'}_{\langle s, t \rangle} = \pi (g). w_{\langle g^{-1}.s, g^{-1}.t \rangle}$. 
On a alors une suite exacte de $G$-modules :
$$
0 \rightarrow C_1 \xrightarrow{\partial_1} C_0  
\xrightarrow{\varepsilon} V \rightarrow 0
$$
\end{propo}

\begin{rmq}\label{RepresentationQuotientInduiteCompacte}
On en déduit que $\pi$ est le quotient d'une induite compacte.\\
En effet, si l'on fixe $s_0$ un sommet de $X_{\mathbb{K}}$, alors $V_{s_0}$ est 
un $\mathcal{R}_1$-module (où $\mathcal{R}_1$ est le normalisateur dans $G$ de 
$\mathcal{A}_{s_0}^{\times}$). 
Soit $s$ un sommet de $X_{\mathbb{K}}$, il existe donc $g$ dans $G$ tel que 
$g.s_0 = s$, alors $g \mathcal{U}_{s_0}^1 g^{-1} = \mathcal{U}_s^1$ et 
$V_s = \pi (g). V_{s_0}$. On en déduit que :
$$
C_0 = \sum_{g \in G}  \pi (g).V_{s_0} 
= \sum_{g \in G / \mathcal{R}_1}  \pi (g).V_{s_0} 
= {\rm c-ind}_{\mathcal{R}_1}^G (V_{s_0})
$$
Ainsi, d'après \ref{ResolutionPi}, on a :
$$
V \simeq \frac{{\rm c-ind}_{\mathcal{R}_1}^G (V_{s_0})}{\partial_1 (C_1)}
$$
\end{rmq}

\subsection{Lien avec la $\mathcal{D}^{\times}$-distinction de $\pi$.}

\begin{propo}
Par exactitude à droite du foncteur ${\rm Hom}$, on a la suite duale :
$$
0 \rightarrow {\rm Hom}_{\mathcal{D}^{\times}} (\pi, \mathds{1}) 
\xrightarrow{\varepsilon^{\ast}} {\rm Hom}_{\mathcal{D}^{\times}} (C_0, \mathds{1}) 
\xrightarrow{\partial_1^{\ast}} {\rm Hom}_{\mathcal{D}^{\times}} (C_1, \mathds{1})
$$
On en déduit que:
$$
{\rm ker} (\partial_1^{\ast}) = {\rm Im} (\varepsilon^{\ast}) 
\simeq {\rm Hom}_{\mathcal{D}^{\times}} (\pi, \mathds{1})
$$
Ainsi, la représentation $\pi$ est $\mathcal{D}^{\times}$-distinguée si et seulement si l'application 
$\partial_1^{\ast}$ n'est pas injective où 
$\partial_1^{\ast} : {\rm Hom}_{\mathcal{D}^{\times}} (C_0, \mathds{1}) 
\rightarrow {\rm Hom}_{\mathcal{D}^{\times}} (C_1, \mathds{1})$ telle 
que pour toute forme linéaire $\mathcal{D}^{\times}$-invariante $f : C_0 \rightarrow \mathbb{C}$, on~a :
$$
\partial_1^{\ast} (f) : C_1 \rightarrow \mathbb{C}, 
w \mapsto f (\partial_1 (w))
$$
\end{propo}

\begin{rmq}
On a un isomorphisme de $\mathbb{C}$-espaces vectoriels:
$$
{\rm Hom}_{\mathbb{C}} (C_0, \mathbb{C}) \rightarrow 
\prod_{s \in X_0} {\rm Hom}_{\mathbb{C}} (V^{\mathcal{U}_s^1}, \mathbb{C}) , \, 
\varphi \mapsto (\varphi_s = \varphi_{\vert V^{\mathcal{U}_s^1}})_{s \in X_0}
$$
\end{rmq}

Ainsi, si $\varphi$ appartient à ${\rm Hom}_{\mathcal{D}^{\times}} (C_0, \mathds{1})$, on peut 
voir $\varphi$ comme une collection de formes linéaires $(\varphi_s)_{s \in X_0}$ où, 
pour tout sommet $s$,  
$\varphi_s$ est la restriction de $\varphi$ à l'espace $V_s$. 
On a les conditions suivantes pour que $\varphi$ appartienne à ${\rm ker} (\partial_1^{\ast})$ :

\begin{propo}\label{ConditionsDistinctionFormeLineaire}
Soit $\varphi = (\varphi_s)_{s \in X_0}$ une forme linéaire sur $C_0$.\\ 
Alors $\varphi \in {\rm Hom}_{\mathcal{D}^{\times}} (C_0, \mathds{1}) \cap {\rm ker} (\partial_1^{\ast})$ 
si et seulement si :
\begin{itemize}
 \item[i)] Pour tous sommets voisins $s$ et $t$, les applications  
$\varphi_{s}$ et $\varphi_{t}$ coïncident sur $V_{\{s, t\}}$.
 \item[ii)] Pour tout $g$ dans $\mathcal{D}^{\times}$, pour tout sommet $s$, on a 
$\varphi_{s} = \varphi_{g.s} \circ \pi (g)$.
\end{itemize}
\end{propo}

\begin{démo}
Soit $\varphi = (\varphi_s)_{s \in X_0}$ une forme linéaire sur $C_0$ 
appartenant à ${\rm Hom}_{\mathcal{D}^{\times}} (C_0, \mathds{1}) \cap {\rm ker} (\partial_1^{\ast})$.
Soit $a = \{ s, t \}$ une arête et $w$ dans $V^{\mathcal{U}_{a}^1}$, alors $w \in C_1$ et 
$\partial_1 (w) = (v_u)_{u \in X_0}$ 
où $v_{t} = w$, $v_{s} = -w$ et $v_u = 0$ sinon. 
Comme $\varphi \in {\rm ker} (\partial_1^{\ast})$, 
$\varphi$ est triviale sur l'image de $\partial_1$. On a donc :
$$
\varphi (\partial_1 (w)) = 0 
= \sum_{u \in X_0} \varphi_u (v_u) = \varphi_{s} (-w) + \varphi_{t} (w) 
= -\varphi_{s} (w) + \varphi_{t} (w)
$$
On en déduit que $\varphi_{s}$ et $\varphi_{t}$ coïncident sur 
$V^{\mathcal{U}_{a}^1}$ (rappelons que 
$V^{\mathcal{U}_{a}^1} \subseteq V^{\mathcal{U}_{s}^1}$ et
$V^{\mathcal{U}_{a}^1} \subseteq V^{\mathcal{U}_{t}^1}$).\\
Soit $g$ dans $\mathcal{D}^{\times}$ 
et $(v_u)_{u \in X_0}$ dans $C_0$, alors, par $\mathcal{D}^{\times}$-équivariance, 
$\varphi (g. (v_u)_{u \in X_0}) = \varphi ((v_u)_{u \in X_0})$. 
Or $\varphi (g. (v_u)_{u \in X_0}) = \varphi ((\pi (g).v_{g^{-1}.u})_{u \in X_0})$.\\
Soit $s$ un sommet et $(v_u)_{u \in X_0}$ tel que $v_u = v \neq 0$ si $u = s$, 
$v_u = 0$ sinon. Alors :
$$
\varphi ((v_u)_{u \in X_0}) = \varphi_{s} (v) \, \, \text{et} \, \, 
\varphi (g. (v_u)_{u \in X_0}) = \varphi ((\pi (g).v_{g^{-1}.u})_{u \in X_0})
$$
or $v_{g^{-1}.u} \neq 0$ si et seulement si $g^{-1}.u = s$, 
i.e $u = g.s$, et alors 
$\pi (g).v_{g^{-1}.u} \in \pi (g).V_{s} = V_{g.s}$. 
On en déduit que 
$\varphi (g. (v_u)_{u \in X_0}) = \varphi_{g.s} (\pi (g).v)$. 
Par conséquent, pour tout $g$ dans $\mathcal{D}^{\times}$ et tout sommet $s$, 
$\varphi_{s} = \varphi_{g.s} \circ \pi (g)$.\\

\end{démo}

\section{Travaux de Silberger et Zink.}

Dans cette partie, on fixe à nouveau $(\pi, V)$ une représentation lisse irréductible de $G$, 
membre de la série discrète, non cuspidale, de niveau $0$. 
\subsection{Paramétrisation de Silberger et Zink de $\pi$.}

\begin{rmq}
Soit $s$ un sommet de $X_{\mathbb{K}}$. Alors 
$V_s = V^{\mathcal{U}_s^1}$ est stable par $\mathcal{A}_s^{\times}$ (car 
$\mathcal{U}_{s}^1 \vartriangleright \mathcal{A}_s^{\times}$) et  
l'action de $\mathcal{U}_{s}^1$ est triviale, donc 
$V_s$ peut être vu comme une représentation de 
$\mathcal{A}_s^{\times} / \mathcal{U}_s^1 \simeq {\rm GL}_2 (k_{\Delta})$. 
On remarque également que $V_s$ est un 
$\mathcal{R}_1 = \langle \varpi_{\Delta} \rangle \mathcal{A}_s^{\times}$-module.
\end{rmq}

\begin{nota}
Pour tout sommet $s$, on notera $\overline{G}_s = \mathcal{A}_s^{\times} / \mathcal{U}_s^1$ 
et $\overline{B}_s$ le sous-groupe de Borel standard de $\overline{G}_s$.
On fixe pour toute la suite $\Phi$ un générateur du groupe de Galois 
${\rm Gal} (k_{\Delta, 2} / k_{\mathbb{K}})$ 
(où $k_{\Delta, 2}$ est une extension quadratique de $k_{\Delta}$). 
\end{nota}

\begin{déf}
\begin{itemize}
\item[i)] Soit $r \in \mathbb{N}^{\ast}$. Soit $\chi$ un caractère de 
$\mathbb{K}_r^{\times}$. On dit que $\chi$ est modéré s'il est trivial sur 
$1+ \mathcal{P}_{\mathbb{K}_r}$. On dit que $\chi$ est $\mathbb{K}$-régulier si 
sa ${\rm Gal} (\mathbb{K}_r / \mathbb{K})$-orbite est de cardinal maximal, $r$.
\item[ii)] On appelle paire admissible de degré $r$ sur $\mathbb{K}$ la donnée de 
la ${\rm Gal} (\overline{\mathbb{K}} / \mathbb{K})$-orbite d'un caractère modéré 
$\mathbb{K}$-régulier de $\mathbb{K}_r^{\times}$.
\end{itemize}

\end{déf}

Nous rappelons ici la paramétrisation de Silberger et Zink de la représentation $\pi$ :

\begin{theo}\label{ParamSilbergerZinkNonCuspidales}
(\cite{SilbergerZink2}, paragraphes 2 et 3).\\ 
Pour $r$ entier naturel non nul, notons 
$\mathcal{I}_r$ l'ensemble des paires admissibles modérées de degré $r$ sur $\mathbb{K}$. 
On a alors une paramétrisation de $\mathcal{R}_0^2 (G)$ par $\sqcup_{r \vert d} \mathcal{I}_r$.\\ 
Soit $(\mathbb{K}_f, \chi_f)$ la paire admissible modérée associée à $\pi$. 
Soit : 
$$
\chi = \chi_f \circ {\rm N}_{\mathbb{K}_d / \mathbb{K}_f}
$$ 
Alors $f$ est la longueur de la 
${\rm Gal} (k_{\Delta,2} / k_{\mathbb{K}})$-orbite de $\overline{\chi}$, $e = d/f$, 
$e^{'} = (e, 2)$ (le PGCD de $e$ et $2$). 
Puisque la représentation $\pi$ n'est pas cuspidale, on a 
$e^{'} \neq 1$, dans notre cas $e^{'} = 2$, $e$ est pair et donc $f$ divise $\delta$. 
Il existe un unique caractère  
$\overline{\chi}_0$ de $k_{\Delta}^{\times}$ qui vérifie :
$$
\overline{\chi} = \overline{\chi}_0 \circ {\rm N}_{k_{\Delta,2} / k_{\Delta}}
$$
Alors, pour tout sommet $s$, la représentation 
$V_s$ de $\overline{G}_s$ se décompose simplement de la 
façon suivante :
$$
V_s \simeq ( \bigoplus_{\nu = 0}^{f-1} \overline{\chi}_0^{\Phi^{\nu}} \otimes {\rm St}_{\overline{G}_s} ) 
\oplus (\bigoplus_{0 \leq \nu_1 < \nu_2 < f} {\rm Ind}_{\overline{B}_s}^{\overline{G}_s} 
(\overline{\chi}_0^{\Phi^{\nu_1}} \otimes \overline{\chi}_0^{\Phi^{\nu_2}} ))  
$$ 
où, pour tout entier $u$, on note $\overline{\chi}_0^{\Phi^u} = \overline{\chi}_0 \circ \Phi^{u}$.
\end{theo}

Les conditions pour qu'une forme linéaire $\varphi \in {\rm Hom}_{\mathcal{D}^{\times}} (C_0, \mathds{1})$ 
appartienne au noyau de $\partial_1^{\ast}$ se lisent sur ses restrictions $\varphi_s$ aux espaces $V_s$. 
Etant donnée la décomposition de $V_s$ comme $\overline{G}_s$-module, nous allons avoir 
besoin de modèles pour les représentations de Steinberg ${\rm St}_{\overline{G}_s}$ et pour les 
induites paraboliques ${\rm Ind}_{\overline{B}_s}^{\overline{G}_s} 
(\overline{\chi}_0^{\Phi^{\nu_1}} \otimes \overline{\chi}_0^{\Phi^{\nu_2}} )$ afin de pouvoir faire les calculs.

\subsection{Modèle pour la représentation de Steinberg au niveau résiduel.}\label{ModeleSteinberg}

\begin{nota}
On fixe pour cette partie $s$ un sommet de $X_{\mathbb{K}}$ et $\nu$ un entier naturel tel que 
$0 \leq \nu < f$. 
On note $\overline{\chi}_1 = \overline{\chi}_0^{\Phi^{\nu}}$,  
$\pi_s^{\nu} = \overline{\chi}_1 \otimes {\rm St}_{\overline{G}_s}$ et $W_s^{\nu}$ l'espace de 
$\pi_s^{\nu}$. 
On notera enfin $\mathbb{P}_s^1 (k_{\Delta})$ l'ensemble des droites d'un $k_{\Delta}$-espace 
vectoriel de dimension $2$ correspondants aux arêtes de $X_{\mathbb{K}}$ dont l'un des sommets est $s$.
\end{nota}

\subsubsection{Représentation de Steinberg.}

On utilisera le modèle suivant, bien connu, pour les représentations de Steinberg (au niveau résiduel) :

\begin{propo}
On choisira le modèle suivant pour la représentation $\pi_s^{\nu}$ :
$$
W_s^{\nu} = \left\lbrace h : \mathbb{P}_s^1 (k_{\Delta}) \rightarrow \mathbb{C} : 
\sum_{x \in \mathbb{P}_s^1 (k_{\Delta})} h(x) = 0 \right\rbrace
$$
et pour tout $g$ dans $\overline{G}_s$, pour tout $h$ dans $W_s^{\nu}$ :
$$
\pi_s^{\nu} (g).h : x \mapsto \overline{\chi}_1 ({\rm det} (g)) h(g^{-1}.x)
$$
\end{propo}

\subsubsection{Dual de la Steinberg.}

On vérifie facilement la propriété suivante :

\begin{propo}
Soit $\varphi$ une forme linéaire sur $W_s^{\nu}$ (i.e. $\varphi \in (\pi_s^{\nu})^{\vee}$). 
Alors, il existe $\widetilde{\varphi}$ dans :
$$
(W_s^{\nu})^{\vee} = \{ h : \mathbb{P}_s^1 (k_{\Delta}) \rightarrow \mathbb{C} \} / 
\{ \text{fonctions constantes} \}
$$
telle que pour tout $h$ dans $(W_s^{\nu})^{\vee}$, pour tout $g$ dans $\overline{G}_s$ :
$$
(\pi_s^{\nu})^{\vee} (g).h = \overline{\chi}_1^{-1} ({\rm det} (g)) {\rm St}_{\overline{G}_s} (g).h
$$
et pour tout $h$ dans $W_s^{\nu}$ :
$$
\varphi (h) = 
\sum_{x \in \overline{G}_s} h(x) \widetilde{\varphi} (x)
$$
\end{propo}

\subsubsection{Modules de Jacquet des représentations de Steinberg.}

On a le résultat classique suivant sur les modules de Jacquet des représentations de Steinberg :

\begin{propo}
 Soit $t$ un sommet voisin de $s$. Notons $a$ l'arête reliant $s$ et $t$, 
et $d_a \in \mathbb{P}_s^1 (k_{\Delta})$ la droite correspondant à l'arête $a$. 
Alors $\overline{U} = \mathcal{U}_a^1 / \mathcal{U}_s^1$ est un radical unipotent d'un sous-groupe de Borel de 
$\overline{G}_s$ qui fixe $d_a$ et agit transitivement sur les éléments de 
$\mathbb{P}_s^1 (k_{\Delta}) \backslash \{ d_a \}$. 
On notera $(W_s^{\nu})_a$ ou encore $(W_s^{\nu})_{d_a}$ le sous-espace $(W_s^{\nu})^{\overline{U}}$ de $W_s^{\nu}$. 
On a alors un isomorphisme de $\mathbb{C}$-espaces vectoriels :
$$
(W_s^{\nu})_a \rightarrow \mathbb{C}, h \mapsto h (d_a)
$$
et :
$$
(W_s^{\nu})_a = \{ h \in W_s^{\nu} : h \, \, 
\text{est constante sur} \, \, \mathbb{P}_s^1 (k_{\Delta}) \backslash \{ d_a \} \}
$$
\end{propo}

\subsection{Modèle pour les induites paraboliques au niveau résiduel.}\label{ModeleInduiteParabolique}

\begin{nota}
Dans cette partie, on fixe $s$ un sommet de $X_{\mathbb{K}}$, $\nu_1$ et $\nu_2$ 
des entiers naturels tels que $0 \leq \nu_1 < \nu_2 \leq f-1$. 
On note $\overline{\chi}_1 = \overline{\chi}_0^{\Phi^{\nu_1}}$, 
$\overline{\chi}_2 = \overline{\chi}_0^{\Phi^{\nu_2}}$, 
$\pi_s^{\nu_1, \nu_2} = 
{\rm Ind}_{\overline{B}_s}^{\overline{G}_s} 
(\overline{\chi}_1 \otimes \overline{\chi}_2)$, et $U_s^{\nu_1, \nu_2}$ 
l'espace de~$\pi_s^{\nu_1, \nu_2}$.
\end{nota}

\subsubsection{L'induite parabolique irréductible.}

On redonne la définition de l'induite parabolique :

\begin{déf}
On a :
$$
U_s^{\nu_1, \nu_2} = \left\lbrace h : \overline{G}_s \rightarrow \mathbb{C} ;  
h(bg) = \overline{\chi}_1 \otimes \overline{\chi}_2 (b) h(g), 
b \in \overline{B}_s, g \in \overline{G}_s \right\rbrace
$$
et :
$$
\pi_s^{\nu_1, \nu_2} (g).h : x \rightarrow h(xg)
$$
pour tout $g$ dans $\overline{G}_s$ et tout $h$ dans $U_s^{\nu_1, \nu_2}$.
\end{déf}

\subsubsection{Dual de l'induite parabolique.}

On montre facilement le résultat suivant :

\begin{propo}
Soit $\varphi$ une forme linéaire sur $U_s^{\nu_1, \nu_2}$ (i.e. $\varphi \in (\pi_s^{\nu_1, \nu_2})^{\vee}$). 
Alors, il existe $\widetilde{\varphi}$ dans 
${\rm Ind}_{\overline{B}_s}^{\overline{G}_s} 
(\overline{\chi}_1^{-1} \otimes \overline{\chi}_2^{-1})$ telle que pour tout $h$ dans $U_s^{\nu_1, \nu_2}$ :
$$
\varphi (h) = 
\sum_{x \in \overline{G}_s} h(x) \widetilde{\varphi} (x)
$$
\end{propo}

\subsubsection{Modules de Jacquet de l'induite parabolique.}

On a le résultat classique suivant sur les modules de Jacquet de $\pi_s^{\nu_1, \nu_2}$ :

\begin{propo}
Soit $\overline{U}_s$ un radical unipotent de $\overline{B}_s$. 
Alors, le module de Jacquet $(\pi_s^{\nu_1, \nu_2})^{\overline{U}_s}$ est un $\mathbb{C}$-espace 
vectoriel de dimension $2$ et 
$(\pi_s^{\nu_1, \nu_2})^{\overline{U}_s} \simeq 
\overline{\chi}_1 \otimes \overline{\chi}_2 
+ \overline{\chi}_2 \otimes \overline{\chi}_1$ comme 
représentation du tore $\overline{T}_s = \overline{B}_s / \overline{U}_s$.
\end{propo}

\subsection{Lien avec le système de coefficients de Schneider et Stuhler.}

\begin{nota}\label{NotationsDecompositionPhi}
On fixe $\varphi = (\varphi_s)_{s \in X_0}$ dans 
${\rm Hom}_{\mathcal{D}^{\times}} (C_0, \mathds{1}) \cap {\rm ker} (\partial_1^{\ast})$. 
Si $s$ est un sommet de $X_{\mathbb{K}}$, alors $\varphi_s$ est une forme linéaire sur :
$$
V_s = ( \bigoplus_{\nu = 0}^{f-1} W_s^{\nu}) \oplus 
(\bigoplus_{0 \leq \nu_1 < \nu_2 < f} U_s^{\nu_1, \nu_2})
$$
On notera $\varphi_s^{\nu}$ (resp. $\varphi_s^{\nu_1, \nu_2}$) la forme linéaire sur 
$V_s$ qui coïncide avec $\varphi_s$ sur $W_s^{\nu}$ (resp. sur $U_s^{\nu_1, \nu_2}$) 
et qui est nulle ailleurs. On a alors :
$$
\varphi_s = \sum_{\nu = 0}^{f-1} \varphi_s^{\nu} + 
\sum_{0 \leq \nu_1 < \nu_2 < f} \varphi_s^{\nu_1, \nu_2}
$$
Enfin, si $g \in \mathcal{A}_s^{\times}$, on notera $\overline{g}$ sa 
réduction dans $\overline{G}_s = \mathcal{A}_s^{\times} / \mathcal{U}_s^1$.
\end{nota}

\begin{propo}
Soit $\varphi = (\varphi_s)_{s \in X_0}$ dans 
${\rm Hom}_{\mathcal{D}^{\times}} (C_0, \mathds{1}) \cap {\rm ker} (\partial_1^{\ast})$. 
Pour tout sommet $s$, pour tout $h$ dans $\mathcal{D}^{\times} \cap \mathcal{A}_s^{\times}$, on a :
$$ 
\varphi_s^{\nu} = \varphi_s^{\nu} \circ \pi_s^{\nu} (\overline{h})
$$
pour tout $\nu$ dans $\{ 0, \cdots, f-1 \}$ 
et :
$$  
\varphi_s^{\nu_1, \nu_2} = \varphi_s^{\nu_1, \nu_2} \circ \pi_s^{\nu_1, \nu_2} (\overline{h})
$$
pour tout $0 \leq \nu_1 < \nu_2 < f$. 
De plus, si $s$ et $t$ sont deux sommets voisins dans $X_{\mathbb{K}}$, en notant 
$a$ l'arête $\{ s, t \}$, on a :
$$
\varphi_{s}^{\nu} = \varphi_t^{\nu} \, \, \text{sur} \, \, 
(W_s^{\nu})^{\mathcal{U}_a^1} \simeq (W_t^{\nu})^{\mathcal{U}_a^1}
$$
pour tout $\nu$ dans $\{ 0, \cdots, f-1 \}$
et :
$$
\varphi_s^{\nu_1, \nu_2} = \varphi_t^{\nu_1, \nu_2} \, \, \text{sur} \, \, 
(U_s^{\nu_1, \nu_2})^{\mathcal{U}_a^1} \simeq (U_t^{\nu_1, \nu_2})^{\mathcal{U}_a^1}
$$
pour tout $0 \leq \nu_1 < \nu_2 < f$.
\end{propo}

\begin{démo}
Soit $\varphi = (\varphi_s)_{s \in X_0}$ dans 
${\rm Hom}_{\mathcal{D}^{\times}} (C_0, \mathds{1}) \cap {\rm ker} (\partial_1^{\ast})$. 
Fixons $s$ un sommet. Soit $v$ dans $W_s^{\nu}$ et $h$ dans 
$\mathcal{D}^{\times} \cap \mathcal{A}_s^{\times}$, alors :
$$
\varphi_s (v) = \varphi_{h.s} \circ \pi (h).v = \varphi_{s} \circ \pi (h).v
$$
Les espaces $U_s^{\nu_1, \nu_2}$ et $W_s^{\mu}$ sont deux à deux non isomorphes 
(en tant que $\overline{G}_s$-modules). Puisque $v \in W_s^{\nu}$, on a 
$\varphi_s (v) = \varphi_s^{\nu} (v)$ et, comme $\pi (h).v \in W_s^{\nu}$, 
$\varphi_{s} \circ \pi (h).v = \varphi_{s} \circ \pi_s^{\nu} (\overline{h}).v 
= \varphi_{s}^{\nu} \circ \pi_s^{\nu} (\overline{h}).v$. On en déduit~que :
$$ 
\varphi_s^{\nu} = \varphi_s^{\nu} \circ \pi_s^{\nu} (\overline{h})
$$
On montre avec un raisonnement analogue que pour tout $0 \leq \nu_1 < \nu_2 < f$ :
$$
\varphi_s^{\nu_1, \nu_2} = \varphi_s^{\nu_1, \nu_2} \circ \pi_s^{\nu_1, \nu_2} (\overline{h})
$$
Considérons à présent $s$ et $t$, deux sommets voisins de $X_{\mathbb{K}}$ et $a = \lbrace s, t \rbrace$ 
l'arête reliant ces deux sommets. 
Alors on a des injections naturelles 
$\mathcal{U}_a \hookrightarrow \mathcal{A}_s^{\times}$, $\mathcal{U}_a \hookrightarrow \mathcal{A}_t^{\times}$ et :
$$
\mathcal{U}_a / \mathcal{U}_a^1 \hookrightarrow \overline{\mathcal{A}}_s^{\times} 
= \mathcal{A}_s^{\times} / \mathcal{U}_{\mathcal{A}_s}^1 , \,
\mathcal{U}_a / \mathcal{U}_a^1 \hookrightarrow \overline{\mathcal{A}}_t^{\times}  
$$
On vérifie facilement que :
$$
(V_s)^{\mathcal{U}_a^1} = (\bigoplus_{\nu = 0}^{f-1} (W_s^{\nu})^{\mathcal{U}_a^1}) \oplus 
(\bigoplus_{0 \leq \nu_1 < \nu_2 < f} (U_s^{\nu_1, \nu_2})^{\mathcal{U}_a^1}) \, \, \text{et} \, \, 
(V_t)^{\mathcal{U}_a^1} = (\bigoplus_{\nu = 0}^{f-1} (W_t^{\nu})^{\mathcal{U}_a^1}) \oplus 
(\bigoplus_{0 \leq \nu_1 < \nu_2 < f} (U_t^{\nu_1, \nu_2})^{\mathcal{U}_a^1})
$$
où les modules de Jacquet $(W_s^{\nu})^{\mathcal{U}_a^1}$ et $(U_s^{\nu_1, \nu_2})^{\mathcal{U}_a^1}$ sont 
deux à deux non isomorphes en tant que $\mathcal{U}_a / \mathcal{U}_a^1 $-modules (de même pour les modules de Jacquet 
$(W_t^{\nu})^{\mathcal{U}_a^1}$ et $(U_t^{\nu_1, \nu_2})^{\mathcal{U}_a^1}$). 
De plus, pour $\nu \in \{ 0, \cdots f-1 \}$, on a un isomorphisme de $\mathcal{U}_a / \mathcal{U}_a^1 $-modules :
$$
(W_s^{\nu})^{\mathcal{U}_a^1} \simeq (W_t^{\nu})^{\mathcal{U}_a^1}
$$
de même, pour $0 \leq \nu_1 < \nu_2 < f$ : 
$$
(U_s^{\nu_1, \nu_2})^{\mathcal{U}_a^1} \simeq (U_t^{\nu_1, \nu_2})^{\mathcal{U}_a^1}
$$
On en déduit que pour tout $\nu \in \{ 0, \cdots f-1 \}$ :
$$
\varphi_{s}^{\nu} = \varphi_t^{\nu} \, \, \text{sur} \, \, 
(W_s^{\nu})^{\mathcal{U}_a^1} \simeq (W_t^{\nu})^{\mathcal{U}_a^1}
$$
et pour $0 \leq \nu_1 < \nu_2 < f$ :
$$
\varphi_s^{\nu_1, \nu_2} = \varphi_t^{\nu_1, \nu_2} \, \, \text{sur} \, \, 
(U_s^{\nu_1, \nu_2})^{\mathcal{U}_a^1} \simeq (U_t^{\nu_1, \nu_2})^{\mathcal{U}_a^1}
$$
\end{démo}

\subsection{Construction d'un type étendu maximal de niveau $0$ pour $\pi$.}\label{TypesEtendusMax}

Dans le théorème \ref{ParamSilbergerZinkNonCuspidales}, nous avons rappelé que les membres de la 
série discrète de niveau $0$ étaient paramétrés par certaines paires admissibles modérées. 
Afin de déterminer des conditions pour qu'une forme linéaire $\varphi$ appartienne au noyau 
de $\partial_1^{\ast}$, nous allons avoir besoin de connaître, pour un sommet~$s$, 
l'action du normalisateur de $\mathcal{A}_s^{\times}$ sur l'espace $V_s$. 
Pour cela, il nous faudra connaître le type 
étendu maximal de niveau $0$ associé à la représentation $\pi$. 
Nous allons donc tout d'abord rappeler la définition 
du type étendu maximal de niveau $0$ d'une représentation, 
au sens de Silberger et Zink. Dans l'article \cite{SilbergerZink2}, on a  
une paramétrisation des types étendus maximaux de niveau $0$ 
par des paires admissibles modérées. 
Nous allons donc rappeler quel est le lien entre la paire admissible 
modérée qui paramétrise une représentation, et celle associée à son type étendu maximal de niveau $0$. 

\subsubsection{Définition d'un type étendu maximal de niveau $0$ pour $\pi$.}

On reprend ici la définition de \cite{SilbergerZink2}, définition 0.8.\\

Un type étendu de niveau $0$ pour $\pi$ est une paire 
$(X, \Sigma)$ où $X$ est un sous-groupe ouvert compact modulo le centre de $G = A^{\times}$ et 
$\Sigma$ est une représentation irréductible de $X$ tels que :
\begin{itemize}
 \item[i)] Il existe un ordre héréditaire $\mathcal{A}$ tel que 
$\mathcal{A}^{\times} \subseteq X$ et 
$\mathds{1}_{\mathcal{U}_{\mathcal{A}}^1} \subset \Sigma_{\vert \mathcal{U}_{\mathcal{A}}^1}$.
\item[ii)] $\Sigma$ apparaît simplement dans $\pi_{\vert X}$.
\item[iii)] Si $\pi^{'} \in \mathcal{R}_0^2 (A^{\times})$ telle que 
$\Sigma \subset \pi^{'}_{\vert X}$, alors $\pi = \pi^{'}$.
\end{itemize}
Soit $\pi \in \mathcal{R}_0^2 (A^{\times})$ et on suppose que $(X, \Sigma)$ 
est un type étendu de niveau $0$ pour $\pi$. Soit $X \subseteq X^{'}$ où $X^{'}$ 
est un sous-groupe ouvert compact modulo le centre de $G$. On suppose que 
$\Sigma^{'}$ est une représentation irréductible de $X^{'}$ telle que 
$\Sigma \subset \Sigma^{'}_{\vert X}$ et $\Sigma^{'} \subset \pi_{\vert X^{'}}$. 
Alors $(X^{'}, \Sigma^{'})$ est aussi un type étendu de niveau $0$ pour $\pi$.\\
En particulier, si ${\rm c-Ind}_X^{X^{'}} \Sigma$ est irréductible alors 
$(X^{'}, {\rm c-Ind}_X^{X^{'}} \Sigma)$ est un type étendu de niveau $0$ pour~$\pi$. 
Puisque chaque sous-groupe ouvert compact modulo le centre $X$ de $G$ est contenu dans 
un sous-groupe ouvert compact modulo le centre maximal $X^{'}$ de $G$, il vient que 
pour chaque sous-groupe ouvert compact modulo le centre maximal $\mathcal{R}$ tel que
 $X \subseteq \mathcal{R}$, il existe un unique type étendu de niveau~$0$ 
$(\mathcal{R}, \widetilde{\Sigma})$ pour $\pi$ tel que 
$\Sigma \subset\widetilde{\Sigma}_{\vert X}$.\\
On dira qu'un type étendu de niveau $0$ pour $\pi$ $(\mathcal{R}, \widetilde{\Sigma})$ 
tel que $\mathcal{R}$ est un sous-groupe ouvert compact modulo le centre maximal est un 
type étendu maximal de niveau $0$ pour $\pi$.

\begin{rmq}\label{LienEntreLesParametrisations}
Plus généralement, si $G_m = {\rm GL}_m (\Delta_{\beta})$ est une forme intérieure de 
${\rm GL}_d (\mathbb{K})$ (i.e $m \times \beta = d$ et $\Delta_{\beta}$ est une 
$\mathbb{K}$-algèbre à division centrale d'indice $\beta$). 
Soit $\mathcal{P}$ l'ensemble des couples $(\tau, \tau_a)$ où 
$\tau$ est un caractère modéré de $\mathbb{K}_d^{\times}$, $a$ est la longueur de la 
${\rm Gal} (k_{\mathbb{K}, d} / k_{\mathbb{K}})$-orbite de $\overline{\tau}$ et 
$\tau_a$ est un caractère $\mathbb{K}$-régulier de $\mathbb{K}_a^{\times}$ qui vérifie 
$\tau = \tau_a \circ {\rm N}_{\mathbb{K}_d / \mathbb{K}_a}$. 
Soit $\mathcal{T} (G_m)$ les classes de conjugaison des types étendus 
maximaux de niveau $0$ des représentations dans $\mathcal{R}_0^2 (G_m)$. Alors, on a une paramétrisation :
$$
{\rm Gal} (\overline{\mathbb{K}} / \mathbb{K}) \setminus \mathcal{P} 
\rightarrow \mathcal{T} (G_m), 
(\tau, \tau_a) \mapsto \widetilde{\Sigma}_{\tau_a}^{G_m}
$$
Soit $\widetilde{\Sigma}_{\tau_a}^{G_m} \in \mathcal{T} (G_m)$  
et $\rho \in \mathcal{R}_0^2 (G_m)$. On note alors 
$\rho = \Pi (\widetilde{\Sigma}_{\tau_a}^{G_m})$ 
si $\widetilde{\Sigma}_{\tau_a}^{G_m}$ est un type étendu maximal de 
niveau $0$ pour $\rho$. 
De plus, on a aussi une paramétrisation :
$$
{\rm Gal} (\overline{\mathbb{K}} / \mathbb{K}) \setminus \mathcal{P} 
\rightarrow \mathcal{R}_0^2 (G_m), 
(\tau, \tau_a) \mapsto \Pi_{\tau_a}^{G_m}
$$
Si $(\tau, \tau_a) \in \mathcal{P}$, on a :
$$
\Pi_{\tau_a}^{G_m} = \Pi (\widetilde{w}^{m-(a,m)} \widetilde{\Sigma}_{\tau_a}^{G_m})
$$
où $w$ est un caractère non ramifié de $\mathbb{K}^{\times}$ 
(i.e. trivial sur $\mathcal{O}_{\mathbb{K}}^{\times}$) tel que 
$w^{2a} = 1$, $w^{a} \neq 1$ et 
\mbox{$\widetilde{w} = w \circ {\rm Nrd}_{{\rm M}_m (\Delta_{\beta})/ \mathbb{K}}$}.\\
Enfin, si $G_{m_1}$, $G_{m_2}$ sont deux formes intérieures de 
${\rm GL}_d (\mathbb{K})$ et si $JL : \mathcal{R}_0^2 (G_{m_1}) \rightarrow \mathcal{R}_0^2 (G_{m_2})$ 
est la correspondance de Jacquet-Langlands, alors, pour tout $(\tau, \tau_a) \in \mathcal{P}$ :
$$
JL (\Pi_{\tau_a}^{G_{m_1}}) = \Pi_{\tau_a}^{G_{m_2}}
$$
(On pourra se référer à \cite{SilbergerZink2}, Théorème 3).

\end{rmq}

\subsubsection{Construction lorsque $f$ est pair.}\label{TypeEtenduMaxfPair}

\begin{rmq}
D'après \ref{LienEntreLesParametrisations}, si $f$ est pair, on a :
$$
\pi = \Pi_{\chi_f}^{{\rm GL}_2 (\Delta_{\delta})} 
= \Pi (\widetilde{w}^{m-(f,m)} \widetilde{\Sigma}_{\chi_f}^{{\rm GL}_2 (\Delta_{\delta})}) 
= \Pi (\widetilde{w}^{2-(f,2)} \widetilde{\Sigma}_{\chi_f}^{{\rm GL}_2 (\Delta_{\delta})}) 
= \Pi (\widetilde{\Sigma}_{\chi_f}^{{\rm GL}_2 (\Delta_{\delta})}) 
$$
\end{rmq}

\noindent
Donnons à présent la construction du type étendu maximal de $\pi$ lorsque $f$ est pair 
(\cite{SilbergerZink2}, paragraphe~4). 
Soit $a$ l'arête $\{ s_0, s_1 \}$. 
Soit $(\sigma, \Sigma) = (\sigma, U_{s_0}^{0, f/2})$ la représentation de 
$\mathcal{A}_{s_0}^{\times} \langle \varpi_{\mathbb{K}} \rangle$ 
définie par :
$$
\sigma = {\rm c-Ind}_{\mathcal{A}_a^{\times} 
\langle \varpi_{\mathbb{K}} \rangle}^{\mathcal{A}_{s_0}^{\times} \langle \varpi_{\mathbb{K}} \rangle} 
(\overline{\chi}_0 \otimes \overline{\chi}_0^{\Phi^{f/2}} \chi_{\vert \mathbb{K}^{\times}}) 
$$
Alors la restriction de $\sigma$ à $\mathcal{A}_{s_0}^{\times}$ est l'induite parabolique 
$\overline{\sigma} = {\rm Ind}_{\overline{B}_{s_0}}^{\overline{G}_{s_0}} 
(\overline{\chi}_0 \otimes \overline{\chi}_0^{\Phi^{f/2}})$ 
et donc $\overline{\sigma}$ est une représentation générique. 
Soit $\overline{U} \subset \overline{\mathcal{A}}_{s_0}^{\times}$, radical unipotent.  
Soient $\varphi$ un caractère non dégénéré de $\overline{U}$ et $\varphi_0$ un caractère 
additif non trivial de $(k_{\mathbb{K}}, +)$ tels que :
$$
\varphi (u) 
= \varphi_0 \circ {\rm tr}_{k_{\Delta} / k_{\mathbb{K}}} (u_{1,2})
$$
pour tout $u = \left( \begin{array}{cc}
                    1 & u_{1,2}\\
                    0 & 1
                   \end{array} \right)
\in \overline{U}$. 
Puisque $\overline{\sigma}$ est générique, il existe un unique vecteur $v$ dans $\Sigma$ (unique 
à multiplication par un scalaire non nul près) tel que pour tout $u$ dans $\overline{U}$ :
$$
\overline{\sigma} (u).v = \varphi (u) v
$$
Comme représentation de $\mathcal{A}_{s_0}^{\times}$, on a 
$\sigma \simeq \sigma \circ {\rm Ad} (\varpi_{\Delta}^{f/2})$. 
Alors, il existe un opérateur d'entrelacement 
$J : \Sigma \rightarrow \Sigma$ tel que, pour tout $x \in \mathcal{A}_{s_0}^{\times}$, 
$\sigma (\varpi_{\Delta}^{f/2} x \varpi_{\Delta}^{-f/2}) 
= J \sigma (x) J^{-1}$. 
On fixe $J$, l'unique opérateur d'entrelacement bijectif de 
$\sigma$ dans $\sigma \circ {\rm Ad} (\varpi_{\Delta}^{f/2})$ qui vérifie 
$J .v = v$.
On définit 
$(\widehat{\sigma}, \widehat{\Sigma}) = 
(\widehat{\sigma}, U_{s_0}^{0, f/2})$
une extension de $(\sigma, \Sigma)$ à 
$N = \langle \varpi_{\Delta}^{f/2} \rangle \mathcal{A}_{s_0}^{\times}$ telle que 
pour tout $i$ dans $\mathbb{Z}$ et tout $h$ dans $\mathcal{A}_{s_0}^{\times}$ :
$$
\widehat{\sigma} ((\varpi_{\Delta}^{f/2})^i h) 
= \zeta^i J^i \sigma (h)
$$
où $\zeta = \chi_f ((-1)^{e-1} \varpi_{\mathbb{K}})$ et $e = d/f$. 
Enfin, on définit :
$$
(\widetilde{\sigma}, \widetilde{\Sigma}) 
= ({\rm c-Ind}_N^{\mathcal{R}_1} \widehat{\sigma}, 
{\rm c-Ind}_N^{\mathcal{R}_1} U_{s_0}^{0, f/2})
$$
représentation de $\mathcal{R}_1 = \langle \varpi_{\Delta} \rangle \mathcal{A}_{s_0}^{\times}$ 
(le normalisateur de $\mathcal{A}_{s_0}^{\times}$ dans $G$). Alors 
$(\widetilde{\sigma}, \widetilde{\Sigma})$ est un type étendu maximal de 
niveau $0$ pour $\pi$.
De plus, $J^e = Id_{\sigma}$.

\section{Conditions de $\mathcal{D}^{\times}$-distinction lorsque $\mathbb{K} / \mathbb{F}$ 
est totalement ramifiée.}

\begin{nota}
On fixe $(\pi, V)$ une représentation lisse et irréductible de $G$, membre de la série discrète, non cuspidale, de niveau $0$. 
Comme dans \ref{ParamSilbergerZinkNonCuspidales}, on note $(\mathbb{K}_f, \chi_f)$ la paire 
admissible modérée associée à $\pi$. 
On supposera dans toute cette partie que 
l'extension $\mathbb{K} / \mathbb{F}$ 
est totalement ramifiée, modérément ramifiée. On notera $k = k_{\mathbb{F}} \simeq k_{\mathbb{K}}$.
\end{nota}

\subsection{Conditions nécessaires de $\mathcal{D}^{\times}$-distinction de $\pi$.}

\begin{nota}\label{NotationHypotheseCNDistinction}
On suppose dans cette partie que la représentation $\pi$ est $\mathcal{D}^{\times}$-distinguée.  
On utilise les mêmes notations que dans \ref{NotationsDecompositionPhi} 
et de \ref{ParamSilbergerZinkNonCuspidales}. 
On fixe 
$\varphi = (\varphi_s)_{s \in X_0}$ dans 
${\rm Hom}_{\mathcal{D}^{\times}} (C_0, \mathds{1}) \cap {\rm ker} (\partial_1^{\ast})$. 
Si $s$ un sommet de $X_{\mathbb{K}}$, on 
notera $\varphi_s^{\nu}$ la restriction à $W_s^{\nu}$ de $\varphi_s$ et 
$\varphi_s^{\nu_1, \nu_2}$ sa restriction à $U_s^{\nu_1, \nu_2}$. On a alors :
$$
\varphi_s = \sum_{\nu = 0}^{f-1} \varphi_s^{\nu} + 
\sum_{0 \leq \nu_1 < \nu_2 < f} \varphi_s^{\nu_1, \nu_2}
$$ 
\end{nota}

Rappelons que $\overline{\chi} = \overline{\chi}_0 \circ {\rm N}_{k_{\Delta,2} / k_{\Delta}}$. 
En écrivant les condition de $\mathcal{D}^{\times}$-distinction de $\pi$ sur les éléments du centre, 
on obtient le résultat suivant :

\begin{propo}
Puisque $\pi$ est $\mathcal{D}^{\times}$-distinguée, son caractère central, $\chi$, 
est trivial sur $\mathbb{F}^{\times}$. Par suite, pour tout $x$ dans $k^{\times}$, on a :
$$
\overline{\chi}_0 (x) \in \{ -1, 1 \}
$$
\end{propo}

\subsubsection{Conditions sur $\varphi_{s_0}$.}\label{EtudeCNSommetS0}

Rappelons que l'image de $X_{\mathbb{F}}$ par $j$ est un sommet, noté $s_0$ et :
$$
V_{s_0} = ( \bigoplus_{\nu = 0}^{f-1} W_{s_0}^{\nu}) \oplus 
(\bigoplus_{0 \leq \nu_1 < \nu_2 < f} U_{s_0}^{\nu_1, \nu_2}) 
= ( \bigoplus_{\nu = 0}^{f-1} \overline{\chi}_0^{\Phi^{\nu}} \otimes {\rm St}_{\overline{G}_{s_0}} ) 
\oplus (\bigoplus_{0 \leq \nu_1 < \nu_2 < f} {\rm Ind}_{\overline{B}_{s_0}}^{\overline{G}_{s_0}} 
(\overline{\chi}_0^{\Phi^{\nu_1}} \otimes \overline{\chi}_0^{\Phi^{\nu_2}} ))
$$

\paragraph{Etude de $\varphi_{s_0}^{\nu}$ pour $0 \leq \nu \leq f-1$.} 

\begin{nota}\label{NotationsS0Steinberg}
 On fixe $\nu$ dans $\{ 0, 1, \cdots, f-1 \}$. 
On utilise le modèle présenté en \ref{ModeleSteinberg}. 
On note 
$\overline{\chi}_1 = \overline{\chi}_0^{\Phi^{\nu}}$. 
Rappelons que $\varphi_{s_0}^{\nu}$ est une forme linéaire sur 
$W_{s_0}^{\nu} = \overline{\chi}_1 \otimes {\rm St}_{\overline{G}_{s_0}}$ 
où $\overline{G}_{s_0} = \overline{\mathcal{A}}_{s_0}^{\times}$. 
Pour tout $g$ dans $\mathcal{U}_{s_0}$, on~a :
$$
\pi_{s_0}^{\nu} (g) = \overline{\chi}_1 \otimes {\rm St}_{\overline{G}_{s_0}} (g) 
= \overline{\chi}_1 ({\rm det} (\overline{g})) {\rm St}_{\overline{G}_{s_0}} (\overline{g})
$$
où $\overline{g}$ est la réduction de $g$ dans le quotient 
$\mathcal{U}_{s_0} / \mathcal{U}_{s_0}^1$. 
On fixe 
$\widetilde{\varphi}_{s_0}^{\nu}$ dans 
$\lbrace h : \mathbb{P}_{s_0}^1(k_{\Delta}) \rightarrow \mathbb{C} \rbrace / 
\lbrace \text{fonctions constantes} \rbrace$ tel que pour tout 
$h$ dans $W_{s_0}^{\nu}$ :
$$
\varphi_{s_0}^{\nu} (h) 
= \sum_{d \in \mathbb{P}_{s_0}^1(k_{\Delta})} h(d) \widetilde{\varphi}_{s_0}^{\nu} (d)
$$
On identifiera le quotient 
${\rm GL}_2 (\mathcal{O}_{\Delta}) / (1+\varpi_{\Delta} {\rm M}_2 (\mathcal{O}_{\Delta}))$ à 
$\overline{G}_{s_0} \simeq {\rm GL}_2 (k_{\Delta})$. 
On a l'inclusion
\mbox{$\mathcal{O}_{\mathcal{D}}^{\times} \subseteq \mathcal{U}_{s_0}$} 
où \mbox{$\mathcal{U}_{s_0} \simeq {\rm GL}_2 (\mathcal{O}_{\Delta})$}, on identifiera donc 
$k_{\mathcal{D}}^{\times}$ avec l'image de $\mathcal{O}_{\mathcal{D}}^{\times}$ dans le quotient 
\mbox{${\rm GL}_2 (\mathcal{O}_{\Delta}) / (1+\varpi_{\Delta} {\rm M}_2 (\mathcal{O}_{\Delta}))$}. 
On fixe $\alpha$ dans $k_{\mathcal{D}} \backslash k_{\Delta}$ tel que 
$\alpha^2 \in k_{\Delta}$. Alors $k_{\mathcal{D}} = k_{\Delta} [\alpha]$.

\end{nota}

\begin{lem}
L'injection de $k_{\mathcal{D}}$ dans ${\rm M}_2 (k_{\Delta})$ comme 
$k_{\Delta}$-algèbre est, 
à conjugaison près, donnée par :
$$
k_{\mathcal{D}} \hookrightarrow {\rm M}_2 (k_{\Delta}), \, 
x+\alpha y \mapsto \left( \begin{array}{cc} 
                     x & \alpha^2 y\\
                     y & x 
                    \end{array} \right)
$$
(où $(x, y) \in k_{\Delta}^2$). Ainsi :
$$
k_{\mathcal{D}} \simeq 
\left\lbrace \left( \begin{array}{cc} 
                     x & \alpha^2 y\\
                     y & x 
                    \end{array} \right) : 
x, y \in k_{\Delta} \right\rbrace
$$
\end{lem}

\begin{propo}
Le quotient $k_{\mathcal{D}}^{\times}  / k_{\Delta}^{\times}$ agit simplement transitivement sur $\mathbb{P}_{s_0}^1(k_{\Delta})$, 
i.e pour toutes droites $d_1$ et $d_2$ dans $\mathbb{P}_{s_0}^1(k_{\Delta})$, il existe un unique élément $x$ dans 
$k_{\mathcal{D}}^{\times}  / k_{\Delta}^{\times}$ tel que $x.d_1 = d_2$.
\end{propo}

\begin{démo}
Soit $x$ dans $k_{\mathcal{D}}^{\times}$ tel qu'il 
existe $d$ dans $\mathbb{P}_{s_0}^1 (k_{\Delta})$ tel que $x.d = d$. 
Alors, $x$ appartient à $k_{\mathcal{D}}^{\times} \subseteq {\rm GL}_2 (k_{\Delta})$ et possède au 
moins une valeur propre. Puisque l'extension $k_{\mathcal{D}} / k_{\Delta}$ est quadratique, 
$x$~appartient à $k_{\Delta}^{\times}$. 
Soit $d$ une droite de $\mathbb{P}_{s_0}^1 (k_{\Delta})$. Notons $\mathcal{O} (d)$ l'orbite de 
$d$ sous l'action de $k_{\mathcal{D}}^{\times} / k_{\Delta}^{\times}$ et 
${\rm Stab} (d)$ l'ensemble $\{ x \in k_{\mathcal{D}}^{\times} / k_{\Delta}^{\times} : x.d = d \}$. 
Alors d'après ce qui précède, ${\rm Stab} (d) = \{ \overline{Id} \}$. Or :
$$
\sharp \mathcal{O} (d) = \sharp ( k_{\mathcal{D}}^{\times} / k_{\Delta}^{\times}) 
/ \sharp ({\rm Stab} (d)) = \sharp ( k_{\mathcal{D}}^{\times} / k_{\Delta}^{\times}) 
= \frac{Q^2 -1}{Q-1} = Q+1 = \sharp \mathbb{P}_{s_0}^1 (k_{\Delta})
$$
On en déduit que l'action est transitive, et même simplement transitive.

\end{démo}

\begin{lem}
Soit $g$ dans $k_{\mathcal{D}}^{\times}$. 
Le fait que ${\rm N}_{k_{\mathcal{D}} / k} (g)$ soit ou non un carré dans $k^{\times}$ ne dépend 
que de la classe de $g$ dans $k_{\mathcal{D}}^{\times} / k_{\Delta}^{\times}$.
\end{lem}

\begin{démo}
Il suffit de remarquer que si $x$ appartient à $k_{\Delta}^{\times}$, alors :
$$
{\rm N}_{k_{\mathcal{D}} / k} (x) 
= {\rm N}_{k_{\Delta} / k} ({\rm N}_{k_{\mathcal{D}} / k_{\Delta}} (x)) 
= ({\rm N}_{k_{\Delta} / k} (x))^2
$$

\end{démo}

\begin{nota}\label{NotationsS0SteinbergDroites}
Notons $t_1, \cdots, t_{Q+1}$ les sommets voisins de $s_0$ (dans $X_{\mathbb{K}}$). 
Pour $i$ dans $\{ 1, \cdots, Q+1 \}$, on notera $a_i$ l'arête $\{ s_0, t_i \}$ et 
$\delta_i^{s_0} \in \mathbb{P}^1_{s_0} (k_{\Delta})$ la droite correspondante. 
On supposera que $t_1 = s_1$, ainsi la droite $\delta_1^{s_0}$ correspond à l'arête 
$\{ s_0, s_1 \}$. 
Puisque  
$k_{\mathcal{D}}^{\times} / k_{\Delta}^{\times}$ agit simplement transitivement sur 
$\mathbb{P}^1_{s_0} (k_{\Delta})$ 
et que $t_1, \cdots, t_{Q+1}$ sont tous dans la même 
$\mathcal{O}_{\mathcal{D}}^{\times}$-orbite, on peut fixer   
$g_i$ dans $\mathcal{O}_{\mathcal{D}}^{\times} \subseteq \mathcal{U}_{s_0}$ et noter $\overline{g}_i$ 
son image dans $\mathcal{U}_{s_0} / \mathcal{U}_{s_0}^1 = \overline{G}_{s_0}$,  
tel que $\overline{g}_i. \delta_i^{s_0} = \delta_1^{s_0}$ (et on a donc 
$g_i.s_0 = s_0$), d'où : 
$$
\varphi_{s_0}^{\nu} = \varphi_{s_0}^{\nu} \circ \pi_{s_0}^{\nu} (g_i)
$$
et :
$$
\pi_{s_0}^{\nu} (g_i) = \overline{\chi}_1 ({\rm N}_{k_{\mathcal{D}} / k} (\overline{g}_i)) 
{\rm St}_{\overline{G}_{s_0}} (\overline{g}_i)
$$
Pour tout $i$ dans $\lbrace 1, \cdots, Q+1 \rbrace$, on notera 
$\alpha_i = \overline{\chi}_1 ({\rm N}_{k_{\mathcal{D}} / k} (\overline{g}_i))$. 
Comme  
${\rm N}_{k_{\mathcal{D}} / k} (\overline{g}_i) \in k^{\times}$, on a 
\mbox{$\alpha_i \in \{ -1, 1 \}$}.
\end{nota}

\begin{propo}
On se place sous les hypothèses décrites en \ref{NotationHypotheseCNDistinction} 
et \ref{NotationsS0Steinberg}. 
Si $\overline{\chi}_0$ est trivial sur $k^{\times}$, 
alors $\varphi_{s_0}^{\nu} = 0$.
\end{propo}

\begin{démo} 
Comme $\overline{\chi}_0$ est trivial sur $k^{\times}$, on a 
$\alpha_i = 1$ pour tout $1 \leq i \leq Q+1$. 
De plus, soit $h_{\delta_i^{s_0}}$ dans $W_{s_0}^{\nu}$ telle que  
$h_{\delta_i^{s_0}} (\delta) = - \frac{1}{Q}$ si $\delta \neq \delta_i^{s_0}$ 
et $h_{\delta_i^{s_0}} (\delta_i^{s_0}) = 1$. Alors :
$$
\varphi_{s_0}^{\nu} \circ \pi (g_i) (h_{\delta_i^{s_0}}) 
= \varphi_{s_0}^{\nu} ({\rm St}_{\overline{G}_{s_0}} (\overline{g}_i).h_{\delta_i^{s_0}}) 
$$
où : 
$$
{\rm St}_{\overline{G}_{s_0}} (\overline{g}_i).h_{\delta_i^{s_0}} : \delta_1^{s_0} \mapsto 1, \, 
\delta \neq \delta_1^{s_0} \mapsto - \frac{1}{Q}
$$
Ainsi :
$$
\varphi_{s_0}^{\nu} \circ \pi (g_i) (h_{\delta_i^{s_0}})
= (1+ \frac{1}{Q}) \widetilde{\varphi}_{s_0}^{\nu} (\delta_1^{s_0}) - \frac{1}{Q} 
\sum_{\delta \in \mathbb{P}_{s_0}^1 (k_{\Delta})} \widetilde{\varphi}_{s_0}^{\nu} (\delta)
$$
et :
$$
\varphi_{s_0}^{\nu} (h_{\delta_i^{s_0}}) 
= \widetilde{\varphi}_{s_0}^{\nu} (\delta_i^{s_0}) - \frac{1}{Q} 
\sum_{\delta \neq \delta_i^{s_0}} \widetilde{\varphi}_{s_0}^{\nu} (\delta) 
= (1+ \frac{1}{Q}) \widetilde{\varphi}_{s_0}^{\nu} (\delta_i^{s_0}) - \frac{1}{Q} 
\sum_{\delta \in \mathbb{P}_{s_0}^1 (k_{\Delta})} \widetilde{\varphi}_{s_0}^{\nu} (\delta) 
$$
On en déduit que pour tout $1 \leq i \leq Q+1$, on a 
$\widetilde{\varphi}_{s_0}^{\nu} (\delta_1^{s_0}) 
= \widetilde{\varphi}_{s_0}^{\nu} (\delta_i^{s_0})$. 
Par suite, $\widetilde{\varphi}_{s_0}^{\nu}$ est une fonction constante, et 
$\varphi_{s_0}^{\nu} = 0$.

\end{démo}

\begin{lem}
On reprend les notations de \ref{NotationsS0SteinbergDroites}. 
 Supposons que $\overline{\chi}_0$ vaut $-1$ sur les éléments de $k^{\times}$ qui ne sont pas 
des carrés et vaut $1$ sur les carrés de $k^{\times}$. 
Soit $r$ le nombre d'indices $i$ tels que $\alpha_i = 1$.\\
On a $r = \frac{Q+1}{2}$.
\end{lem}

\begin{démo}
Rappelons que $k_{\mathcal{D}}^{\times} / k_{\Delta}^{\times}$ agit simplement transitivement sur 
$\mathbb{P}_{s_0}^1 (k_{\Delta})$. Le quotient $k_{\mathcal{D}}^{\times} / k_{\Delta}^{\times}$ 
est de cardinal $Q+1$ et, l'action de $k_{\mathcal{D}}^{\times} / k_{\Delta}^{\times}$ 
étant simplement transitive, les $\overline{g}_i$ sont deux à deux distincts. Ainsi :
$$
k_{\mathcal{D}}^{\times} / k_{\Delta}^{\times} 
= \{ \overline{g}_i : 1 \leq i \leq Q+1 \}
$$ 
La norme ${\rm N}_{k_{\mathcal{D}} / k} : k_{\mathcal{D}}^{\times} \rightarrow k^{\times}$ 
est surjective. Soit $(k^{\times})^2$ l'ensemble des carrés de $k^{\times}$, alors l'application :
$$
\overline{N}_0 : k_{\mathcal{D}}^{\times} \rightarrow k^{\times} / (k^{\times})^2, 
x \mapsto {\rm N}_{k_{\mathcal{D}} / k} (x) (k^{\times})^2
$$
est encore surjective. De plus, pour tout $x$ dans $k_{\Delta}^{\times}$,
${\rm N}_{k_{\mathcal{D}} / k} (x) = ({\rm N}_{k_{\Delta} / k} (x))^2$ est un carré, et donc 
$\overline{N}_0$ passe au quotient. On en déduit une surjection :
$$
\overline{N} : k_{\mathcal{D}}^{\times} / k_{\Delta}^{\times} \rightarrow 
k^{\times} / (k^{\times})^2, \overline{g}_i \mapsto {\rm N}_{k_{\mathcal{D}} / k} (\overline{g}_i) (k^{\times})^2
$$
Ainsi ${\rm ker} (\overline{N})$ est de cardinal $\frac{Q+1}{2}$.  
De plus, $\overline{g}_i \in {\rm ker} (\overline{N})$ si et seulement si 
${\rm N}_{k_{\mathcal{D}} / k} (\overline{g}_i)$ est un carré dans $k^{\times}$, 
i.e, si $\alpha_i = 1$. Par conséquent $r$ est le cardinal de ${\rm ker} (\overline{N})$.

\end{démo}

\begin{propo}
Supposons que $\overline{\chi}_0$ vaut $-1$ sur les éléments de $k^{\times}$ qui ne sont pas 
des carrés et vaut $1$ sur les carrés de $k^{\times}$. \\
Si $\widetilde{\varphi}_{s_0}^{\nu} (\delta_1^{s_0}) 
= \widetilde{\varphi}_{s_0}^{\nu} (\delta_{Q+1}^{s_0})$ alors 
$\varphi_{s_0}^{\nu} = 0$.\\
Si $\widetilde{\varphi}_{s_0}^{\nu} (\delta_1^{s_0}) \neq \widetilde{\varphi}_{s_0}^{\nu} (\delta_{Q+1}^{s_0})$ 
alors pour tout $h$ dans $W_{s_0}^{\nu}$ :
$$
\varphi_{s_0}^{\nu} (h) = \left( \sum_{i=1}^{(Q+1)/2} h(\delta_i^{s_0}) \right)
\times (\widetilde{\varphi}_{s_0}^{\nu} (\delta_1^{s_0}) - \widetilde{\varphi}_{s_0}^{\nu} (\delta_{Q+1}^{s_0}))
$$
et : 
$$
{\rm ker} (\varphi_{s_0}^{\nu}) 
= \lbrace h \in W_{s_0}^{\nu} : \sum_{i=1}^{(Q+1)/2} h(\delta_i^{s_0}) = 0 \rbrace
$$
\end{propo}

\begin{démo}
Quitte à réordonner les droites, on peut supposer que 
$\alpha_1 = \alpha_2 = \cdots = \alpha_{\frac{Q+1}{2}} = 1$. 
Alors, avec des calculs analogues au cas précédent, on a :
$$
\widetilde{\varphi}_{s_0}^{\nu} (\delta_1^{s_0}) = \widetilde{\varphi}_{s_0}^{\nu} (\delta_2^{s_0}) 
= \cdots = \widetilde{\varphi}_{s_0}^{\nu} (\delta_{\frac{Q+1}{2}}^{s_0})
$$
(et, comme $\overline{g}_1 = \overline{Id}$ 
dans $k_{\mathcal{D}}^{\times} / k_{\Delta}^{\times}$, 
on peut encore supposer, après ré-indexation, que $\delta_1^{s_0}$ correspond à l'arête 
$\{ s_0, s_1 \}$).
Puis, pour $i$ dans $\lbrace \frac{Q+3}{2}, \cdots, Q+1 \rbrace$ :
$$
\varphi_{s_0}^{\nu} (h_{\delta_i^{s_0}}) = 
(1+ \frac{1}{Q}) \widetilde{\varphi}_{s_0}^{\nu} (\delta_i^{s_0}) - \frac{1}{Q} 
\sum_{\delta \in \mathbb{P}_{s_0}^1 (k_{\Delta})} \widetilde{\varphi}_{s_0}^{\nu} (\delta) 
$$
et :
$$
\varphi_{s_0}^{\nu} \circ \pi (g_i) (h_{\delta_i^{s_0}}) = 
-(1+ \frac{1}{Q}) \widetilde{\varphi}_{s_0}^{\nu} (\delta_1^{s_0}) + \frac{1}{Q} 
\sum_{\delta \in \mathbb{P}_{s_0}^1 (k_{\Delta})} \widetilde{\varphi}_{s_0}^{\nu} (\delta) 
$$
Ainsi :
$$
( \frac{Q+1}{Q}) \widetilde{\varphi}_{s_0}^{\nu} (\delta_i^{s_0}) 
= - ( \frac{Q+1}{Q}) \widetilde{\varphi}_{s_0}^{\nu} (\delta_1^{s_0}) 
+ \frac{2}{Q} \sum_{\delta \in \mathbb{P}_{s_0}^1 (k_{\Delta})} \widetilde{\varphi}_{s_0}^{\nu} (\delta) 
$$
d'où :
$$
\widetilde{\varphi}_{s_0}^{\nu} (\delta_i^{s_0}) 
= - \widetilde{\varphi}_{s_0}^{\nu} (\delta_1^{s_0}) 
+ \frac{2}{Q+1} \sum_{\delta \in \mathbb{P}_{s_0}^1 (k_{\Delta})} \widetilde{\varphi}_{s_0}^{\nu} (\delta)
$$
On en déduit que :
$$
\widetilde{\varphi}_{s_0}^{\nu} (\delta_{\frac{Q+3}{2}}^{s_0}) = \cdots 
= \widetilde{\varphi}_{s_0}^{\nu} (\delta_{Q+1}^{s_0}) =  - \widetilde{\varphi}_{s_0}^{\nu} (\delta_1^{s_0}) 
+ \frac{2}{Q+1} \sum_{\delta \in \mathbb{P}_{s_0}^1 (k_{\Delta})} \widetilde{\varphi}_{s_0}^{\nu} (\delta)
$$
Ainsi, pour tout $h$ dans $W_{s_0}^{\nu}$ :
$$
\varphi_{s_0}^{\nu} (h) 
= \sum_{i=1}^{Q+1} h(\delta_i^{s_0}) \widetilde{\varphi}_{s_0}^{\nu} (\delta_i^{s_0}) 
= \widetilde{\varphi}_{s_0}^{\nu} (\delta_1^{s_0}) \times \sum_{i=1}^{(Q+1)/2} h(\delta_i^{s_0}) 
+ \widetilde{\varphi}_{s_0}^{\nu} (\delta_{Q+1}^{s_0}) \times \sum_{i=(Q+3)/2}^{Q+1} h(\delta_i^{s_0})
$$
Or : 
$$
\sum_{i=1}^{Q+1} h(\delta_i^{s_0}) = 0 
\Rightarrow \sum_{i=(Q+3)/2}^{Q+1} h(\delta_i^{s_0}) = -  \sum_{i=1}^{(Q+1)/2} h(\delta_i^{s_0}) 
$$
Par suite :
$$
\varphi_{s_0}^{\nu} (h) = \left( \sum_{i=1}^{(Q+1)/2} h(\delta_i^{s_0}) \right)
\times (\widetilde{\varphi}_{s_0}^{\nu} 
(\delta_1^{s_0}) - \widetilde{\varphi}_{s_0}^{\nu} (\delta_{Q+1}^{s_0}))
$$
Finalement, ou bien $\widetilde{\varphi}_{s_0}^{\nu} (\delta_1^{s_0}) 
= \widetilde{\varphi}_{s_0}^{\nu} (\delta_{Q+1}^{s_0})$, et dans ce cas 
$\varphi_{s_0}^{\nu} = 0$; ou bien  
$\widetilde{\varphi}_{s_0}^{\nu} (\delta_1^{s_0}) \neq \widetilde{\varphi}_{s_0}^{\nu} (\delta_{Q+1}^{s_0})$ 
et alors, pour tout $h$ dans $W_{s_0}^{\nu}$, on a :
$$
\varphi_{s_0}^{\nu} (h) = \left( \sum_{i=1}^{(Q+1)/2} h(\delta_i^{s_0}) \right)
\times (\widetilde{\varphi}_{s_0}^{\nu} 
(\delta_1^{s_0}) - \widetilde{\varphi}_{s_0}^{\nu} (\delta_{Q+1}^{s_0}))
$$
et : 
$$
{\rm ker} (\varphi_{s_0}^{\nu}) 
= \lbrace h \in W_{s_0}^{\nu} : \sum_{i=1}^{(Q+1)/2} h(\delta_i^{s_0}) = 0 \rbrace
$$
\end{démo}

\begin{rmq}
On peut réécrire cette formule : notons $\mathcal{O}_1$ l'ensemble des droites 
$\delta_i^{s_0}$ telles que $\alpha_i = {\rm N}_{k_{\mathcal{D}} / k} (\overline{g}_i)$ est un carré 
dans $k^{\times}$ et $\mathcal{O}_2 = \mathbb{P}_{s_0}^1 (k_{\Delta}) \backslash \mathcal{O}_1$. 
Alors $\varphi_{s_0}^{\nu}$ est proportionnelle à :
$$
h \mapsto \sum_{d \in \mathcal{O}_1} h(d) = - \sum_{d \in \mathcal{O}_2} h(d)
$$
\end{rmq}

\begin{propo}\label{DimensionEntrelacementSteinberg}
On se place sous les hypothèses \ref{NotationHypotheseCNDistinction} avec les notations 
\ref{NotationsS0Steinberg} et \ref{NotationsS0SteinbergDroites}. 
L'espace d'entrelacement 
${\rm Hom}_{k_{\mathcal{D}}^{\times}} (\pi_{s_0}^{0}, \mathds{1})$ est de dimension au plus $1$. 
Plus précisément, si $\overline{\chi}_0$ est trivial sur $k^{\times}$, 
${\rm Hom}_{k_{\mathcal{D}}^{\times}} (\pi_{s_0}^{0}, \mathds{1}) = \lbrace 0 \rbrace$, sinon, 
${\rm Hom}_{k_{\mathcal{D}}^{\times}} (\pi_{s_0}^{0}, \mathds{1})$ est de dimension $1$.
\end{propo}

\begin{démo}
On a :
$$
{\rm dim}_{\mathbb{C}} ( {\rm Hom}_{k_{\mathcal{D}}^{\times}} (\pi_{s_0}^{0}, \mathds{1})) 
= \langle \chi_{\pi_{s_0}^{0}}, \chi_{\mathds{1}_{k_{\mathcal{D}}^{\times}}} \rangle
$$
où $\chi_{\pi_{s_0}^{0}}$ est le caractère de la représentation $\pi_{s_0}^{0}$.  
Rappelons que :
$$
\chi_{{\rm St}_{\overline{G}_{s_0}}} 
= \chi_{{\rm Ind}_{\overline{B}_{s_0}}^{\overline{G}_{s_0}} \mathds{1}} 
- \mathds{1}_{\overline{G}_{s_0}}
$$
Remarquons que $k_{\mathcal{D}}^{\times} 
= k_{\Delta}^{\times} \sqcup (k_{\mathcal{D}}^{\times} \setminus k_{\Delta}^{\times})$ 
où tous les éléments de $k_{\Delta}^{\times}$ sont centraux, et tous les éléments de 
$k_{\mathcal{D}}^{\times} \setminus k_{\Delta}^{\times}$ sont elliptiques réguliers. 
Par conséquent, si $x$ appartient à $k_{\Delta}^{\times}$, on a :
$$
\chi_{{\rm St}_{\overline{G}_{s_0}}} (x) 
= {\rm dim} ({\rm St}_{\overline{G}_{s_0}}) 
= Q+1-1 = Q
$$
et, si $x$ appartient à $k_{\mathcal{D}}^{\times} \setminus k_{\Delta}^{\times}$, on a :
$$
\chi_{{\rm St}_{\overline{G}_{s_0}}} (x) 
= \chi_{{\rm Ind}_{\overline{B}_{s_0}}^{\overline{G}_{s_0}} \mathds{1}} (x)
- \mathds{1}_{\overline{G}_{s_0}} (x) 
= 0 - 1 = -1
$$
On en déduit que :
\begin{eqnarray*}
\langle \chi_{\pi_{s_0}^{0}}, \chi_{\mathds{1}_{k_{\mathcal{D}}^{\times}}} \rangle 
& = & \frac{1}{\vert k_{\mathcal{D}}^{\times} \vert} 
       \sum_{x \in k_{\mathcal{D}}^{\times}} \chi_{\overline{\chi}_0 \otimes {\rm St}_{\overline{G}_{s_0}}} (x) 
       = \frac{1}{\vert k_{\mathcal{D}}^{\times} \vert} ( \sum_{x \in k_{\Delta}^{\times}} 
          \overline{\chi}_0 \circ {\rm N}_{k_{\mathcal{D}} / k} (x) \times Q 
          - \sum_{x \in k_{\mathcal{D}}^{\times} \setminus k_{\Delta}^{\times}} 
          \overline{\chi}_0 \circ {\rm N}_{k_{\mathcal{D}} / k} (x))\\
& = &  \frac{1}{\vert k_{\mathcal{D}}^{\times} \vert} ( (Q+1) \times \sum_{x \in k_{\Delta}^{\times}} 
          \overline{\chi}_0 \circ {\rm N}_{k_{\mathcal{D}} / k} (x) 
          - \sum_{x \in k_{\mathcal{D}}^{\times}} 
          \overline{\chi}_0 \circ {\rm N}_{k_{\mathcal{D}} / k} (x))
\end{eqnarray*}
Par hypothèse, $\overline{\chi}_0^2$ est trivial sur $k^{\times}$, ainsi :
$$
\sum_{x \in k_{\Delta}^{\times}} 
          \overline{\chi}_0 \circ {\rm N}_{k_{\mathcal{D}} / k} (x) 
= \sum_{x \in k_{\Delta}^{\times}} 
          \overline{\chi}_0 (x^2) = \vert k_{\Delta}^{\times} \vert = Q-1
$$
On distingue alors deux cas.\\
Supposons dans un premier temps que $\overline{\chi}_0$ est trivial sur $k^{\times}$, alors :
$$
\langle \chi_{\pi_{s_0}^{0}}, \chi_{\mathds{1}_{k_{\mathcal{D}}^{\times}}} \rangle 
= \frac{1}{\vert k_{\mathcal{D}}^{\times} \vert} ((Q+1)(Q-1)- \vert k_{\mathcal{D}}^{\times} \vert) 
= \frac{1}{\vert k_{\mathcal{D}}^{\times} \vert} ((Q+1)(Q-1)- (Q^2 -1)) = 0
$$
et :
$$
{\rm Hom}_{k_{\mathcal{D}}^{\times}} (\pi_{s_0}^{0}, \mathds{1}) 
= \lbrace 0 \rbrace
$$
Sinon, $\overline{\chi}_0$ est trivial sur les carrés de $k^{\times}$ et vaut $-1$ sur les éléments de 
$k^{\times}$ qui ne sont pas des carrés. Alors $ \overline{\chi}_0 \circ {\rm N}_{k_{\mathcal{D}} / k}$ est non trivial sur 
$k_{\mathcal{D}}^{\times}$ et :
$$
\langle \chi_{\pi_{s_0}^{0}}, \chi_{\mathds{1}_{k_{\mathcal{D}}^{\times}}} \rangle 
= \frac{1}{\vert k_{\mathcal{D}}^{\times} \vert} ((Q+1)(Q-1)- 0) = 1 
= {\rm dim}_{\mathbb{C}} ({\rm Hom}_{k_{\mathcal{D}}^{\times}} (\pi_{s_0}^{0}, \mathds{1}) )
$$

\end{démo}

\paragraph{Etude de $\varphi_{s_0}^{\nu_1, \nu_2}$ pour $0 \leq \nu_1 < \nu_2 \leq f-1$.}

\begin{nota}\label{NotationsS0Induites}
On se place à nouveau sous les hypothèses de \ref{NotationHypotheseCNDistinction}. 
On fixe deux entiers naturels $\nu_1$ et $\nu_2$ tels que 
$0 \leq \nu_1 < \nu_2 \leq f-1$. 
Nous noterons 
$\overline{\chi}_1 = \overline{\chi}_0^{\Phi^{\nu_1}}$, 
$\overline{\chi}_2 = \overline{\chi}_0^{\Phi^{\nu_2}}$. 
Soit \mbox{$\pi_{s_0}^{\nu_1, \nu_2} = U_{s_0}^{\nu_1, \nu_2} 
= {\rm Ind}_{\overline{B}_{s_0}}^{\overline{G}_{s_0}} 
(\overline{\chi}_1 \otimes \overline{\chi}_2 )$}. 
Puisque $\overline{\chi}_0^2$ est trivial sur $k^{\times}$, nécessairement, 
$\overline{\chi}_1 \times \overline{\chi}_2$ est également trivial sur $k^{\times}$. 
On utilise le modèle présenté en \ref{ModeleInduiteParabolique}. 
On fixe 
\mbox{$\widetilde{\varphi}_{s_0}^{\nu_1, \nu_2} \in {\rm Ind}_{\overline{B}_{s_0}}^{\overline{G}_{s_0}} 
(\overline{\chi}_1^{-1} \otimes \overline{\chi}_2^{-1} )$} tel que :
$$ 
\varphi_{s_0}^{\nu_1, \nu_2} (h) 
= \sum_{x \in \overline{G}_{s_0}} \widetilde{\varphi}_{s_0}^{\nu_1, \nu_2} (x) h(x)
$$
pour tout $h$ dans $\pi_{s_0}^{\nu_1, \nu_2}$.
\end{nota}

\begin{propo}
Si 
$\overline{\chi}_1 \times \overline{\chi}_2$
est non trivial sur $k_{\Delta}^{\times}$ alors 
$\varphi_{s_0}^{\nu_1, \nu_2} = 0$.
\end{propo}

\begin{démo}
Pour tout $g$ dans $k_{\mathcal{D}}^{\times}$, on a 
$\varphi_{s_0}^{\nu_1, \nu_2} 
= \varphi_{s_0}^{\nu_1, \nu_2} \circ \pi_{s_0}^{\nu_1, \nu_2} (g)$. 
Ainsi 
$\varphi_{s_0}^{\nu_1, \nu_2}$ appartient à  
${\rm Hom}_{k_{\mathcal{D}}^{\times}} (\pi_{s_0}^{\nu_1, \nu_2}, \mathds{1})$. 
Nous allons déterminer les conditions pour que l'espace d'entrelacement 
${\rm Hom}_{k_{\mathcal{D}}^{\times}} (\pi_{s_0}^{\nu_1, \nu_2}, \mathds{1})$ 
soit non trivial, i.e pour que 
$\mathds{1}_{k_{\mathcal{D}}^{\times}} \subseteq (\pi_{s_0}^{\nu_1, \nu_2})_{\vert k_{\mathcal{D}}^{\times}}$. 
Le caractère $\mathds{1}_{k_{\mathcal{D}}^{\times}}$ est une sous-représentation de 
$(\pi_{s_0}^{\nu_1, \nu_2})_{\vert k_{\mathcal{D}}^{\times}}$ si et seulement si le produit scalaire des 
caractères $\langle \chi_{\pi_{s_0}^{\nu_1, \nu_2}}, \mathds{1}_{k_{\mathcal{D}}^{\times}} \rangle$ 
est non nul. 
Rappelons que pour tout $g$~dans~$\overline{G}_{s_0}$, on~a :
$$
\chi_{\pi_{s_0}^{\nu_1, \nu_2}} (g) 
= \frac{1}{\vert \overline{B}_{s_0} \vert} 
\sum_{\substack{x \in \overline{G}_{s_0} \\ x^{-1} g x \in \overline{B}_{s_0}}} 
\overline{\chi}_1 \otimes \overline{\chi}_2 (x^{-1} g x)
$$
Ainsi : 
\begin{eqnarray*}
\langle \chi_{\pi_{s_0}^{\nu_1, \nu_2}}, \mathds{1}_{k_{\mathcal{D}}^{\times}} \rangle 
& = & \frac{1}{\vert k_{\mathcal{D}}^{\times} \vert} 
       \sum_{g \in k_{\mathcal{D}}^{\times}} \chi_{\pi_{s_0}^{\nu_1, \nu_2}} (g)\\
& = & \frac{1}{\vert k_{\mathcal{D}}^{\times} \vert \times \vert \overline{B}_{s_0} \vert} 
      \sum_{g \in k_{\mathcal{D}}^{\times}} 
      \sum_{\substack{x \in \overline{G}_{s_0} \\ x^{-1} g x \in \overline{B}_{s_0}}} 
       \overline{\chi}_1 \otimes \overline{\chi}_2 (x^{-1} g x)\\
\end{eqnarray*}
On remarque que $k_{\mathcal{D}}^{\times} 
= (k_{\mathcal{D}}^{\times} \backslash k_{\Delta}^{\times}) 
\sqcup k_{\Delta}^{\times}$ où tous les éléments de 
$k_{\mathcal{D}}^{\times} \backslash k_{\Delta}^{\times}$ sont elliptiques réguliers 
(donc aucun de leurs conjugués n'est contenu dans un sous-groupe de Borel de $\overline{G}_{s_0}$) et 
tous les éléments de $k_{\Delta}^{\times}$ sont centraux dans $\overline{G}_{s_0}$. 
Par conséquent :
\begin{eqnarray*}
\langle \chi_{\pi_{s_0}^{\nu_1, \nu_2}}, \mathds{1}_{k_{\mathcal{D}}^{\times}} \rangle 
& = & \frac{1}{\vert k_{\mathcal{D}}^{\times} \vert \times \vert \overline{B}_{s_0} \vert} 
      \sum_{g \in k_{\Delta}^{\times}} 
      \sum_{\substack{x \in \overline{G}_{s_0} \\ x^{-1} g x \in \overline{B}_{s_0}}} 
       \overline{\chi}_1 \otimes \overline{\chi}_2 (x^{-1} g x)\\
& = & \frac{\vert \overline{G}_{s_0} \vert}{\vert k_{\mathcal{D}}^{\times} \vert \times \vert \overline{B}_{s_0} \vert} 
      \sum_{g \in k_{\Delta}^{\times}}  
       \overline{\chi}_1 (g) \times \overline{\chi}_2 (g)\\
\end{eqnarray*}
On en déduit que 
$\langle \chi_{\pi_{s_0}^{\nu_1, \nu_2}}, \mathds{1}_{k_{\mathcal{D}}^{\times}} \rangle 
\neq 0$ si et seulement si $(\overline{\chi}_1 \times \overline{\chi}_2)_{\vert k_{\Delta}^{\times}} 
= \mathds{1}_{\vert k_{\Delta}^{\times}}$.

\end{démo}

\begin{propo}\label{DimensionEntrelacementInduite}
On utilise les notations de \ref{NotationsS0Induites}. 
Supposons que $\overline{\chi}_1 \times \overline{\chi}_2$ 
est trivial sur $k_{\Delta}^{\times}$. 
L'espace d'entrelacement 
${\rm Hom}_{k_{\mathcal{D}}^{\times}} (\pi_{s_0}^{\nu_1, \nu_2}, \mathds{1})$ est de dimension $1$. 
Par conséquent, 
la forme linéaire $\varphi_{s_0}^{\nu_1, \nu_2}$, qui appartient à 
${\rm Hom}_{k_{\mathcal{D}}^{\times}} (\pi_{s_0}^{\nu_1, \nu_2}, \mathds{1})$, 
est unique, à multiplication par un scalaire près.
\end{propo}

\begin{démo}
Avec un calcul rapide, on remarque que 
$\frac{\vert \overline{G}_{s_0} \vert}{\vert \overline{B}_{s_0} \vert} 
= Q+1 = \frac{\vert k_{\mathcal{D}}^{\times} \vert}{\vert k_{\Delta}^{\times} \vert}$. 
Par conséquent, si $\overline{\chi}_1 \times \overline{\chi}_2$ 
est trivial sur $k_{\Delta}^{\times}$, on a 
$\langle \chi_{\pi_{s_0}^{\nu_1, \nu_2}}, \mathds{1}_{k_{\mathcal{D}}^{\times}} \rangle 
= \frac{\vert \overline{G}_{s_0} \vert}{\vert \overline{B}_{s_0} \vert} \times 
\frac{\vert k_{\Delta}^{\times} \vert}{\vert k_{\mathcal{D}}^{\times} \vert} = 1$.

\end{démo}

Rappelons que $\overline{\chi} = \overline{\chi}_0 \circ {\rm N}_{k_{\Delta, 2} / k_{\Delta}}$ et que la 
${\rm Gal} (k_{\Delta, 2} / k)$-orbite de $\overline{\chi}$ est de longueur $f$. On montre facilement 
le lemme suivant :

\begin{lem}
Soit $\langle \Phi \rangle = {\rm Gal} (k_{\Delta, 2} / k)$. 
Alors, tout comme $\overline{\chi}$, la $\langle \Phi \rangle$-orbite de 
$\overline{\chi}_0$ est de longueur $f$.
\end{lem}

\begin{propo}
On utilise les notations introduites en \ref{NotationsS0Induites}. 
Si le caractère $\overline{\chi}_1 \times \overline{\chi}_2$ est 
trivial sur $k_{\Delta}^{\times}$ alors $f$ est pair et $\nu_2 = \nu_1+\frac{f}{2}$.\\
Réciproquement, supposons que $\nu_2 = \nu_1+\frac{f}{2}$. 
Si $\overline{\chi}_1 \times \overline{\chi}_2$ est 
trivial sur $k_{\Delta}^{\times}$ alors $\overline{\chi}_0$ est trivial sur $k_{f/2}^{\times}$, 
où $k_{f/2}$ est une extension 
de degré $f/2$ de $k$ contenue dans $k_{\Delta}$. En particulier, $\overline{\chi}_0$ est trivial sur $k^{\times}$.\\
Ainsi, si $\varphi_{s_0}^{\nu_1, \nu_2} \neq 0$, le caractère $\overline{\chi}_0$ est trivial sur $k^{\times}$ 
et $f$ est pair (donc l'indice $\delta$ de la $\mathbb{K}$-algèbre à division $\Delta$ est pair).
\end{propo}

\begin{démo}
Supposons que $\overline{\chi}_1 \times \overline{\chi}_2$ soit trivial sur $k_{\Delta}^{\times}$, 
alors :
\begin{eqnarray*}
\forall x \in k_{\Delta}^{\times}, \overline{\chi}_0 (\Phi^{\nu_1} (x) \Phi^{\nu_2} (x)) = 1
& \Leftrightarrow & \forall x \in k_{\Delta}^{\times}, 
                   \overline{\chi}_0 (x \Phi^{\nu_2- \nu_1} (x)) = 1\\
& \Leftrightarrow & \forall x \in k_{\Delta}^{\times}, 
                   \overline{\chi}_0 \circ \Phi^{\nu_2- \nu_1} (x) = (\overline{\chi}_0 (x))^{-1}\\
& \Leftrightarrow & \forall x \in k_{\Delta}^{\times}, 
                       \overline{\chi}_0 \circ \Phi^{\nu_2- \nu_1} (x) = \overline{\chi}_0 (x^{-1})
\end{eqnarray*}
Soit $x$ dans $k_{\Delta}^{\times}$, alors :
$$
\overline{\chi}_0 (x) = \overline{\chi}_0 ((x^{-1})^{-1}) 
= \overline{\chi}_0 ( \Phi^{\nu_2- \nu_1} (x^{-1})) 
= \overline{\chi}_0 ( (\Phi^{\nu_2- \nu_1} (x))^{-1}) 
=   \overline{\chi}_0 \circ \Phi^{\nu_2- \nu_1} ( \Phi^{\nu_2- \nu_1} (x))
$$
On en déduit que 
$\overline{\chi}_0 = \overline{\chi}_0 \circ \Phi^{2 (\nu_2-\nu_1)}$. 
Ainsi $f$ divise $2 (\nu_2-\nu_1)$, or $2 \leq 2 (\nu_2-\nu_1) < 2f$, donc nécessairement 
$f = 2 (\nu_2-\nu_1)$, $f$ est pair et $\nu_2 = \nu_1+\frac{f}{2}$.\\
Réciproquement, supposons que $\nu_2 = \nu_1+\frac{f}{2}$. Alors 
$\overline{\chi}_1 \times \overline{\chi}_2$ est trivial sur $k_{\Delta}^{\times}$ si et seulement si :
$$
\forall x \in k_{\Delta}^{\times}, \overline{\chi}_0 (\Phi^{\nu_1} (x) \Phi^{\nu_1+f/2} (x)) = 1
 \Leftrightarrow  \forall x \in k_{\Delta}^{\times}, \overline{\chi}_0 (x \Phi^{f/2} (x)) = 1
$$
Notons $k_f$ et $k_{f/2}$ les extensions de corps de $k$ de degré $f$ et $f/2$ respectivement,
 contenues dans $k_{\Delta}$. Alors $k_{f/2}$ est l'ensemble des points fixes de $\Phi^{f/2}$. 
 Par restriction à $k_f$, on a $\langle \Phi^{f/2} \rangle = {\rm Gal} (k_f / k_{f/2})$ et :
 $$
 {\rm N}_{k_f / k_{f/2}} : k_f^{\times} \twoheadrightarrow k_{f/2}, x \mapsto x \Phi^{f/2} (x)
 $$
On en déduit que si $\overline{\chi}_1 \times \overline{\chi}_2$ est trivial sur $k_{\Delta}^{\times}$ alors :
$$
\forall x \in k_{f}^{\times}, \overline{\chi}_0 (x \Phi^{f/2} (x)) = 1 = \overline{\chi}_0 ( {\rm N}_{k_f / k_{f/2}} (x)) 
\Leftrightarrow \forall x \in k_{f/2}^{\times}, \overline{\chi}_0 (x) = 1
$$

\end{démo}

\subsubsection{Conditions sur $\varphi_{s_1}$.}\label{EtudeCNSommetS1}

\begin{nota}\label{NotationsS1}
On se place sous les hypothèses de \ref{NotationHypotheseCNDistinction}. 
Rappelons que $s_1$ est un sommet de $X_{\mathbb{K}}$ tel que $d(s_0, s_1) = 1$ où 
$s_0 = j(X_{\mathbb{F}})$ et :
$$
V_{s_1} = \left( \bigoplus_{\nu = 0}^{f-1} W_{s_1}^{\nu} \right) 
\oplus \left( \bigoplus_{0 \leq \nu_1 < \nu_2 \leq f-1} 
U_{s_1}^{\nu_1, \nu_2} \right)
$$
où :
$$
W_{s_1}^{\nu} \simeq \overline{\chi_0}^{\Phi^{\nu}} \otimes {\rm St}_{\overline{G}_{s_1}} 
\, \, \text{et} \, \, 
U_{s_1}^{\nu_1, \nu_2} 
\simeq {\rm Ind}_{\overline{B}_{s_1}}^{\overline{G}_{s_1}} 
\overline{\chi_0}^{\Phi^{\nu_1}} \otimes \overline{\chi_0}^{\Phi^{\nu_2}}
$$
Comme précédemment, on identifie  
$\mathbb{P}_{s_1}^1 (k_{\Delta})$ avec l'ensemble des droites correspondant aux arêtes de l'arbre 
$X_{\mathbb{K}}$ dont l'un des sommets est $s_1$. 
On notera $\delta_{Q+1}^{s_1}$ la droite correspondant à l'arête $\lbrace s_0, s_1 \rbrace$. 
Soient $u_1, \cdots, u_Q$ les sommets voisins de $s_1$ distincts de $s_0$ et 
$\delta_1^{s_1}, \cdots, \delta_Q^{s_1}$ les droites correspondantes. 
\end{nota}

\begin{lem}\label{LemmeFCourtes}
Pour tout $i$ dans 
$\lbrace 1, \cdots, Q \rbrace$, il existe $g_i$ dans $\mathcal{O}_{\mathcal{D}}^{\times}$ 
tel que $g_i.s_0 = s_0$, $g_i.s_1 = s_1$ et $g_i.u_i = u_1$. 
\end{lem}

\begin{démo}
On a démontré en \ref{OrbitesCasTotRamifie} que les sommets $u_1, \cdots, u_Q$ sont 
tous dans la même $\mathcal{O}_{\mathcal{D}}^{\times}$-orbite. 
Pour tout $i$ dans 
$\lbrace 1, \cdots, Q \rbrace$, on peut donc fixer $g_i$ dans $\mathcal{O}_{\mathcal{D}}^{\times}$ 
tel que :
$$
g_i.u_i = u_1
$$
Comme $g_i \in \mathcal{O}_{\mathcal{D}}^{\times} \subset \mathcal{D}^{\times}$, on a 
$g_i.s_0 = s_0$. Le groupe $\mathcal{D}^{\times}$ agit via des isométries sur les sommets. 
On en déduit facilement que $g_i.s_1$ est un sommet voisin de $s_0$. De même, on remarque que 
$g_i.s_1$ est un sommet voisin de $u_1 \neq s_0$. On a donc deux chemins géodésiques reliant 
$s_0$ à $u_1$ : 
$$
[s_0, s_1, u_1] \, \, \text{et} \, \, 
[s_0, g_i.s_1, u_1]
$$
Par unicité d'un tel chemin, on a bien 
$g_i.s_1 = s_1$.

\end{démo}

\begin{nota}
 Comme $u_1, \cdots, u_Q$ sont dans la même $\mathcal{O}_{\mathcal{D}}^{\times}$-orbite et, 
par lemme \ref{LemmeFCourtes}, pour tout $i$ dans 
$\lbrace 1, \cdots, Q \rbrace$, il existe $g_i$ dans $\mathcal{O}_{\mathcal{D}}^{\times}$ 
tel que $g_i.s_0 = s_0$, $g_i.s_1 = s_1$ et $g_i.\delta_i^{s_1} = \delta_1^{s_1}$. 
\end{nota}

\begin{lem}\label{CalculsReductionGiSommetS1}
Notons $\widetilde{g}_i$ la réduction de $g_i$ dans $\mathcal{U}_{s_1} / \mathcal{U}_{s_1}^1$. 
Alors, pour tout $i$ dans $\{ 1, \cdots, Q \}$, ${\rm det} (\widetilde{g}_i)$~est un carré dans 
$k^{\times}$, par conséquent, pour tout $\nu$ dans $\{ 0, \cdots, f-1 \}$ :
$$
\overline{\chi}_0^{\Phi^{\nu}} ({\rm det} (\widetilde{g}_i)) {\rm St}_{\overline{G}_{s_1}} (\widetilde{g}_i) 
= {\rm St}_{\overline{G}_{s_1}} (\widetilde{g}_i)
$$
\end{lem}

\begin{démo}
Soit $i$ dans $\lbrace 1, \cdots, Q \rbrace$. 
On a $g_i \in \mathcal{O}_{\mathcal{D}}^{\times} \subseteq \mathcal{U}_{s_0}$.  
Notons $\overline{g_i}$ la réduction de $g_i$ dans 
\mbox{$k_{\mathcal{D}}^{\times} \subseteq \mathcal{U}_{s_0} / \mathcal{U}_{s_0}^1$}. 
Alors $\overline{g}_i$ fixe la droite $\delta_{Q+1}^{s_1} \in \mathbb{P}_{s_0}^1 (k_{\Delta})$ 
(car fixe $s_0$ et $s_1$). Puisque 
$k_{\mathcal{D}}^{\times} / k_{\Delta}^{\times}$ agit simplement transitivement sur 
$\mathbb{P}_{s_0}^1 (k_{\Delta})$, on a 
$\overline{g}_i \in k_{\Delta}^{\times}$. 
Par suite :
$$
g_i \in (\mathcal{O}_{\Delta}^{\times} \mathcal{U}_{\mathcal{D}}^1 ) \cap \mathcal{U}_{s_1} 
\subseteq (\mathcal{O}_{\Delta}^{\times} \mathcal{U}_{s_0}^1 ) \cap \mathcal{U}_{s_1}
$$
Nous allons déterminer l'intersection $(\mathcal{O}_{\Delta}^{\times} \mathcal{U}_{s_0}^1 ) \cap \mathcal{U}_{s_1}$ 
puis sa réduction dans $\mathcal{U}_{s_1} / \mathcal{U}_{s_1}^1$. 
On identifie $\mathcal{U}_{s_0}$ à ${\rm GL}_2 (\mathcal{O}_{\Delta})$. 
Notons :
$$
w_0 = \left( \begin{array}{cc} 
                 1 & 0 \\
                 0 & \varpi_{\Delta} 
         \end{array} \right)
$$
Alors $\mathcal{U}_{s_1} = w_0 \mathcal{U}_{s_0} w_0^{-1}$ donc :
$$
\mathcal{U}_{s_0}^1 = \left( \begin{array}{cc} 
                1+ \mathcal{P}_{\Delta} & \mathcal{P}_{\Delta} \\
                 \mathcal{P}_{\Delta} & 1+\mathcal{P}_{\Delta} 
         \end{array} \right) \, \, \text{et} \, \, 
\mathcal{O}_{\Delta}^{\times} \mathcal{U}_{s_0} ^1
= \left\lbrace \left( \begin{array}{cc} 
                 u & \varpi_{\Delta} v \\
                 \varpi_{\Delta} w & u
         \end{array} \right) : u \in \mathcal{O}_{\Delta}^{\times}, v, w \in \mathcal{O}_{\Delta}  \right\rbrace
$$
de plus :
$$
\mathcal{U}_{s_1} = \left( \begin{array}{cc} 
                 \mathcal{O}_{\Delta} & \mathcal{P}_{\Delta}^{-1} \\
                 \mathcal{P}_{\Delta} & \mathcal{O}_{\Delta} 
         \end{array} \right) \, \, \text{d'où} \, \, 
(\mathcal{O}_{\Delta}^{\times} \mathcal{U}_{s_0} ^1) \cap \mathcal{U}_{s_1} 
=  \left\lbrace \left( \begin{array}{cc} 
                 u & \varpi_{\Delta} v \\
                 \varpi_{\Delta} w & u
         \end{array} \right) : u \in \mathcal{O}_{\Delta}^{\times}, v, w \in \mathcal{O}_{\Delta}  \right\rbrace
$$
On a :
$$
\mathcal{U}_{s_1}^1 = w_0 \mathcal{U}_{s_0}^1 w_0^{-1} 
= \left( \begin{array}{cc} 
                 1+\mathcal{P}_{\Delta} & \mathcal{O}_{\Delta} \\
                 \mathcal{P}_{\Delta}^2 & 1+\mathcal{P}_{\Delta} 
         \end{array} \right) \, \, \text{et} \, \, 
\mathcal{P}_{\Delta} / \mathcal{P}_{\Delta}^2 \simeq (k_{\Delta}, +)
$$
Soient $u_i$ dans $\mathcal{O}_{\Delta}^{\times}$, 
et $v_i, w_i$ dans $\mathcal{O}_{\Delta}$ tels que :
$$
g_i = \left( \begin{array}{cc} 
                 u_i & \varpi_{\Delta} v_i \\
                 \varpi_{\Delta} w_i & u_i
         \end{array} \right)
$$
alors :
$$
\widetilde{g}_i = \left( \begin{array}{cc} 
                 \overline{u}_i & 0 \\
                 \overline{w}_i & \overline{u}_i
         \end{array} \right)
= \overline{u}_i \left( \begin{array}{cc} 
                 1 & 0 \\
                 \overline{u}_i^{-1} \overline{w}_i & 1
         \end{array} \right)
$$
où $\overline{u}_i$ et $\overline{w}_i$ sont les réductions de $u_i$ et $w_i$ dans 
$k_{\Delta}^{\times}$. Ainsi :
$$
{\rm det} (\widetilde{g}_i) = {\rm det} \left( \overline{u}_i \left( \begin{array}{cc} 
                 1 & 0 \\
                 \overline{u}_i^{-1} \overline{w}_i & 1
         \end{array} \right) \right) 
= {\rm det} (\overline{u}_i) \times 1
$$
Or ${\rm det} (\overline{u}_i) = {\rm N}_{k_{\mathcal{D}} / k} (\overline{u}_i)$ 
est un carré dans $k^{\times}$, donc 
$\overline{\chi}_0^{\Phi^{\nu}} ({\rm det} (\widetilde{g}_i)) = 1$.

\end{démo}

\paragraph{Etude de $\varphi_{s_1}^{\nu}$ pour $0 \leq \nu \leq f-1$.}

\begin{nota}\label{NotationsS1Steinberg}
On fixe $\nu$ dans $\lbrace 0, \cdots, f-1 \rbrace$, et 
$\widetilde{\varphi}_{s_1}^{\nu}$ dans 
$\lbrace h : \mathbb{P}_{s_1}^1 (k_{\Delta}) \rightarrow \mathbb{C} \rbrace / 
\lbrace \text{fonctions constantes} \rbrace$ telle que pour tout $h$ dans $W_{s_1}^{\nu}$ :
$$
\varphi_{s_1}^{\nu} (h) 
= \sum_{d \in \mathbb{P}_{s_1}^{1} (k_{\Delta})} h(d) \widetilde{\varphi}_{s_1}^{\nu} (d)
$$
Notons $\overline{\chi}_1 = \overline{\chi}_0^{\nu}$.
\end{nota}

\begin{propo}
On se place sous les hypothèses de \ref{NotationsS1} et 
\ref{NotationsS1Steinberg}. 
Supposons que $\overline{\chi}_0$ est trivial sur $k^{\times}$. 
Alors $\varphi_{s_1}^{\nu} = 0$.
\end{propo}

\begin{démo}
Supposons que $\overline{\chi}_0$ est trivial sur $k^{\times}$. 
Dans ce cas, nous avons déjà vu que $\varphi_{s_0}^{\nu} = 0$. 
Fixons $h$ dans $W_{s_1}^{\nu}$ telle que $h(\delta_{Q+1}^{s_1}) = 1$ et 
$h(\delta_i^{s_1}) = - \frac{1}{Q}$ si $i \neq Q+1$. Alors 
$h$ est dans le module de Jacquet $(W_{s_1}^{\nu})^{\mathcal{U}_a^1}$ correspondant 
à l'arête $a = \{ s_0, s_1 \}$. Les formes linéaires 
$\varphi_{s_0}^{\nu}$ et $\varphi_{s_1}^{\nu}$ coïncident sur $(W_{s_1}^{\nu})^{\mathcal{U}_a^1}$, donc  
$\varphi_{s_1}^{\nu} (h) = \varphi_{s_0}^{\nu} (h) = 0$. 
Or, $\varphi_{s_1}^{\nu} (h) = \widetilde{\varphi}_{s_1}^{\nu} (\delta_{Q+1}^{s_1}) 
- \frac{1}{Q} \sum_{i=1}^Q \widetilde{\varphi}_{s_1}^{\nu} (\delta_i^{s_1})$. 
On en déduit que  
\mbox{$\widetilde{\varphi}_{s_1}^{\nu} (\delta_{Q+1}^{s_1}) 
= \frac{1}{Q} \sum_{i=1}^Q \widetilde{\varphi}_{s_1}^{\nu} (\delta_i^{s_1})$}. 
Puis, en remarquant que pour tout $i$ dans $\lbrace 1, \cdots, Q \rbrace$, on a 
\mbox{$\varphi_{s_1}^{\nu} = \varphi_{s_1}^{\nu} \circ \pi (g_i)$}, et avec des calculs analogues à l'étude 
de $\varphi_{s_0}^{\nu}$, on montre que 
$\widetilde{\varphi}_{s_1}^{\nu} (\delta_{1}^{s_1}) = \cdots 
= \widetilde{\varphi}_{s_1}^{\nu} (\delta_{Q}^{s_1})$. 
Par conséquent, $\widetilde{\varphi}_{s_1}^{\nu}$ est constante et 
$\varphi_{s_1}^{\nu} = 0$.

\end{démo}

\begin{propo}
On utilise les notations de \ref{NotationsS1} et 
\ref{NotationsS1Steinberg}. 
Supposons que $\overline{\chi}_0$ n'est pas trivial sur $k^{\times}$ (mais uniquement trivial 
sur les carrés de $k^{\times}$).\\
Si $\varphi_{s_0}^{\nu} = 0$, alors $\varphi_{s_1}^{\nu} = 0$.\\
Sinon, pour tout $h$ dans $W_{s_1}^{\nu}$ :
$$
\varphi_{s_1}^{\nu} (h) 
= \frac{Q+1}{2Q} 
\times (\widetilde{\varphi}_{s_0}^{\nu} (\delta_1^{s_0}) 
- \widetilde{\varphi}_{s_0}^{\nu} (\delta_{Q+1}^{s_0})) \times h(\delta_{Q+1}^{s_1})
$$
et : 
$$
{\rm ker} ( \varphi_{s_1}^{\nu}) = 
\lbrace h \in W_{s_1}^{\nu} : h(\delta_{Q+1}^{s_1}) = 0 \rbrace
$$
\end{propo}

\begin{démo}
Supposons que $\overline{\chi}_0$ n'est pas trivial sur $k^{\times}$. 
Il nous faut distinguer deux cas. 
Tout d'abord, si $\widetilde{\varphi}_{s_0}^{\nu} (\delta_1^{s_0}) 
= \widetilde{\varphi}_{s_0}^{\nu} (\delta_{Q+1}^{s_0})$, alors 
$\varphi_{s_0}^{\nu} = 0$ et, avec des calculs analogues au cas précédent, on montre 
que $\varphi_{s_1}^{\nu} = 0$. 
Sinon, rappelons que pour tout $h$ dans $W_{s_0}^{\nu}$ :
$$ 
\varphi_{s_0}^{\nu} (h) = \left( \sum_{i=1}^{(Q+1)/2} h(\delta_i^{s_0}) \right) 
\times (\widetilde{\varphi}_{s_0}^{\nu} (\delta_1^{s_0}) 
- \widetilde{\varphi}_{s_0}^{\nu} (\delta_{Q+1}^{s_0}))
$$
Les applications $\varphi_{s_0}^{\nu}$ et $\varphi_{s_1}^{\nu}$ coïncident sur le module de Jacquet 
$(W_{s_1}^{\nu})^{\delta_{Q+1}^{s_1}}$ (les fonctions dans $W_{s_1}^{\nu}$ qui sont 
constantes sur $\mathbb{P}_{s_1}^1 (k_{\Delta}) \backslash \{ \delta_{Q+1}^{s_1} \}$). 
Fixons $h$ dans $(W_{s_1}^{\nu})^{\delta_{Q+1}^{s_1}}$ telle que $h(\delta_{Q+1}^{s_1}) = 1$ et 
$h(\delta) = - \frac{1}{Q}$ sinon. Alors 
$\varphi_{s_1}^{\nu} (h) = \varphi_{s_0}^{\nu} (h)$, ainsi :
$$
\widetilde{\varphi}_{s_1}^{\nu} (\delta_{Q+1}^{s_1}) - \frac{1}{Q} 
\sum_{i=1}^Q \widetilde{\varphi}_{s_1}^{\nu} (\delta_{i}^{s_1}) 
= \left( \sum_{i=1}^{(Q+1)/2} h(\delta_i^{s_0}) \right) 
\times (\widetilde{\varphi}_{s_0}^{\nu} (\delta_1^{s_0}) 
- \widetilde{\varphi}_{s_0}^{\nu} (\delta_{Q+1}^{s_0}))
$$
Rappelons que la droite $\delta_{Q+1}^{s_1}$ correspond à l'arête 
$\{ s_0, s_1 \}$ et donc à la droite $\delta_1^{s_0}$.  
Ainsi :
$$
\sum_{i=1}^{(Q+1)/2} h(\delta_i^{s_0}) = 1 - \frac{Q-1}{2} \times \frac{1}{Q} 
= \frac{Q+1}{2Q}
$$
Par suite 
$\widetilde{\varphi}_{s_1}^{\nu} (\delta_{Q+1}^{s_1}) - \frac{1}{Q} 
\sum_{i=1}^Q \widetilde{\varphi}_{s_1}^{\nu} (\delta_{i}^{s_1}) 
= \frac{Q+1}{2Q} 
\times (\widetilde{\varphi}_{s_0}^{\nu} (\delta_1^{s_0}) 
- \widetilde{\varphi}_{s_0}^{\nu} (\delta_{Q+1}^{s_0}))$. 
Comme dans le cas précédent, on montre que 
$\widetilde{\varphi}_{s_1}^{\nu} (\delta_{1}^{s_1})  = \cdots 
= \widetilde{\varphi}_{s_1}^{\nu} (\delta_{Q}^{s_1})$. 
On en déduit que 
\mbox{$\widetilde{\varphi}_{s_1}^{\nu} (\delta_{Q+1}^{s_1}) 
= \widetilde{\varphi}_{s_1}^{\nu} (\delta_{1}^{s_1}) 
+ \frac{Q+1}{2Q} 
\times (\widetilde{\varphi}_{s_0}^{\nu} (\delta_1^{s_0}) 
- \widetilde{\varphi}_{s_0}^{\nu} (\delta_{Q+1}^{s_0}))$}. 
Ainsi, pour tout $h$ dans $W_{s_1}^{\nu}$, on a :
\begin{eqnarray*}
\varphi_{s_1}^{\nu} (h) 
& = & \widetilde{\varphi}_{s_1}^{\nu} (\delta_{Q+1}^{s_1}) h(\delta_{Q+1}^{s_1}) 
+ \widetilde{\varphi}_{s_1}^{\nu} (\delta_{1}^{s_1}) \times \left( \sum_{i=1}^Q 
h (\delta_i^{s_1}) \right)\\
& = & \widetilde{\varphi}_{s_1}^{\nu} (\delta_1^{s_1}) \times h(\delta_{Q+1}^{s_1})
+ \frac{Q+1}{2Q} 
\times (\widetilde{\varphi}_{s_0}^{\nu} (\delta_1^{s_0}) 
- \widetilde{\varphi}_{s_0}^{\nu} (\delta_{Q+1}^{s_0})) \times h(\delta_{Q+1}^{s_1}) 
+ \widetilde{\varphi}_{s_1}^{\nu} (\delta_{1}^{s_1}) \times \left( \sum_{i=1}^Q 
h (\delta_i^{s_1}) \right)\\
& = & \widetilde{\varphi}_{s_1}^{\nu} (\delta_{1}^{s_1}) \times \left( \sum_{i=1}^{Q+1} 
h (\delta_i^{s_1}) \right) + \frac{Q+1}{2Q} 
\times (\widetilde{\varphi}_{s_0}^{\nu} (\delta_1^{s_0}) 
- \widetilde{\varphi}_{s_0}^{\nu} (\delta_{Q+1}^{s_0})) \times h(\delta_{Q+1}^{s_1}) \\
& = & \frac{Q+1}{2Q} 
\times (\widetilde{\varphi}_{s_0}^{\nu} (\delta_1^{s_0}) 
- \widetilde{\varphi}_{s_0}^{\nu} (\delta_{Q+1}^{s_0})) \times h(\delta_{Q+1}^{s_1})
\end{eqnarray*}

\end{démo}

\paragraph{Etude de $\varphi_{s_1}^{\nu_1, \nu_2}$ pour $0 \leq \nu_1 < \nu_2 \leq f-1$.}

\begin{nota}\label{NotationsS1Induites}
On se place sous les hypothèses de \ref{NotationsS1}. 
On fixe deux entiers naturels $\nu_1$ et $\nu_2$ tels que 
$0 \leq \nu_1 < \nu_2 \leq f-1$. 
Nous noterons 
$\overline{\chi}_1 = \overline{\chi}_0^{\Phi^{\nu_1}}$, 
$\overline{\chi}_2 = \overline{\chi}_0^{\Phi^{\nu_2}}$. 
Soit \mbox{$\pi_{s_1}^{\nu_1, \nu_2} = U_{s_1}^{\nu_1, \nu_2} 
= {\rm Ind}_{\overline{B}_{s_1}}^{\overline{G}_{s_1}} 
(\overline{\chi}_1 \otimes \overline{\chi}_2 )$}. 
On utilise le modèle présenté en \ref{ModeleInduiteParabolique}. 
On fixe $\widetilde{\varphi}_{s_1}^{\nu_1, \nu_2}$ dans 
${\rm Ind}_{\overline{B}_{s_1}}^{\overline{G}_{s_1}} 
(\overline{\chi}_1^{-1} \otimes \overline{\chi}_2^{-1} )$ tel que pour tout $h$ 
dans~$\pi_{s_1}^{\nu_1, \nu_2}$ :
$$ 
\varphi_{s_1}^{\nu_1, \nu_2} (h) 
= \sum_{x \in \overline{G}_{s_1}} \widetilde{\varphi}_{s_1}^{\nu_1, \nu_2} (x) h(x)
$$
\end{nota}

\begin{lem}
On a :
$$
\varphi_{s_1}^{\nu_1, \nu_2} \in {\rm Hom}_{K_{s_1}} (\pi_{s_1}^{\nu_1, \nu_2}, \mathds{1})
$$
où $K_{s_1}$ est la réduction de $\mathcal{O}_{\mathcal{D}}^{\times} \cap \mathcal{U}_{s_1}$ 
dans $\overline{G}_{s_1}$. Alors :
$$
\langle \chi_{\pi_{s_1}^{\nu_1, \nu_2}}, \mathds{1}_{K_{s_1}} \rangle 
= \frac{2}{Q-1} \times (\sum_{u \in k_{\Delta}^{\times}} \overline{\chi}_1 \times \overline{\chi}_2 (u))
$$
\end{lem}

\begin{démo}
Soit $g$ dans $\mathcal{O}_{\mathcal{D}}^{\times} \cap \mathcal{U}_{s_1}$ et $\widetilde{g}$ la 
réduction de $g$ dans $\overline{G}_{s_1} = \mathcal{U}_{s_1} / \mathcal{U}_{s_1}^1$. On a :
$$
\varphi_{s_1}^{\nu_1, \nu_2} = \varphi_{s_1}^{\nu_1, \nu_2} \circ \pi_{s_1}^{\nu_1, \nu_2} (\widetilde{g})
$$
Notons $K_{s_1}$ la réduction de $\mathcal{O}_{\mathcal{D}}^{\times} \cap \mathcal{U}_{s_1}$ 
dans $\overline{G}_{s_1}$. D'après les calculs précédents, on a :
$$
K_{s_1} = \left\lbrace u \left( \begin{array}{cc}
                               1 & 0\\
                               w & 1
                               \end{array} \right) : 
u \in k_{\Delta}^{\times}, w \in k_{\Delta} \right\rbrace
$$
et 
$\varphi_{s_1}^{\nu_1, \nu_2} \in {\rm Hom}_{K_{s_1}} (\pi_{s_1}^{\nu_1, \nu_2}, \mathds{1})$. 
L'espace d'entrelacement ${\rm Hom}_{K_{s_1}} (\pi_{s_1}^{\nu_1, \nu_2}, \mathds{1})$ 
est non trivial si et seulement si 
$\langle \chi_{\pi_{s_1}^{\nu_1, \nu_2}}, \mathds{1}_{K_{s_1}} \rangle \neq 0$, or :
$$
\langle \chi_{\pi_{s_1}^{\nu_1, \nu_2}}, \mathds{1}_{K_{s_1}} \rangle 
= \frac{1}{\vert K_{s_1} \vert \times \vert \overline{B}_{s_1} \vert} 
\sum_{g \in K_{s_1}} 
\sum_{\substack{x \in \overline{G}_{s_1} \\ x^{-1} g x \in \overline{B}_{s_1}}} 
\overline{\chi}_1 \otimes \overline{\chi}_{2} (x^{-1} g x)
$$
On montre facilement que si $x \in \overline{G}_{s_1}$ et $g \in K_{s_1}$, alors 
$x^{-1} g x \in \overline{B}_{s_1}$ si et seulement si $x \in \overline{B}_{s_1}$. 
On peut à présent calculer 
$\langle \chi_{\pi_{s_1}^{\nu_1, \nu_2}}, \mathds{1}_{K_{s_1}} \rangle $.
\begin{eqnarray*}
\langle \chi_{\pi_{s_1}^{\nu_1, \nu_2}}, \mathds{1}_{K_{s_1}} \rangle 
& = & \frac{1}{\vert K_{s_1} \vert \times \vert \overline{B}_{s_1} \vert} 
\sum_{g \in K_{s_1}} 
\sum_{\substack{x \in \overline{G}_{s_1} \\ x^{-1} g x \in \overline{B}_{s_1}}} 
\overline{\chi}_1 \otimes \overline{\chi}_{2} (x^{-1} g x)\\
& = & \frac{\vert \overline{G}_{s_1} \vert}{\vert K_{s_1} \vert \times \vert \overline{B}_{s_1} \vert} 
\sum_{g \in k_{\Delta}^{\times}}  
\overline{\chi}_1 \times \overline{\chi}_{2} (g)\\
& + & \frac{1}{\vert K_{s_1} \vert \times \vert \overline{B}_{s_1} \vert} 
\sum_{u \in k_{\Delta}^{\times}} \sum_{w \in k_{\Delta}^{\times}} \sum_{x \in \overline{B}_{s_1}} 
\overline{\chi}_1 \otimes \overline{\chi}_{2} \left(u \times \left( x^{-1} \left( \begin{array}{cc}
                               1 & 0\\
                               w & 1
                               \end{array} \right) x \right)\right)\\
& = & \frac{\vert \overline{G}_{s_1} \vert}{\vert K_{s_1} \vert \times \vert \overline{B}_{s_1} \vert} 
\sum_{g \in k_{\Delta}^{\times}}  
\overline{\chi}_1 \times \overline{\chi}_{2} (g)
  +  \frac{\vert k_{\Delta}^{\times} \vert \times \vert \overline{B}_{s_1} \vert}{\vert K_{s_1} \vert \times \vert \overline{B}_{s_1} \vert} 
\sum_{u \in k_{\Delta}^{\times}}  
\overline{\chi}_1 \times \overline{\chi}_{2} (u)\\
& = & \frac{2}{Q-1} \times \left( \sum_{u \in k_{\Delta}^{\times}}  
\overline{\chi}_1 \times \overline{\chi}_{2} (u) \right)
\end{eqnarray*}
\end{démo}

On en déduit immédiatement les résultats suivants :

\begin{propo}
On utilise les notations de \ref{NotationsS1Induites}. 
Si $\overline{\chi}_1 \times \overline{\chi}_2$ n'est pas trivial sur $k_{\Delta}^{\times}$, 
on a $\varphi_{s_1}^{\nu_1, \nu_2} = 0$.
\end{propo}

\begin{propo}
Si $\overline{\chi}_1 \times \overline{\chi}_2$ est trivial sur $k_{\Delta}^{\times}$ 
(c'est le cas si $\overline{\chi}_0$ est trivial sur $k_{f/2}^{\times}$, si $f$ est pair 
et $\nu_2 = \nu_1+ f/2$), alors :
$$
\langle \chi_{\pi_{s_1}^{\nu_1, \nu_2}}, \mathds{1}_{K_{s_1}} \rangle  
= 2
$$
et $\varphi_{s_1}^{\nu_1, \nu_2} \in 
{\rm Hom}_{K_{s_1}} (\pi_{s_1}^{\nu_1, \nu_2}, \mathds{1}) $ où 
${\rm Hom}_{K_{s_1}} (\pi_{s_1}^{\nu_1, \nu_2}, \mathds{1})$ est un espace de dimension $2$. 
\end{propo}

\subsubsection{Conditions sur $\varphi_{s_k}$ pour $k \geq 1$.}\label{EtudeCNSommetSk}

\begin{nota}
On se place sous les hypothèses de \ref{NotationHypotheseCNDistinction}. 
Rappelons que $s_0 = j(X_{\mathbb{F}})$. 
On fixe $k$ un entier naturel supérieur ou égal à $1$. 
Alors $s_1, \cdots, s_k$ sont des sommets de $X_{\mathbb{K}}$ tels que 
$d(s_j, s_{j+1}) = 1$ pour tout $j$ dans $\{ 0, \cdots, k-1 \}$. 
On a $\overline{G}_{s_k} = \mathcal{U}_{s_k} / \mathcal{U}_{s_k}^1$. 
Soit $\mathbb{P}_{s_k}^1 (k_{\Delta}) = \{ \delta_1^{s_k}, \cdots, \delta_{Q+1}^{s_k} \}$ les droites 
correspondants aux arêtes de $X_{\mathbb{K}}$ dont $s_k$ est l'un des sommets où 
$\delta_{Q+1}^{s_k}$ correspond à l'arête $\{ s_{k-1}, s_k \}$. 
Notons $u_1, \cdots, u_Q$ les sommets voisins de $s_k$ distincts de $s_{k-1}$. 
\end{nota}

\begin{lem}\label{LemmeFCourtesVersion2}
Pour tout $i$ dans $\{ 1, \cdots, Q \}$, il existe 
$g_i$ dans $\mathcal{O}_{\mathcal{D}}^{\times} \cap \mathcal{U}_{s_1} \cap \cdots \cap \mathcal{U}_{s_k}$ tel que 
$g_i. u_i = u_1$.
\end{lem}

\begin{démo}
On reprend le même raisonnement que dans la démonstration de \ref{LemmeFCourtes}. 
Pour $i$ dans $\{ 1, \cdots, Q \}$, il existe $g_i$ dans $\mathcal{O}_{\mathcal{D}}^{\times}$ 
tel que $g_i.u_i = u_1$. On a alors deux chemins géodésiques reliant $s_0$ et~$u_1$ :
$$
[s_0, s_1, \cdots, s_k, u_1 ] \, \, \text{et} \, \, 
[s_0, g_i.s_1, \cdots, g_i.s_k, u_1 ]
$$
Par unicité d'un tel chemin, on a bien 
$g_i.s_1 = s_1, g_i.s_2 = s_2, \cdots, g_i.s_k = s_k$.

\end{démo}

\begin{nota}
D'après le lemme \ref{LemmeFCourtesVersion2}, on sait que 
$\mathcal{O}_{\mathcal{D}}^{\times} \cap \mathcal{U}_{s_1} \cap \cdots \cap \mathcal{U}_{s_k}$ agit transitivement 
sur les sommets $u_1, \cdots, u_Q$ et donc, pour tout $i$ dans $\{ 1, \cdots, Q \}$, il existe 
$g_i$ dans $\mathcal{O}_{\mathcal{D}}^{\times} \cap \mathcal{U}_{s_1} \cap \cdots \cap \mathcal{U}_{s_k}$ tel que 
$g_i. \delta_i^{s_k} = \delta_1^{s_k}$.
\end{nota}

Avec des calculs analogues à ceux de la démonstration de \ref{CalculsReductionGiSommetS1}, on montre 
le résultat suivant :

\begin{lem}
Notons $\widetilde{g}_i$ la réduction de $g_i$ dans $\mathcal{U}_{s_k} / \mathcal{U}_{s_k}^1$. 
Alors :
$$
\overline{\chi}_0^{\Phi^{\nu}} ({\rm det} (\widetilde{g}_i)) = 1
$$
pour tout $\nu$ dans $\{ 0, \cdots, f-1 \}$.
\end{lem}

On montre facilement par récurrence sur $k$ les propriétés suivantes :

\begin{propo}
Soit $\nu$ un entier naturel tel que $0 \leq \nu \leq f-1$, alors :

\begin{itemize}
\item[i)] Si le caractère $\overline{\chi}_0$ est trivial sur $k^{\times}$, 
$\varphi_{s_k}^{\nu} = 0$.
\item[ii)] Si le caractère $\overline{\chi}_0$ n'est pas trivial sur $k^{\times}$, on distingue deux cas :
ou bien $\varphi_{s_0}^{\nu} = 0$ et alors $\varphi_{s_k}^{\nu} = 0$; ou bien 
$\varphi_{s_0}^{\nu} \neq 0$, et alors :
$$
\varphi_{s_k}^{\nu} (h) 
= \frac{(-1)^{k-1} (Q+1)}{2 Q^{k}} \times 
(\widetilde{\varphi}_{s_0}^{\nu} (\delta_1^{s_0}) 
- \widetilde{\varphi}_{s_0}^{\nu} (\delta_{Q+1}^{s_0})) \times h(\delta_{Q+1}^{s_k})
$$
pour tout $h$ dans $W_{s_k}^{\nu}$.
\end{itemize}

\end{propo}

\begin{propo}
Si $\overline{\chi}_1 \times \overline{\chi}_2$ n'est pas trivial sur 
$k_{\Delta}^{\times}$, on a 
$\varphi_{s_k}^{\nu_1, \nu_2} = 0$.
\end{propo}

\begin{démo}
Soit $g$ dans $\mathcal{O}_{\mathcal{D}}^{\times} \cap \mathcal{U}_{s_1} 
\cap \cdots \cap \mathcal{U}_{s_k}$ et $\widetilde{g}$ la 
réduction de $g$ dans $\overline{G}_{s_k} = \mathcal{U}_{s_k} / \mathcal{U}_{s_k}^1$. On a 
$\varphi_{s_k}^{\nu_1, \nu_2} = \varphi_{s_k}^{\nu_1, \nu_2} \circ \pi_{s_k}^{\nu_1, \nu_2} (\widetilde{g})$. 
Notons $K_{s_k}$ la réduction de $\mathcal{O}_{\mathcal{D}}^{\times} \cap \mathcal{U}_{s_1} \cap \cdots \cap \mathcal{U}_{s_k}$ 
dans $\overline{G}_{s_k}$, alors :
$$
K_{s_k} = \left\lbrace u \left( \begin{array}{cc}
                               1 & 0\\
                               w & 1
                               \end{array} \right) : 
u \in k_{\Delta}^{\times}, w \in k_{\Delta} \right\rbrace
$$
et 
$\varphi_{s_k}^{\nu_1, \nu_2} \in {\rm Hom}_{K_{s_k}} (\pi_{s_k}^{\nu_1, \nu_2}, \mathds{1})$. 
Or :
$$
\langle \chi_{\pi_{s_k}^{\nu_1, \nu_2}}, \mathds{1}_{K_{s_k}} \rangle 
= \frac{1}{\vert K_{s_k} \vert \times \vert \overline{B}_{s_k} \vert} 
\sum_{g \in K_{s_k}} 
\sum_{\substack{x \in \overline{G}_{s_k} \\ x^{-1} g x \in \overline{B}_{s_k}}} 
\overline{\chi}_1 \otimes \overline{\chi}_{2} (x^{-1} g x)
$$
Soit $g = u \left( \begin{array}{cc}
                               1 & 0\\
                               w & 1
                               \end{array} \right) \in K_{s_k}$ 
et $x \in \overline{G}_{s_k}$. Si 
$x^{-1} g x$ appartient à $\overline{B}_{s_k}$, on a :
$$
\overline{\chi}_1 \otimes \overline{\chi}_{2} (x^{-1} g x) 
= (\overline{\chi}_1 \times \overline{\chi}_2 (u)) \times 1
$$ 
On en déduit que $\langle \chi_{\pi_{s_k}^{\nu_1, \nu_2}}, \mathds{1}_{K_{s_k}} \rangle$ est 
un multiple de 
$\sum_{u \in k_{\Delta}^{\times}} \overline{\chi}_1 \times \overline{\chi}_2 (u) = 0$
(car $\overline{\chi}_1 \times \overline{\chi}_2$ n'est pas trivial sur $k_{\Delta}^{\times}$). 
Par conséquent 
$\varphi_{s_k}^{\nu_1, \nu_2} = 0$.

\end{démo}

\subsubsection{Première conclusion sur les conditions nécessaires de distinction.}

On utilise les notations de \ref{NotationHypotheseCNDistinction}, \ref{NotationsS0Steinberg}.

\begin{theo}\label{SyntheseCNDistinctionTotRam}
Soit $(\pi, V)$ une représentation membre de la série discrète de $G$, de niveau $0$, 
non cuspidale et $(\mathbb{K}_f, \chi_f)$ la paire admissible modérée associée à $\pi$ 
définie en \ref{ParamSilbergerZinkNonCuspidales}. On définit comme dans 
\ref{DefinitionSystCoeffShneiderEtStuhler} un système de coefficients sur $X_{\mathbb{K}}$ d'espace vectoriel $V$. 
Soit $\varphi = (\varphi_s)_{s \in X_0}$ dans 
\mbox{${\rm ker} (\partial_1^{\ast}) \cap {\rm Hom}_{\mathcal{D}^{\times}} (C_0, \mathds{1})$}.\\
Si $f$ est impair et si $\overline{\chi}_0$ est trivial sur $k^{\times}$, alors 
$\pi$ ne peut pas être $\mathcal{D}^{\times}$-distinguée.\\
Si $f$ est pair, si $\overline{\chi}_0$ est trivial sur $k^{\times}$ sans être trivial 
sur $k_{f/2}^{\times}$, alors 
$\pi$ ne peut pas être $\mathcal{D}^{\times}$-distinguée.\\
Supposons que $\pi$ soit $\mathcal{D}^{\times}$-distinguée. Alors :
\begin{itemize}
\item[i)] Si $\overline{\chi}_0$ est non trivial sur $k^{\times}$ mais trivial sur les carrés de 
$k^{\times}$, pour tout sommet $s$, pour tout $0 \leq \nu_1 < \nu_2 \leq f-1$, on a :
$$
\varphi_s^{\nu_1, \nu_2} = 0
$$
et pour tout $\nu \in \{ 0, \cdots, f-1 \}$, 
pour tout sommet $s$ distinct de $s_0$, on a :
$$
\varphi_{s}^{\nu} (h) 
= \frac{(-1)^{k-1} (Q+1)}{2 Q^{k}} \times 
(\widetilde{\varphi}_{s_0}^{\nu} (\delta_1^{s_0}) 
- \widetilde{\varphi}_{s_0}^{\nu} (\delta_{Q+1}^{s_0})) \times h(\delta_{Q+1}^{s}) \, \,  
\text{pour tout} \, \, h \in W_{s}^{\nu} 
$$
où $k = d (s, s_0)$, $\delta_{Q+1}^{s}$ est la droite correspondant à l'arête 
$\{ \widetilde{s}, s \}$, avec $[ s_0, \cdots, \widetilde{s}, s]$ chemin géodésique reliant $s$ 
à $s_0$. 
\item[ii)] Si $f$ est pair, si $\overline{\chi}_0$ est trivial sur $k_{f/2}^{\times}$ 
(donc en particulier, $\overline{\chi}_0$ est trivial sur $k^{\times}$), alors, pour tout sommet $s$ 
et tout $0 \leq \nu \leq f-1$ :
$$
\varphi_s^{\nu} = 0
$$
et $\varphi_s^{\nu_1, \nu_2} = 0$ dès que 
$\nu_2 \neq \nu_1 + \frac{f}{2}$.
\end{itemize}

\end{theo}

\subsubsection{Multiplicité $1$.}

On utilise encore dans cette partie les notations et hypothèses de \ref{NotationHypotheseCNDistinction}. 
On suppose en particulier que $\pi$ est $\mathcal{D}^{\times}$-distinguée. 
On fixe $\varpi_{\mathbb{K}}$ et $\varpi_{\mathbb{F}}$ telles que 
$\varpi_{\mathbb{K}}^2 = \varpi_{\mathbb{F}}$. 

\begin{propo}\label{ResultatMultiplicite1Steinberg}
Soit $Y_{{\rm St}}$ l'ensemble les formes linéaires 
$\varphi$ dans ${\rm ker} (\partial_1^{\ast})$ telles que  
pour tout sommet $s$, pour tout $0 \leq \nu_1 < \nu_2 \leq f-1$, on a $\varphi_s^{\nu_1, \nu_2} = 0$. 
L'application :
$$
Y_{{\rm St}} \rightarrow {\rm Hom}_{k_{\mathcal{D}}^{\times}} (\pi_{s_0}^0, \mathds{1}), 
\varphi \mapsto \varphi_{s_0}^0
$$  
est injective.
\end{propo}

\begin{démo}
Fixons $\varphi$ dans $Y_{{\rm St}}$. 
En utilisant les résultats de \ref{SyntheseCNDistinctionTotRam}, 
on remarque que, pour $\nu$ dans 
$\{ 0, \cdots, f-1 \}$, et pour tout sommet $s$, la forme linéaire 
$\varphi_s^{\nu}$ est entièrement déterminée par $\varphi_{s_0}^{\nu}$. Ainsi, l'application :
$$
(\varphi_s^{\nu})_{s \in X_0} \mapsto \varphi_{s_0}^{\nu}
$$
est injective. 
Rappelons que pour tout $s$ dans $X_0$ (i.e $s$ sommet de $X_{\mathbb{K}}$) et 
pour tout $x$ dans $\mathcal{D}^{\times}$, on a  
$\varphi_s = \varphi_{x.s} \circ \pi (x)$. 
En particulier, si $g \in \mathcal{O}_{\mathcal{D}}^{\times}$, alors 
$\varpi_{\mathcal{D}} g \varpi_{\mathcal{D}}^{-1} \in \mathcal{D}^{\times}$ et :
$$
\varphi_{s_0} = 
\varphi_{(\varpi_{\mathcal{D}} g \varpi_{\mathcal{D}}^{-1}).s_0} \circ 
\pi (\varpi_{\mathcal{D}} g \varpi_{\mathcal{D}}^{-1}) 
= \varphi_{s_0} \circ \pi (\varpi_{\mathcal{D}} g \varpi_{\mathcal{D}}^{-1})
$$
Que dire de l'opérateur $\pi (\varpi_{\mathcal{D}} g \varpi_{\mathcal{D}}^{-1})$ 
pour $g$ dans $\mathcal{O}_{\mathcal{D}}^{\times}$?\\
On remarque que 
$\varpi_{\mathcal{D}}^{2d} = \varpi_{\mathbb{F}} 
= \varpi_{\mathbb{K}}^2 = \varpi_{\Delta}^{2d}$. 
Ainsi $\varpi_{\mathcal{D}} \in \varpi_{\Delta} {\rm GL}_2 (\mathcal{O}_{\Delta})$. 
De plus, la conjugaison par $\varpi_{\Delta}$ engendre le groupe de Galois 
$\langle \Phi \rangle = {\rm Gal} (k_{\Delta} / k)$. 
Pour $g$ dans $\mathcal{U}_{s_0} \simeq {\rm GL}_2 (\mathcal{O}_{\Delta})$, on note $g^{\Phi} = \Phi (g)$ 
(action sur les coefficients). 
On notera $\widetilde{\pi}$ la restriction de $\pi$ à ${\rm GL}_2 (\mathcal{O}_{\Delta})$ et 
$\widetilde{\pi}^{\Phi} : 
g \in {\rm GL}_2 (\mathcal{O}_{\Delta}) \mapsto \widetilde{\pi} (g^{\Phi})$. 
Pour tout $g$ dans ${\rm GL}_2 (\mathcal{O}_{\Delta})$, on a :
$$
\widetilde{\pi} (\varpi_{\mathcal{D}} g \varpi_{\mathcal{D}}^{-1}) 
= \widetilde{\pi} (g^{\Phi})
$$
(car la conjugaison par un élément de $\mathcal{U}_{s_0}$ transforme $\widetilde{\pi}$ en 
une représentation équivalente). 
Pour tout $g$ dans $\mathcal{O}_{\mathcal{D}}^{\times}$, on a 
$\varphi_{s_0} = \varphi_{s_0} \circ \widetilde{\pi}^{\Phi} (g)$. 
Fixons $g$ dans $\mathcal{O}_{\mathcal{D}}^{\times}$. 
Soit $h$ dans $W_{s_0}^{\nu}$, alors :
$$
\widetilde{\pi}^{\Phi} (g).h = \widetilde{\pi} ({\Phi} (g)).h 
= \overline{\chi}_0^{\Phi^{\nu+1}} ({\rm det} (g)) {\rm St}_{\overline{G}_{s_0}} (\Phi (g)).h
$$
On notera ${\rm St}_{\overline{G}_{s_0}}^{\Phi}$ la représentation 
${\rm St}_{\overline{G}_{s_0}}^{\Phi} : x \mapsto {\rm St}_{\overline{G}_{s_0}} (\Phi (x))$. 
Nous allons montrer que ${\rm St}_{\overline{G}_{s_0}}$ et ${\rm St}_{\overline{G}_{s_0}}^{\Phi}$ 
sont des représentations équivalentes de $\overline{G}_{s_0}$. 
On choisit le modèle suivant pour ${\rm St}_{\overline{G}_{s_0}}$ : 
$$
{\rm St}_{\overline{G}_{s_0}} 
= \left( {\rm Ind}_{\overline{B}_{s_0}}^{\overline{G}_{s_0}} \mathds{1}_{\overline{B}_{s_0}} \right) / 
\mathds{1}_{\overline{G}_{s_0}}
$$
On remarque immédiatement que 
$\mathds{1}_{\overline{G}_{s_0}}^{\Phi}$ et $\mathds{1}_{\overline{G}_{s_0}}$ sont équivalentes. 
Soit :
$$
\psi : {\rm Ind}_{\overline{B}_{s_0}}^{\overline{G}_{s_0}} \mathds{1}_{\overline{B}_{s_0}} 
\rightarrow {\rm Ind}_{\overline{B}_{s_0}}^{\overline{G}_{s_0}} \mathds{1}_{\overline{B}_{s_0}} , 
h \mapsto h \circ \Phi^{-1}
$$
On vérifie alors facilement que $\psi$ est un opérateur d'entrelacement bijectif 
entre ${\rm St}_{\overline{G}_{s_0}}$ et ${\rm St}_{\overline{G}_{s_0}}^{\Phi}$.\\
On en déduit que, pour tout $g$ dans $\mathcal{O}_{\mathcal{D}}^{\times}$, on a :
$$
\forall h \in W_{s_0}^{\nu}, \widetilde{\pi}^{\Phi} (g).h  
= \overline{\chi}_0^{\Phi^{\nu+1}} \otimes {\rm St}_{\overline{G}_{s_0}} (g).h 
\Rightarrow \forall h \in W_{s_0}^{\nu}, \widetilde{\pi}^{\Phi} (g).h \in W_{s_0}^{\nu+1}
$$
Par suite, pour tout $h$ dans $W_{s_0}^{\nu}$, on a :
$$
\varphi_{s_0} (h) = \varphi_{s_0}^{\nu} (h) 
= \varphi_{s_0} (\widetilde{\pi}^{\Phi} (g).h) 
= \varphi_{s_0}^{\nu+1} (\widetilde{\pi}^{\Phi} (g).h)
$$
En particulier, pour $g = 1$ dans $\mathcal{O}_{\mathcal{D}}^{\times}$ et pour tout 
$h$ dans $W_{s_0}^{\nu}$, on a   
$\varphi_{s_0}^{\nu} (h) 
= \varphi_{s_0}^{\nu+1} (h)$. 
On en déduit que l'application $(\varphi_s^{\nu})_{s \in X_0, 0 \leq \nu \leq f-1} \mapsto \varphi_{s_0}^0$ 
est injective.

\end{démo}

\begin{lem}\label{ModeleFormeLineaireS0Induite}
Supposons que $\varphi_{s_0}^{\nu_1, \nu_2}$ est non nulle, en particulier, 
$0 \leq \nu_1 < \nu_2 \leq f-1$, $f$ est pair et $\nu_2 = \nu_1 + \frac{f}{2}$. 
Alors, on peut supposer que  
pour tout $x$ dans $k_{\mathcal{D}}^{\times}$ et tout $b$ dans $\overline{B}_{s_0}$ :
$$
\widetilde{\varphi}_{s_0}^{\nu_1, \nu_2} (bx) 
= \overline{\chi}_1^{-1} \otimes \overline{\chi}_2^{-1} (b)
$$
(et ceci ne dépend pas du représentant $bx$ dans 
$\overline{G}_{s_0} = \overline{B}_{s_0} k_{\mathcal{D}}^{\times}$).
\end{lem}

\begin{démo}
On suppose ici que $\overline{\chi}_1 \times \overline{\chi}_2$ est trivial sur 
$k_{\Delta}^{\times}$ et que l'on peut trouver $\varphi_{s_0}^{\nu_1, \nu_2}$ non triviale. 
Puisque $\varphi_{s_0}^{\nu_1, \nu_2}$ est $k_{\mathcal{D}}^{\times}$-équivariante, 
l'application $\widetilde{\varphi}_{s_0}^{\nu_1, \nu_2}$ 
est aussi $k_{\mathcal{D}}^{\times}$-invariante (à droite). 
Le groupe $\overline{G}_{s_0} = {\rm GL}_2 (k_{\Delta})$ agit transitivement sur 
$\mathbb{P}_{s_0}^1 (k_{\Delta})$ et $\overline{B}_{s_0}$ est le stabilisateur de la droite 
$[1 : 0]$. On en déduit la bijection suivante :
$$
\mathbb{P}_{s_0}^1 (k_{\Delta}) \rightarrow \overline{G}_{s_0} / \overline{B}_{s_0}, 
g.[1: 0] \mapsto g \overline{B}_{s_0}
$$
On peut donc voir l'action de $k_{\mathcal{D}}^{\times}$ sur $\mathbb{P}_{s_0}^1 (k_{\Delta})$ 
comme une action sur $\overline{G}_{s_0} / \overline{B}_{s_0}$ : si 
$d \in \mathbb{P}_{s_0}^1 (k_{\Delta})$, il existe $g$ dans $\overline{G}_{s_0}$ tel que 
$d = g.[1:0]$, alors, pour $x$ dans $k_{\mathcal{D}}^{\times}$, $x.d \in \mathbb{P}_{s_0}^1 (k_{\Delta})$ 
s'identifie à $(xg) \overline{B}_{s_0}$ dans $\overline{G}_{s_0} / \overline{B}_{s_0}$. 
L'action de $k_{\mathcal{D}}^{\times}$ étant transitive, pour tout $g$ dans 
$\overline{G}_{s_0}$, il existe $x$ dans $k_{\mathcal{D}}^{\times}$ tel que 
$(xg) \overline{B}_{s_0} = \overline{B}_{s_0}$, ainsi $g \in k_{\mathcal{D}}^{\times} \overline{B}_{s_0}$. 
On a donc 
$\overline{G}_{s_0} = k_{\mathcal{D}}^{\times} \overline{B}_{s_0} 
= \overline{B}_{s_0} k_{\mathcal{D}}^{\times}$
avec $k_{\mathcal{D}}^{\times} \cap \overline{B}_{s_0} = k_{\Delta}^{\times}$. 
Ainsi, pour tout $x$ dans $k_{\mathcal{D}}^{\times}$ et tout $b$ dans $\overline{B}_{s_0}$, 
$\widetilde{\varphi}_{s_0}^{\nu_1, \nu_2} (bx) 
= \overline{\chi}_1^{-1} \otimes \overline{\chi}_2^{-1} (b)$ 
(et ceci ne dépend pas du représentant $bx$ dans 
$\overline{G}_{s_0} = \overline{B}_{s_0} k_{\mathcal{D}}^{\times}$ car 
$\overline{\chi}_1^{-1} \otimes \overline{\chi}_2^{-1}$ est trivial sur 
$k_{\Delta}^{\times}$).

\end{démo}

\begin{lem}\label{Lemme1Mult1Induites}
Supposons que $\varphi_{s_0}^{\nu_1, \nu_2}$ soit non nulle. Alors, pour tout $k$ dans 
$\mathbb {Z}$, la forme linéaire $\varphi_{s_k}^{\nu_1, \nu_2}$ est entièrement déterminée par 
$\varphi_{s_0}^{\nu_1, \nu_2}$.
\end{lem}

\begin{démo}
Nous allons raisonner par récurrence sur $k$. 
On fixe $0 \leq \nu_1 < \nu_2 \leq f-1$ tels que 
$\varphi_{s_0}^{\nu_1, \nu_2}$ est non nulle. On suppose que 
pour $x \in k_{\mathcal{D}}^{\times}$ et $b \in \overline{B}_{s_0}$,  
$\widetilde{\varphi}_{s_0}^{\nu_1, \nu_2} (bx) 
= \overline{\chi}_1^{-1} \otimes \overline{\chi}_2^{-1} (b)$. 
Notons $a_0$ l'arête $\{ s_0, s_1 \}$, alors 
$\varphi_{s_0}^{\nu_1, \nu_2}$ et $\varphi_{s_1}^{\nu_1, \nu_2}$ 
coïncident sur  
$J_0^{\nu_1, \nu_2} = (U_{s_0}^{\nu_1, \nu_1})^{\mathcal{U}_{a_0}^1} 
\simeq (U_{s_1}^{\nu_1, \nu_1})^{\mathcal{U}_{a_0}^1}$. 
Soit $h$ dans $J_0^{\nu_1, \nu_2}$ dont le support est contenu dans 
$\overline{B}_{s_0}$, alors :
$$
\varphi_{s_0}^{\nu_1, \nu_2} (h) 
= \sum_{x \in \overline{B}_{s_0}} h(x) \widetilde{\varphi}_{s_0}^{\nu_1, \nu_2} (x) 
= \sum_{x \in \overline{B}_{s_0}} 
 \overline{\chi}_1 \otimes \overline{\chi}_2 (x) h(1) 
\overline{\chi}_1^{-1} \otimes \overline{\chi}_2^{-1} (x) 
= \sharp (\overline{B}_{s_0}) \times h(1)
$$
Soit :
\begin{eqnarray*}
D_0 
& = & \{ h \in J_0^{\nu_1, \nu_2} \vert {\rm Supp} (h) \subseteq \overline{B}_{s_0} \}\\
& = &  \{ h : \overline{G}_{s_0} \rightarrow \mathbb{C} \vert 
{\rm Supp} (h) \subseteq \overline{B}_{s_0} \, \, \text{et pour tout} \, \, 
b \in \overline{B}_{s_0}, h(b) = \overline{\chi}_1 \otimes \overline{\chi}_2 (b) h(1) \}
\end{eqnarray*}
Alors $D_0$ est une droite vectorielle. On note $f_0$ l'élément de $D_0$ tel que 
$f_0(1) = 1$, alors :
$$
D_0 = \mathbb{C} f_0 \, \, \text{et} \, \, 
\varphi_{s_0}^{\nu_1, \nu_2} (f_0) = \mu_0 = \sharp (\overline{B}_{s_0}) \neq 0
$$
On en déduit que $\varphi_{s_0}^{\nu_1, \nu_2}$ est non nulle sur le module de Jacquet 
$J_0^{\nu_1, \nu_2}$. 
Par suite, $\varphi_{s_1}^{\nu_1, \nu_2}$ est non nulle sur le module de Jacquet 
$J_0^{\nu_1, \nu_2}$ et 
$\varphi_{s_1}^{\nu_1, \nu_2} (f_0) = \mu_0 
= \varphi_{s_0}^{\nu_1, \nu_2} (f_0)$. 
Ainsi, $\varphi_{s_1}^{\nu_1, \nu_2}$ est une forme linéaire non nulle de 
$U_{s_1}^{\nu_1, \nu_2}$. Notons $H_{s_1}^{\nu_1, \nu_2}$ son noyau. Il s'agit d'un hyperplan de 
$U_{s_1}^{\nu_1, \nu_2}$. De plus, si $x$ appartient à $H_{s_1}^{\nu_1, \nu_2} \cap D_0$, alors 
il existe $\lambda$ dans $\mathbb{C}$ tel que $x = \lambda f_0$ et  
$\varphi_{s_1}^{\nu_1, \nu_2} (x) = 0 = \lambda \mu_0$, donc $\lambda = 0$ et $x=0$. 
On en déduit que :
$$
U_{s_1}^{\nu_1, \nu_2} = D_0 \oplus H_{s_1}^{\nu_1, \nu_2} 
= D_0 \oplus {\rm ker} (\varphi_{s_1}^{\nu_1, \nu_2})
$$
et $\varphi_{s_1}^{\nu_1, \nu_2}$ est entièrement déterminé par 
$\varphi_{s_1}^{\nu_1, \nu_2} (f_0) = \varphi_{s_0}^{\nu_1, \nu_2} (f_0) = \mu_0$. 
Puis, on remarque que 
$\varphi_{s_1}^{\nu_1, \nu_2}$ et $\varphi_{s_2}^{\nu_1, \nu_2}$ 
coïncident sur le module de Jacquet 
$J_1^{\nu_1, \nu_2} = (U_{s_1}^{\nu_1, \nu_1})^{\mathcal{U}_{a_1}^1} 
\simeq (U_{s_2}^{\nu_1, \nu_1})^{\mathcal{U}_{a_1}^1}$ où 
$a_1$ est l'arête $\{ s_1, s_2 \}$ et que pour tout 
$h$ dans $J_1^{\nu_1, \nu_2}$ dont le support est contenu dans $\overline{B}_{s_1}$, on a 
$\varphi_{s_1}^{\nu_1, \nu_2} (h) = \mu_0 \times h(1)$.\\
Avec un raisonnement analogue, on montre que 
$U_{s_2}^{\nu_1, \nu_2} = D_1 \oplus H_{s_2}^{\nu_1, \nu_2} 
= D_1 \oplus {\rm ker} (\varphi_{s_2}^{\nu_1, \nu_2})$ 
où $D_1 = \mathbb{C} f_1$ et $f_1 \in J_1^{\nu_1, \nu_2}$ 
est la fonction dont le support est contenu dans $\overline{B}_{s_1}$ et telle que 
$f_1 (1) = 1$. Ainsi, $\varphi_{s_2}^{\nu_1, \nu_2}$ est entièrement déterminée par 
$\varphi_{s_1}^{\nu_1, \nu_2} (f_1) = \mu_0 = \varphi_{s_0}^{\nu_1, \nu_2} (f_0)$.\\
Une récurrence immédiate nous permet d'obtenir le résultat.

\end{démo}

Nous allons compléter les résultats des parties \ref{EtudeCNSommetS0}, \ref{EtudeCNSommetS1} et 
\ref{EtudeCNSommetSk} en étudiant l'action de l'uniformisante de $\mathcal{D}$. Montrons le résultat suivant :

\begin{propo}\label{ActionUniformisantefPair}
Supposons que $f$ est pair. Alors, pour tout $\nu_1$ 
dans $\lbrace 0, \cdots, (f/2)-1 \rbrace$, on a :
$$
\pi (\varpi_{\mathcal{D}}^{-1}). U_{s_0}^{\nu_1, \nu_1 + f/2} 
= U_{s_0}^{\nu_1 +1, \nu_1 + (f/2) +1}
$$
\end{propo}

\begin{démo}
Nous allons utiliser plusieurs lemmes. 
Pour cela, rappelons quelques notations. 
On suppose ici que $f$ est pair. 
On note $\Phi$ un générateur de ${\rm Gal} (k_{\Delta} / k)$. On peut supposer que 
l'action de $\Phi$ sur $\mathcal{A}_{s_0}^{\times} = {\rm GL}_2 (\mathcal{O}_{\Delta})$ est 
induite par la conjugaison 
par $\varpi_{\Delta}$. 
Soit $\mathcal{K}_{s_0} = \mathcal{R}_1 
= \langle \varpi_{\Delta} \rangle \mathcal{A}_{s_0}^{\times}$. 
Le $\mathcal{K}_{s_0}$-module $V^{\mathcal{U}_{s_0}^1}$ 
se décompose simplement de la façon suivante :
$$
V_{s_0} = V^{\mathcal{U}_{s_0}^1} 
= (\bigoplus_{\nu = 0}^{f-1} W_{s_0}^{\nu}) \oplus (\bigoplus_{0 \leq \nu_1 < \nu_2 \leq f-1} 
U_{s_0}^{\nu_1, \nu_2})
$$
où :
$$
W_{s_0}^{\nu} \simeq \overline{\chi}_0^{\Phi^{\nu}} \otimes {\rm St}_{\overline{G}_{s_0}} 
\, \, \text{et} \, \, 
U_{s_0}^{\nu_1, \nu_2} 
\simeq {\rm Ind}_{\overline{B}_{s_0}}^{\overline{G}_{s_0}} 
(\overline{\chi}_0^{\Phi^{\nu_1}} \otimes \overline{\chi}_0^{\Phi^{\nu_2}})
$$
On note 
$(\sigma, \Sigma) = 
({\rm Ind}_{\overline{B}_{s_0}}^{\overline{G}_{s_0}} 
(\overline{\chi}_0 \otimes \overline{\chi}_0^{\Phi^{f/2}}), 
U_{s_0}^{0, f/2} )$
vue comme représentation de $\mathcal{A}_{s_0}^{\times} \langle \varpi_{\mathbb{K}} \rangle$ 
(en rajoutant $\chi_{\vert \mathbb{K}^{\times}}$ sur le centre). 
Alors :
$$
N_{\mathcal{R}_1} (\sigma) 
= \{ g \in \mathcal{R}_1 : \sigma^g 
\simeq_{\mathcal{A}_{s_0}^{\times} \langle \varpi_{\mathbb{K}} \rangle - module} \sigma \} 
= \langle \varpi_{\Delta}^{f/2} \rangle \mathcal{A}_{s_0}^{\times}
$$
Ainsi $\sigma^{\varpi_{\Delta}^{f/2}} \simeq \sigma$, et il existe un opérateur d'entrelacement bijectif 
$J : U_{s_0}^{0, f/2} \rightarrow U_{s_0}^{0, f/2}$ tel que pour tout  
$h \in \mathcal{A}_{s_0}^{\times}$, on a $\sigma (\varpi_{\Delta}^{f/2} h \varpi_{\Delta}^{-f/2}) 
= J \sigma (h) J^{-1}$. 
On fixe $J$ comme dans \ref{TypeEtenduMaxfPair}. 
On définit \mbox{$(\widehat{\sigma}, \widehat{\Sigma}) = 
(\widehat{\sigma}, U_{s_0}^{0, f/2})$}
une extension de $(\sigma, \Sigma)$ à $N = N_{\mathcal{R}_1} (\sigma)$ telle que 
pour tout $i$ dans $\mathbb{Z}$ et tout $h$~dans~$\mathcal{A}_{s_0}^{\times}$, on a 
\mbox{$\widehat{\sigma} ((\varpi_{\Delta}^{f/2})^i h) 
= \zeta^i J^i \sigma (h)$}, 
où $\zeta = \chi_f ((-1)^{e-1} \varpi_{\mathbb{K}})$ et $e = d/f$. 
Enfin, on définit :
$$
(\widetilde{\sigma}, \widetilde{\Sigma}) 
= ({\rm c-Ind}_N^{\mathcal{R}_1} \widehat{\sigma}, 
{\rm c-Ind}_N^{\mathcal{R}_1} U_{s_0}^{0, f/2})
$$
représentation de $\mathcal{R}_1$. Alors 
$(\widetilde{\sigma}, \widetilde{\Sigma})$ est un type étendu maximal de 
niveau $0$ pour $\pi$. 
Ainsi :
$$
{\rm Hom}_{\mathcal{K}_{s_0}} (\widetilde{\sigma}, \pi) \neq 0
$$ 
Il existe donc une injection $I : \widetilde{\Sigma} \hookrightarrow V_{\pi}$ telle que 
pour tout $v \in \widetilde{\Sigma}$ et tout $ h \in 
\mathcal{K}_{s_0} = \mathcal{R}_1$, on a 
$I(\widetilde{\sigma} (h).v) = \pi (h). I(v)$. 
Par construction, $\widetilde{\Sigma} \subseteq V_{s_0}$. On peut supposer que 
$I$ est une inclusion, $\widetilde{\Sigma} \subseteq V_{s_0}$ et donc 
$\widetilde{\sigma} (h).v = \pi (h). v$ 
pour tout $v$ dans $\widetilde{\Sigma} \subseteq V_{s_0}$ et pour tout $h$ dans  
$\mathcal{K}_{s_0}$.\\

On démontre facilement les trois lemmes suivants :

\begin{lem}
On a une injection $\mathcal{A}_{s_0}^{\times}$-équivariante :
$$
\Sigma = U_{s_0}^{0, f/2} \hookrightarrow 
\widetilde{\Sigma} = {\rm c-Ind}_N^{\mathcal{R}_1} U_{s_0}^{0, f/2} , 
w \mapsto f_w
$$
où $f_w : \mathcal{R}_1 \rightarrow U_{s_0}^{0, f/2}$ a pour support 
$N = N_{\mathcal{R}_1} (\sigma)$ et 
pour tout $n$ dans $N$, $f_w (n) = \widehat{\sigma} (n).w$.
\end{lem}

\begin{lem}
Pour tout $w$ dans $U_{s_0}^{0, f/2}$,
$\widetilde{\sigma} (\varpi_{\Delta}^{-1}). f_w$ est une fonction à support 
dans $N \varpi_{\Delta}$, et pour tout $n$ dans $N$ :
$$
\widetilde{\sigma} (\varpi_{\Delta}^{-1}). f_w 
(n \varpi_{\Delta}) = \widehat{\sigma} (n).w
$$
\end{lem}

\begin{lem}
Pour tout $w^{'}$ dans $U_{s_0}^{1, (f/2) +1}$, on a 
$w^{'} \circ \Phi^{-1} \in U_{s_0}^{0, f/2}$.
\end{lem}

\begin{lem}
On a une injection $\mathcal{A}_{s_0}^{\times}$-équivariante :
$$
\Sigma = U_{s_0}^{1, (f/2)+1} \hookrightarrow 
\widetilde{\Sigma} = {\rm c-Ind}_N^{\mathcal{R}_1} U_{s_0}^{0, f/2} , 
w^{'} \mapsto g_{w^{'}}
$$
où $g_{w^{'}} : \mathcal{R}_1 \rightarrow U_{s_0}^{0, f/2}$ a pour support $N \varpi_{\Delta}$ et 
pour tout $n$ dans $N$, $g_{w^{'}} (n \varpi_{\Delta}) = \widehat{\sigma} (n).(w^{'} \circ \Phi^{-1})$.
\end{lem}

\begin{démo}
Il faut montrer que l'application $w^{'} \mapsto g_{w^{'}}$ 
est bien $\mathcal{A}_{s_0}^{\times}$-équivariante, c'est-à-dire que pour tout 
$h \in \mathcal{A}_{s_0}^{\times}$ et tout 
$w^{'} \in U_{s_0}^{1, (f/2)+1}$, on a :
$$ 
g_{{\rm Ind}_{\overline{B}_{s_0}}^{\overline{G}_{s_0}} 
(\overline{\chi}_0^{\Phi} \otimes \overline{\chi}_0^{\Phi^{(f/2) +1}}) (h).w^{'}} 
= \widetilde{\sigma} (h). g_{w^{'}}
$$
Soient $h$ dans $\mathcal{A}_{s_0}^{\times}$ et $w^{'}$ dans 
$U_{s_0}^{1, (f/2)+1}$. 
Alors $g_{{\rm Ind}_{\overline{B}_{s_0}}^{\overline{G}_{s_0}} 
(\overline{\chi}_0^{\Phi} \otimes \overline{\chi}_0^{\Phi^{(f/2) +1}}) (h).w^{'}}$ 
est une fonction à support dans~$N \varpi_{\Delta}$. Soit $n$ dans $N$. 
Soit $i$ dans $\mathbb{Z}$ et $n_0$ dans $\mathcal{A}_{s_0}^{\times}$ tels que 
$n = (\varpi_{\Delta}^{f/2})^i n_0$. On a :
$$
g_{{\rm Ind}_{\overline{B}_{s_0}}^{\overline{G}_{s_0}} 
(\overline{\chi}_0^{\Phi} \otimes \overline{\chi}_0^{\Phi^{(f/2) +1}}) (h).w^{'}} (n \varpi_{\Delta}) 
= \widehat{\sigma} (n) . (({\rm Ind}_{\overline{B}_{s_0}}^{\overline{G}_{s_0}} 
(\overline{\chi}_0^{\Phi} \otimes \overline{\chi}_0^{\Phi^{(f/2) +1}}) (h).w^{'}) \circ \Phi^{-1})
$$
or 
$\widehat{\sigma} (n) = \zeta^i J^i \sigma (n_0)$ 
et :
$$
\sigma (n_0) . (({\rm Ind}_{\overline{B}_{s_0}}^{\overline{G}_{s_0}} 
(\overline{\chi}_0^{\Phi} \otimes \overline{\chi}_0^{\Phi^{(f/2) +1}}) (h).w^{'}) \circ \Phi^{-1}) : 
x \mapsto 
({\rm Ind}_{\overline{B}_{s_0}}^{\overline{G}_{s_0}} 
(\overline{\chi}_0^{\Phi} \otimes \overline{\chi}_0^{\Phi^{(f/2) +1}}) (h).w^{'}) \circ \Phi^{-1} (x n_0)
$$
où :
$$
({\rm Ind}_{\overline{B}_{s_0}}^{\overline{G}_{s_0}} 
(\overline{\chi}_0^{\Phi} \otimes \overline{\chi}_0^{\Phi^{(f/2) +1}}) (h).w^{'}) \circ \Phi^{-1} :
 X \mapsto 
({\rm Ind}_{\overline{B}_{s_0}}^{\overline{G}_{s_0}} 
(\overline{\chi}_0^{\Phi} \otimes \overline{\chi}_0^{\Phi^{(f/2) +1}}) (h).w^{'}) (\Phi^{-1} (X)) 
= w^{'} (\Phi^{-1} (X) h)
$$
On en déduit que :
$$
g_{{\rm Ind}_{\overline{B}_{s_0}}^{\overline{G}_{s_0}} 
(\overline{\chi}_0^{\Phi} \otimes \overline{\chi}_0^{\Phi^{(f/2) +1}}) (h).w^{'}} (n \varpi_{\Delta}) 
= \zeta^i J^i \sigma (n_0). 
(({\rm Ind}_{\overline{B}_{s_0}}^{\overline{G}_{s_0}} 
(\overline{\chi}_0^{\Phi} \otimes \overline{\chi}_0^{\Phi^{(f/2) +1}}) (h).w^{'}) \circ \Phi^{-1})
$$
avec :
$$
\sigma (n_0). 
(({\rm Ind}_{\overline{B}_{s_0}}^{\overline{G}_{s_0}} 
(\overline{\chi}_0^{\Phi} \otimes \overline{\chi}_0^{\Phi^{(f/2) +1}}) (h).w^{'}) \circ \Phi^{-1}) :
\overline{G_{s_0}} \rightarrow \mathbb{C}, 
x \mapsto w^{'} (\Phi^{-1} (x) \Phi^{-1} (n_0) h)
$$
On a 
$\widetilde{\sigma} (h). g_{w^{'}} : 
x \mapsto g_{w^{'}} (xh)$. 
On remarque que pour tout $n$ dans $N$, 
$n \varpi_{\Delta} h = n \Phi (h) \varpi_{\Delta}$. 
Ainsi $\widetilde{\sigma} (h). g_{w^{'}}$ est une fonction à support dans 
$N \varpi_{\Delta}$. Soit $n$ dans $N$. 
Soit $i$ dans $\mathbb{Z}$ et $n_0$ dans $\mathcal{A}_{s_0}^{\times}$ tels que 
$n = (\varpi_{\Delta}^{f/2})^i n_0$. Alors :
$$
\widetilde{\sigma} (h). g_{w^{'}} (n \varpi_{\Delta}) 
= g_{w^{'}} (n \varpi_{\Delta} h) 
= g_{w^{'}} (n \Phi (h) \varpi_{\Delta}) 
= \widehat{\sigma} (n \Phi (h)).(w^{'} \circ \Phi^{-1}) 
= \zeta^i J^i \sigma (n_0 \Phi (h)). (w^{'} \circ \Phi^{-1})
$$
avec 
$\sigma (n_0 \Phi (h)). (w^{'} \circ \Phi^{-1}) : x 
\mapsto w^{'} \circ \Phi^{-1} (x n_0 \Phi (h)) 
= w^{'} (\Phi^{-1} (x) \Phi^{-1} (n_0) h)$. 
Par suite :
$$
g_{{\rm Ind}_{\overline{B}_{s_0}}^{\overline{G}_{s_0}} 
(\overline{\chi}_0^{\Phi} \otimes \overline{\chi}_0^{\Phi^{(f/2) +1}}) (h).w^{'}} 
= \widetilde{\sigma} (h). g_{w^{'}}
$$

\end{démo}

\begin{lem}
 On a :
$$
\widetilde{\sigma} (\varpi_{\Delta}^{-1}). U_{s_0}^{0, f/2} 
= U_{s_0}^{1, (f/2) +1}
$$
\end{lem}

\begin{démo}
En utilisant les notations précédentes, on remarque que pour tout 
$w$ dans $U_{s_0}^{0, f/2}$, on a 
$\widetilde{\sigma} (\varpi_{\Delta}^{-1}). f_w 
= g_{w^{'}}$, 
où $w^{'} = w \circ \Phi \in U_{s_0}^{1, (f/2)+1}$. 
Par conséquent, $\widetilde{\sigma} (\varpi_{\Delta}^{-1}). U_{s_0}^{0, f/2} 
\subseteq U_{s_0}^{1, (f/2) +1}$. Les deux espaces 
$U_{s_0}^{0, f/2}$ et $U_{s_0}^{1, (f/2)+1}$ sont de même dimension (finie), par suite :
$$
\widetilde{\sigma} (\varpi_{\Delta}^{-1}). U_{s_0}^{0, f/2} 
= U_{s_0}^{1, (f/2) +1}
$$
\end{démo}

Par récurrence, on montre de façon analogue que 
pour tout $0 \leq \nu_1 \leq (f/2)-1$, on a :
$$
\widetilde{\sigma} (\varpi_{\Delta}^{-1}). U_{s_0}^{\nu_1, \nu_1 + f/2} 
= U_{s_0}^{\nu_1 +1, \nu_1 +(f/2) +1}
$$
Enfin, on peut montrer le résultat \ref{ActionUniformisantefPair}. 
En effet, on remarque que 
$\varpi_{\mathcal{D}}^{d} = \varpi_{\mathbb{F}} 
= \varpi_{\mathbb{K}}^2 =  
\varpi_{\Delta}^{d}$. 
Ainsi, il existe $u_0$ dans $\mathcal{A}_{s_0}^{\times} = {\rm GL}_2 (\mathcal{O}_{\Delta})$ tel que 
$\varpi_{\mathcal{D}} = \varpi_{\Delta} u_0$ et $\varpi_{\mathcal{D}} \in \mathcal{R}_1$. 
Par conséquent :
$$
\widetilde{\sigma} (\varpi_{\mathcal{D}}^{-1}) 
= \widetilde{\sigma} (u_0^{-1}) \widetilde{\sigma} (\varpi_{\Delta}^{-1})
$$
Soit $0 \leq \nu_1 \leq (f/2)-1$. Alors 
$U_{s_0}^{\nu_1 +1, \nu_1 +(f/2) +1}$ est un $\mathcal{A}_{s_0}^{\times}$-module, donc :
$$
\widetilde{\sigma} (u_0^{-1}) . U_{s_0}^{\nu_1 +1, \nu_1 +(f/2) +1}
= U_{s_0}^{\nu_1 +1, \nu_1 +(f/2) +1}
$$
et :
$$
\widetilde{\sigma} (\varpi_{\mathcal{D}}^{-1}). U_{s_0}^{\nu_1, \nu_1 + f/2} 
= \widetilde{\sigma} (u_0^{-1}) \widetilde{\sigma} (\varpi_{\Delta}^{-1}) . 
U_{s_0}^{\nu_1, \nu_1 + f/2} = \widetilde{\sigma} (u_0^{-1}) . U_{s_0}^{\nu_1 +1, \nu_1 +(f/2) +1}
= U_{s_0}^{\nu_1 +1, \nu_1 +(f/2) +1}
$$

\end{démo}

\begin{propo}\label{ResultatMult1Induite}
Soit $Y_{{\rm Ind}}$ l'ensemble des formes linéaires $\varphi$ dans ${\rm ker} (\partial_1^{\ast})$ 
telles que pour tout sommet $s$, pour tout $0 \leq \nu \leq f-1$, on a $\varphi_s^{\nu} = 0$. 
Si $f$ est impair, $Y_{{\rm Ind}} = \{ 0 \}$.\\
Supposons que $f$ est pair. Alors 
l'application :
$$
Y_{{\rm Ind}} \rightarrow {\rm Hom}_{k_{\mathcal{D}}^{\times}} (\pi_{s_0}^{0, f/2}, \mathds{1}), 
\varphi \mapsto \varphi_{s_0}^{0, f/2}
$$  
est injective.
\end{propo}

\begin{démo}
Si $f$ est impair, alors d'après \ref{SyntheseCNDistinctionTotRam}, on a bien 
$Y_{{\rm Ind}} = \{ 0 \}$.\\
Supposons que $f$ est pair. Toujours d'après \ref{SyntheseCNDistinctionTotRam}, on sait que 
$\varphi_{s}^{\nu_1, \nu_2} = 0$ dès que $\nu_2 \neq \nu_1 + (f/2)$. 
On a 
$\varphi_{s_0} = \varphi_{\varpi_{\mathcal{D}}^{-1}. s_0} \circ \pi (\varpi_{\mathcal{D}}^{-1}) 
= \varphi_{s_0} \circ \pi (\varpi_{\mathcal{D}}^{-1})$. 
Soit $\nu_1$ dans $\{ 0, \cdots, (f/2)-1 \}$ et $s$ un sommet de $X_{\mathbb{K}}$. 
Soit $w$ dans $U_{s_0}^{\nu_1, \nu_1 + f/2}$, alors 
$\varphi_{s_0} (w) = 
\varphi_{s_0}^{\nu_1, \nu_1+ (f/2)} (w)$ 
et, d'après \ref{ActionUniformisantefPair} :
$$
\pi (\varpi_{\mathcal{D}}^{-1}).w = \widetilde{\sigma} (\varpi_{\mathcal{D}}^{-1}).w \in 
 U_{s_0}^{\nu_1 +1, \nu_1 +(f/2) +1}
$$
On en déduit que pour tout $w$ dans $U_{s_0}^{\nu_1, \nu_1 + f/2}$ :
$$
\varphi_{s_0}^{\nu_1, \nu_1+ (f/2)} (w) 
= \varphi_{s_0}^{\nu_1 +1, \nu_1+1+ (f/2)} (\widetilde{\sigma} (\varpi_{\mathcal{D}}^{-1}).w)
$$
Enfin, le lemme \ref{Lemme1Mult1Induites} nous permet de conclure. 

\end{démo}

\begin{theo}\label{TheoremeMultiplicite1NonCusp}
Soit $(\pi, V)$ une représentation lisse et irréductible de $G$, membre de la 
série discrète, de niveau $0$, non cuspidale. Alors :
$$
{\rm dim}_{\mathbb{C}} ({\rm Hom}_{\mathcal{D}^{\times}} (\pi, \mathds{1})) \leq 1
$$
\end{theo}

\begin{démo}
Si $\pi$ n'est pas $\mathcal{D}^{\times}$-distinguée, alors 
${\rm dim}_{\mathbb{C}} ({\rm Hom}_{\mathcal{D}^{\times}} (\pi, \mathds{1})) = 0$.\\
Sinon, si $\pi$ est $\mathcal{D}^{\times}$-distinguée, alors 
${\rm dim}_{\mathbb{C}} ({\rm Hom}_{\mathcal{D}^{\times}} (\pi, \mathds{1})) 
= {\rm dim}_{\mathbb{C}} ({\rm ker} (\partial_1^{\ast}))$. 
D'après \ref{SyntheseCNDistinctionTotRam}, si $\varphi \in {\rm ker} (\partial_1^{\ast})$, alors, 
ou bien $\varphi \in Y_{{\rm St}}$, ou bien $\varphi \in Y_{{\rm Ind}}$. 
Or, d'après \ref{ResultatMultiplicite1Steinberg}, \ref{DimensionEntrelacementSteinberg}, \ref{ResultatMult1Induite} 
et \ref{DimensionEntrelacementInduite}, les espaces 
$Y_{{\rm St}}$ et $Y_{{\rm Ind}}$ sont de dimension au plus $1$. 
On en déduit que 
${\rm dim}_{\mathbb{C}} ({\rm Hom}_{\mathcal{D}^{\times}} (\pi, \mathds{1})) \leq 1$.\\

\end{démo}
\subsection{Synthèse sur les conditions de $\mathcal{D}^{\times}$-distinction de $\pi$.}

Nous allons démontrer que, si $f$ est pair, alors la représentation $\pi$ ne peut pas être 
$\mathcal{D}^{\times}$-distinguée. 
On utilise les notations introduites en \ref{NotationHypotheseCNDistinction} et 
\ref{NotationsS0Induites}. 

\begin{nota}
On suppose ici que $f$ est pair. 
On note $\chi_1 = \overline{\chi}_0$, $\chi_2 = \overline{\chi}_0^{\Phi^{f/2}}$. 
On suppose donc que $\chi_1 \times \chi_2$ est trivial sur $k_{\Delta}^{\times}$. 
Alors $\chi = \chi_1 \otimes \chi_2$ est une représentation du tore $\overline{T}_{s_0}$ vue comme une 
représentation de $\overline{B}_{s_0}$ triviale sur $\overline{U}_{s_0}$.  
On notera 
$w = \left(\begin{array}{cc}
       0 & 1\\
       1 & 0
       \end{array} \right)$. 
Soit $\chi = \chi_1 \otimes \chi_2$. 
On remarque que $\chi \circ {\rm Ad} (w) = \chi_2 \otimes \chi_1$. 
On notera $\chi^{w} = \chi \circ {\rm Ad} (w)$. 
Soit $\psi = {\rm Ad} (\varpi_{\Delta}^{f/2})$, alors $\chi_1 \circ \psi = \chi_2$ et 
$\chi_2 \circ \psi = \chi_1$. 
On note $(\sigma, {\rm Ind}_{\overline{B}_{s_0}}^{\overline{G}_{s_0}} \chi)$. 
\end{nota}

\begin{déf}
On définit :
$$
D = \lbrace F : \overline{G}_{s_0} \rightarrow \mathbb{C}, 
F(b_1 g b_2) = \chi^w (b_1) \chi (b_2) F(g); b_1, b_2 \in \overline{B}_{s_0}, g \in \overline{G}_{s_0} \rbrace
$$
Pour $F$ dans $D$ et pour $h : \overline{G}_{s_0} \rightarrow V_1 = \mathbb{C}$, on note :
$$
F \ast h : \overline{G}_{s_0} \rightarrow V_2 = \mathbb{C}, 
x \mapsto \frac{1}{\vert \overline{G}_{s_0} \vert} 
\sum_{g \in \overline{G}_{s_0}} F(x g^{-1}) h(g)
$$

\end{déf}

\begin{rmq}
Si $F \in D$, on a $F(b) = 0$ pour tout $b$ dans $\overline{B}_{s_0}$. 
On vérifie facilement que si $f_1$ appartient à 
${\rm Ind}_{\overline{B}_{s_0}}^{\overline{G}_{s_0}} \chi_1 \otimes \chi_2$, on a 
$F \ast f_1 \in {\rm Ind}_{\overline{B}_{s_0}}^{\overline{G}_{s_0}} \chi_2 \otimes \chi_1$.
\end{rmq}

\begin{déf}
 Soit $F$ dans $D$, on définit 
$L_F : {\rm Ind}_{\overline{B}_{s_0}}^{\overline{G}_{s_0}} \chi \rightarrow {\rm Ind}_{\overline{B}_{s_0}}^{\overline{G}_{s_0}} \chi^w$ 
par $L_F (f_1) = F \ast f_1$.
\end{déf}

On a le résultat classique suivant :

\begin{propo}
(Mackey) Si $F \in D$ alors 
$L_F \in {\rm Hom}_{\overline{G}_{s_0}} ({\rm Ind}_{\overline{B}_{s_0}}^{\overline{G}_{s_0}} \chi, 
{\rm Ind}_{\overline{B}_{s_0}}^{\overline{G}_{s_0}} \chi^w)$ 
et on a un isomorphisme de $\mathbb{C}$-espaces vectoriels :
$$
D \rightarrow {\rm Hom}_{\overline{G}_{s_0}} 
({\rm Ind}_{\overline{B}_{s_0}}^{\overline{G}_{s_0}} \chi, 
{\rm Ind}_{\overline{B}_{s_0}}^{\overline{G}_{s_0}} \chi^w), 
F \mapsto L_F
$$
\end{propo}

\begin{rmq}
Un élément $F$ dans $D$ est entièrement déterminé par la valeur $F(w)$. Ainsi, 
si \mbox{$F(w) \neq 0$}, $L_F$ est un isomorphisme.
\end{rmq}

\begin{déf}
Soit :
$$
I : {\rm Ind}_{\overline{B}_{s_0}}^{\overline{G}_{s_0}} \chi^w 
\rightarrow (\sigma \circ \psi, {\rm Ind}_{\overline{B}_{s_0}}^{\overline{G}_{s_0}} \chi), 
f \mapsto f \circ \psi^{-1}
$$
alors $I$ est un isomorphisme de $\overline{G}_{s_0}$-modules.
\end{déf}

\begin{lem}
Soit $F$ dans $D$ tel que $F(w) \neq 0$. Posons 
$J = I \circ L_F$. Alors 
$J : {\rm Ind}_{\overline{B}_{s_0}}^{\overline{G}_{s_0}} \chi 
\rightarrow {\rm Ind}_{\overline{B}_{s_0}}^{\overline{G}_{s_0}} \chi$ 
est un opérateur d'entrelacement bijectif tel que, pour tout $h$ dans $\mathcal{A}_{s_0}^{\times}$ :
$$
J \circ \sigma (h) = \sigma (\varpi_{\Delta}^{f/2} h \varpi_{\Delta}^{-f/2}) \circ J
$$
\end{lem}

\begin{propo}\label{NonDistinctionCasfPair}
Supposons que $\overline{\chi}_0$ est trivial sur $k^{\times}$. Alors 
$\varphi_{s_0}^{0, f/2}$ est nulle.
\end{propo}

\begin{démo}
Supposons que $\overline{\chi}_0$ est trivial sur $k^{\times}$. 
Raisonnons par l'absurde pour montrer que \mbox{$\varphi_{s_0}^{0, f/2} = 0$}. 
Supposons donc que $\varphi_{s_0}^{0, f/2} \neq 0$. 
D'après le lemme \ref{ModeleFormeLineaireS0Induite}, on peut supposer que, 
pour tout $h$ dans 
\mbox{$U_{s_0}^{0, f/2} = {\rm Ind}_{\overline{B}_{s_0}}^{\overline{G}_{s_0}} \chi$}, on a 
$\varphi_{s_0}^{0, f/2} (h)
= \sum_{g \in \overline{G}_{s_0}} h(g) \widetilde{\varphi}_{s_0}^{0, f/2} (g)$ 
où, pour tout $b$ dans $\overline{B}_{s_0}$ et tout $x$ dans $k_{\mathcal{D}}^{\times}$, 
\mbox{$\widetilde{\varphi}_{s_0}^{0, f/2} (xb) 
= \chi_1^{-1} \otimes \chi_2^{-1} (b)$}. 
Nous allons fixer $J$, un opérateur d'entrelacement bijectif 
\mbox{$J : {\rm Ind}_{\overline{B}_{s_0}}^{\overline{G}_{s_0}} \chi 
\rightarrow {\rm Ind}_{\overline{B}_{s_0}}^{\overline{G}_{s_0}} \chi$} comme dans la construction 
de \ref{TypeEtenduMaxfPair}.\\
On commence par fixer $\mu_0 : (k, +) \rightarrow \mathbb{C}^{\times}$ un caractère additif non trivial de 
$k$. 
On cherche $f_0$, non nulle, dans ${\rm Ind}_{\overline{B}_{s_0}}^{\overline{G}_{s_0}} \chi$ telle que, 
pour tout $u = \left(\begin{array}{cc}
                   1 & u_0\\
                   0 & 1
                   \end{array} \right)$ :
$$
\sigma (u).f_0 = \mu (u) f_0, \, \, \text{où} \, \, 
\mu (u) = \mu_0 \circ {\rm tr}_{k_{\Delta} / k} (u_0)
$$
Rappelons que 
$\overline{G}_{s_0} = \overline{B}_{s_0} \sqcup \overline{T}_{s_0} \overline{U}_{s_0} w \overline{U}_{s_0}$ 
et si $x \in \overline{G}_{s_0} \setminus \overline{B}_{s_0}$, 
alors $x$ s'écrit de façon unique sous la forme $x = tu_1xu_2$ avec $t \in \overline{T}_{s_0}$ et 
$u_1, u_2 \in \overline{U}_{s_0}$. On définit $f_0 : \overline{G}_{s_0} \rightarrow \mathbb{C}$ telle que 
$f_0 (b) = 0$ pour tout $b$ dans $\overline{B}_{s_0}$ et, pour $t \in \overline{T}_{s_0}$ et 
$u_1, u_2 \in \overline{U}_{s_0}$, $f_0 (t u_1 w u_2) = \chi (t) \mu (u_2)$.\\
Alors $f_0$ est bien définie, $f_0 \in {\rm Ind}_{\overline{B}_{s_0}}^{\overline{G}_{s_0}} \chi$ et, pour tout 
$u$ dans $\overline{U}_{s_0}$, $\sigma (u).f_0 = \mu (u) f_0$.\\
On vérifie facilement que si $F \in D$ vérifie $F(w) = -(Q+1)$, alors 
$J.f_0 = f_0$ où $J = I \circ L_F$.\\
On fixe donc $J = I \circ L_F$ où $F(w) = -(Q+1)$.\\
D'après \ref{ActionUniformisantefPair}, on sait que, pour tout $h$ dans 
$U_{s_0}^{0, f/2}$, on a 
$\pi (\varpi_{\Delta}^{-f/2}).h = \sigma (\varpi_{\Delta}^{-f/2}).h 
\in U_{s_0}^{0, f/2}$ et  
$\varphi_{s_0}^{0, f/2} (\pi (\varpi_{\Delta}^{-f/2}).h) 
= \varphi_{\varpi_{\Delta}^{-f/2}.s_0} \circ \pi (\varpi_{\Delta}^{-f/2}) (h) 
= \varphi_{s_0} (h)
= \varphi_{s_0}^{0, f/2} (h)$. 
Or $\pi (\varpi_{\Delta}^{-f/2}).h = \zeta J.h$ 
(où $\zeta$ est défini comme dans \ref{TypeEtenduMaxfPair}) donc 
$\varphi_{s_0}^{0, f/2} (h) = \zeta \varphi_{s_0}^{0, f/2} (J.h)$. \\
Notons $\widetilde{f} \in {\rm Ind}_{\overline{B}_{s_0}}^{\overline{G}_{s_0}} 
\overline{\chi}_0 \otimes \overline{\chi}_0^{\Phi^{f/2}} 
= {\rm Ind}_{\overline{B}_{s_0}}^{\overline{G}_{s_0}} \chi$ telle que 
$\widetilde{f} (g) = 0$ pour $g \notin \overline{B}_{s_0}$ et, pour 
$g \in \overline{B}_{s_0}$ 
\mbox{$\widetilde{f} (g) = \frac{1}{\vert \overline{B}_{s_0} \vert} \chi (g)$}. 
Alors 
$\varphi_{s_0}^{0, f/2} (\widetilde{f}) 
= \sum_{b \in \overline{B}_{s_0}} \widetilde{\varphi}_{s_0}^{0, f/2} (b) \widetilde{f} (b) 
= \frac{1}{\vert \overline{B}_{s_0} \vert} \sum_{b \in \overline{B}_{s_0}} \chi^{-1} (b) \chi (b) 
= 1$. 
On en déduit que 
\mbox{$\varphi_{s_0}^{0, f/2} (\widetilde{f}) = 1 
= \zeta \varphi_{s_0}^{0, f/2} (J.\widetilde{f})$}. 

\begin{itemize}
\item[1.] Calculons tout d'abord $J. \widetilde{f}$.\\
Un calcul simple nous permet de vérifier que $L_F (\widetilde{f}) = \frac{1}{\vert \overline{G}_{s_0} \vert} F$. Ainsi, 
pour tout $b \in \overline{B}_{s_0}$, on a $J.\widetilde{f} (b) = 0$ et si $b_1, b_2 \in \overline{B}_{s_0}$ :
$$
J.\widetilde{f} (b_1 w b_2) = \frac{F(w)}{\vert \overline{G}_{s_0} \vert} 
\chi (b_1) \chi^w (b_2)
$$

\item[2.] Calculons $\varphi_{s_0}^{0, f/2} (J. \widetilde{f})$.\\
On a :
$$
\varphi_{s_0}^{0, f/2} (J. \widetilde{f}) 
= \sum_{g \in \overline{B}_{s_0} w \overline{B}_{s_0}} \widetilde{\varphi}_{s_0}^{0, f/2} (g) 
J. \widetilde{f} (g) 
= \sum_{g \in (k_{\mathcal{D}}^{\times} \backslash k_{\Delta}^{\times}) \overline{B}_{s_0}} 
\widetilde{\varphi}_{s_0}^{0, f/2} (g) J. \widetilde{f} (g) 
$$
Notons :
$$
X_1 = \left\lbrace \left( \begin{array}{cc}
             0 & B\\
             C & D
             \end{array} \right) : B, C \in k_{\Delta}^{\times}, D \in k_{\Delta} \right\rbrace 
\, \, \text{et} \, \, 
X_2 = \left\lbrace \left( \begin{array}{cc}
             A & B\\
             C & D
             \end{array} \right) : A, C \in k_{\Delta}^{\times}, B, D \in k_{\Delta} : 
AD-BC \neq 0 \right\rbrace
$$
Alors $\overline{B}_{s_0} w \overline{B}_{s_0} = \overline{G}_{s_0} \backslash \overline{B}_{s_0} = X_1 \sqcup X_2$. 
Fixons $X = \left( \begin{array}{cc}
             0 & B\\
             C & D
             \end{array} \right) \in X_1$. Alors :
$$
X = \left( \begin{array}{cc}
             0 & B\\
             C & D
             \end{array} \right)
= \left( \begin{array}{cc}
             0 & \alpha^2\\
             1 & 0
             \end{array} \right) 
\left( \begin{array}{cc}
             C & D\\
             0 & B/\alpha^2
             \end{array} \right) 
= \left( \begin{array}{cc}
             B & 0\\
             0 & C
             \end{array} \right) w 
\left( \begin{array}{cc}
             1 & D/C\\
             0 & 1
             \end{array} \right)
$$
Donc :
$$
\widetilde{\varphi}_{s_0}^{0, f/2} (X) 
= \chi_1^{-1} (C) \chi_2^{-1} (B/\alpha^2) 
= \chi_1^{-1} (C) \chi_2^{-1} (B) \chi_2 (\alpha^2)
$$
et :
$$
J. \widetilde{f} (X) 
= \frac{F(w)}{\vert \overline{G}_{s_0} \vert} 
\chi_1 (B) \chi_2 (C) 
$$
Ainsi :
\begin{eqnarray*}
S_1 
& = & \sum_{X \in X_1} \widetilde{\varphi}_{s_0}^{0, f/2} (X) 
   J. \widetilde{f} (X)\\
& = & \frac{F(w)}{\vert \overline{G}_{s_0} \vert} \times \chi_2 (\alpha^2) \times 
   (\sum_{B \in k_{\Delta}^{\times}} \chi_1 (B) \chi_2^{-1} (B)) \times 
   (\sum_{C \in k_{\Delta}^{\times}} \chi_2 (C) \chi_1^{-1} (C))
\end{eqnarray*}
Or 
$\sum_{B \in k_{\Delta}^{\times}} \chi_1 (B) \chi_2^{-1} (B) 
= \vert  k_{\Delta}^{\times} \vert \langle \chi_1, \chi_2 \rangle = 0$ car $\chi_1 \neq \chi_2$, donc 
$S_1 = 0$.\\

Fixons à présent $Y = \left( \begin{array}{cc}
             A & B\\
             C & D
             \end{array} \right) \in X_2$. Alors :
$$
Y 
= \left( \begin{array}{cc}
             A/C & \alpha^2\\
             1 & A/C
             \end{array} \right) 
\left( \begin{array}{cc}
             C & D - A \frac{DA-BC}{A^2 - \alpha^2 C^2}\\
             0 & C \frac{BC-AD}{\alpha^2 C^2 - A^2}
             \end{array} \right) 
=  \left( \begin{array}{cc}
             \frac{BC-AD}{C} & A\\
             0 & C
             \end{array} \right) w 
\left( \begin{array}{cc}
             1 & D/C\\
             0 & 1
             \end{array} \right)
$$
Donc :
$$
\widetilde{\varphi}_{s_0}^{0, f/2} (Y) 
= \chi_1^{-1} (C) \chi_2^{-1} (C) \chi_2^{-1} (BC-AD) \chi_2 (\alpha^2 C^2 - A^2)
$$
et :
$$
J.\widetilde{f} (Y) 
= \frac{F(w)}{\vert \overline{G}_{s_0} \vert} 
\chi_1 (BC-AD) \chi_1^{-1} (C) \chi_2 (C)
$$
On en déduit que :
\begin{eqnarray*}
S_2 
& = & \sum_{Y \in X_2} \widetilde{\varphi}_{s_0}^{0, f/2} (Y) 
   J. \widetilde{f} (Y)\\
& = & \frac{F(w)}{\vert \overline{G}_{s_0} \vert} \times  
   \sum_{A, C \in k_{\Delta}^{\times}} \chi_1^{-1} (C^2) \chi_2 (\alpha^2 C^2 - A^2) \times 
   (\sum_{\substack{B, D \in k_{\Delta} \\ AD-BC \neq 0}} \chi_1 (BC-AD) \chi_2^{-1} (BC-AD))
\end{eqnarray*}
Soient $A$ et $C$ dans $k_{\Delta}^{\times}$. Alors :
$$
\sum_{\substack{B, D \in k_{\Delta} \\ AD-BC \neq 0}} \chi_1 (BC-AD) \chi_2^{-1} (BC-AD) 
= \sum_{s \in k_{\Delta}^{\times}} 
\sum_{\substack{B, D \in k_{\Delta} \\ AD-BC = s}} \chi_1 (s) \chi_2^{-1} (s)
$$
Or l'ensemble 
$\{ (B,D) \in k_{\Delta} \times k_{\Delta} : AD-BC = s \}$ est de cardinal $Q$. 
Ainsi :
$$
\sum_{\substack{B, D \in k_{\Delta} \\ AD-BC \neq 0}} \chi_1 (BC-AD) \chi_2^{-1} (BC-AD) 
= Q \times \sum_{s \in k_{\Delta}^{\times}} \chi_1 (s) \chi_2^{-1} (s) 
= \langle \chi_1, \chi_2 \rangle = 0
$$
Finalement, $S_2 = 0$.

\end{itemize}
D'après les calculs précédents, on a 
$\varphi_{s_0}^{0, f/2} (J. \widetilde{f}) 
= S_1 + S_2 = 0 
= \zeta^{-1}$, 
ce qui est impossible car \mbox{$\zeta^e = \chi (\varpi_{\mathbb{K}}) \neq 0$}.

\end{démo}

On utilise les notations de \ref{NotationHypotheseCNDistinction}, \ref{NotationsS0Steinberg}. 
Le résultat suivant découle directement de \ref{SyntheseCNDistinctionTotRam}, 
\ref{NonDistinctionCasfPair} et \ref{ResultatMult1Induite} :

\begin{theo}\label{ConclusionCNDistinction}
Soit $(\pi, V)$ une représentation membre de la série discrète de $G$, de niveau $0$, 
non cuspidale et $(\mathbb{K}_f, \chi_f)$ la paire admissible modérée associée à $\pi$ 
définie en \ref{ParamSilbergerZinkNonCuspidales}. On définit comme dans 
\ref{DefinitionSystCoeffShneiderEtStuhler} un système de coefficients sur $X_{\mathbb{K}}$ d'espace vectoriel $V$. 
 Soit $\varphi = (\varphi_s)_{s \in X_0}$ dans 
\mbox{${\rm ker} (\partial_1^{\ast}) \cap {\rm Hom}_{\mathcal{D}^{\times}} (C_0, \mathds{1})$}.\\
Si $f$ est pair, alors 
$\pi$ n'est pas $\mathcal{D}^{\times}$-distinguée.\\
Supposons que $\pi$ soit $\mathcal{D}^{\times}$-distinguée. Alors $f$ est impair, 
$\overline{\chi}_0$ est non trivial sur $k^{\times}$ mais trivial sur les carrés de 
$k^{\times}$, pour tout sommet $s$, pour tout $0 \leq \nu_1 < \nu_2 \leq f-1$, on a :
$$
\varphi_s^{\nu_1, \nu_2} = 0
$$
et pour tout $\nu \in \{ 0, \cdots, f-1 \}$, pour tout sommet $s$ distinct de $s_0$, on a :
$$
\varphi_{s}^{\nu} (h) 
= \frac{(-1)^{k-1} (Q+1)}{2 Q^{k}} \times 
(\widetilde{\varphi}_{s_0}^{\nu} (\delta_1^{s_0}) 
- \widetilde{\varphi}_{s_0}^{\nu} (\delta_{Q+1}^{s_0})) \times h(\delta_{Q+1}^{s}) 
\, \, \text{pour tout} \, \,  
h \in W_{s}^{\nu} 
$$
où $k = d (s, s_0)$, $\delta_{Q+1}^{s}$ est la droite correspondant à l'arête 
$\{ \widetilde{s}, s \}$, avec $[ s_0, \cdots, \widetilde{s}, s]$ chemin géodésique reliant $s$ 
à $s_0$. 
 
\end{theo}

\subsection{Distinction et correspondance de Jacquet-Langlands.}\label{DistinctionEtJLNonCusp}

\begin{nota}
Soit $n = 2 \delta$, soit $f$ un diviseur de $\delta$. 
On note $r = (2 \delta)/f$, alors $n = r \times f$. 
Soit $\Delta$ une $\mathbb{K}$-algèbre à division centrale d'indice $\delta$. 
On suppose que l'extension quadratique $\mathbb{K} / \mathbb{F}$ est totalement ramifiée, modérément ramifiée. 
On fixe $\varpi_{\mathbb{K}}$ et $\varpi_{\mathbb{F}}$ 
des uniformisantes de $\mathbb{K}$ et $\mathbb{F}$ 
telles que 
$\varpi_{\mathbb{K}}^2 = \varpi_{\mathbb{F}}$. 
On note $JL$ la correspondance de Jacquet-Langlands :
$$
JL : \mathcal{R}_0^2 ({\rm GL}_2 (\Delta)) 
\rightarrow \mathcal{R}_0^2 ({\rm GL}_n (\mathbb{K}))  
$$
Soit $(\chi_f, \mathbb{K}_f)$ une paire admissible modérée (i.e $\mathbb{K}_f$ est une extension non ramifiée de 
degré $f$ de $\mathbb{K}$ et $\chi_f$ est un caractère modéré, $\mathbb{K}$-régulier, de $\mathbb{K}_f^{\times}$). 
\end{nota}

\begin{nota}
Soit $\rho \in \mathcal{R}_0^2 ({\rm GL}_2 (\Delta))$ de paire admissible modérée associée 
$(\chi_f, \mathbb{K}_f)$. Alors d'après \ref{ConclusionCNDistinction}, la représentation 
$\rho$ n'est pas $\mathcal{D}^{\times}$-distinguée. 
Nous noterons $\Pi = JL (\rho)$. Alors $\Pi$ est une série discrète de niveau $0$ 
de paire admissible associée 
$(\chi_f, \mathbb{K}_f)$ (cf. \cite{SilbergerZink2}). 
\end{nota}

Commençons par quelques rappels sur les membres de la série discrète de ${\rm GL}_n (\mathbb{K})$ 
(on pourra retrouver ces résultats dans \cite{BushnellHenniartArticle}, section 3). 
Introduisons quelques notations. 
On fixe $e$ un diviseur de $n$ et $b$ un entier naturel tel que $n = e \times b$. 
Soit $P$ le sous-groupe parabolique standard (triangulaire supérieur par blocs) de 
${\rm GL}_n (\mathbb{K})$ de radical unipotent $N$ tel que 
$P/N = L \simeq {\rm GL}_b (\mathbb{K})^{\times e}$. 
Si $\tau$ est une représentation cuspidale de ${\rm GL}_b (\mathbb{K})$ on note, pour tout 
réel $a$, $\tau^a$ la représentation de ${\rm GL}_b (\mathbb{K})$ définie par :
$$
\tau^a : {\rm GL}_b (\mathbb{K}) \rightarrow V_{\tau}, 
x \mapsto \vert \vert {\rm det} (x) \vert \vert_{\mathbb{K}}^a \tau (x)
$$
et $\tau_L$ la représentation de $L$ définie par 
$\tau_L = \tau^{(1-e)/2} \otimes \tau^{(3-e)/2} \otimes \cdots 
\otimes \tau^{(e-1)/2}$. 
Enfin, on note $I_P^{{\rm GL}_n (\mathbb{K})} \tau_L$ l'induite parabolique normalisée de $\tau_L$. 
On a les propriétés suivantes :

\begin{theo}
\begin{itemize}
 \item[i)] La représentation $I_P^{{\rm GL}_n (\mathbb{K})} \tau_L$ admet un unique quotient 
irréductible. On note ${\rm St}_e (\tau)$, et on appelle Steinberg généralisée de ${\rm GL}_n (\mathbb{K})$ 
de support cuspidal $\tau^{\otimes e}$, cette représentation. 
Alors ${\rm St}_e (\tau)$ est membre de la série discrète de ${\rm GL}_n (\mathbb{K})$.
 \item[ii)] Si $\rho$ est membre de la série discrète de ${\rm GL}_n (\mathbb{K})$, il existe un diviseur $e$ de 
$n$ et une représentation cuspidale $\tau$ de ${\rm GL}_{n/e} (\mathbb{K})$ tels que 
$\rho \simeq {\rm St}_e (\tau)$. De plus, le couple $(e, \tau)$ est entièrement déterminé par 
la classe d'isomorphisme de la représentation $\rho$.
\end{itemize}

\end{theo}

Ainsi, la représentation $\Pi$ est isomorphe à une représentation de Steinberg généralisée. 
Dans leur article \cite{BushnellHenniartArticle}, C. J. Bushnell et G. Henniart font le lien 
entre la paire admissible modérée associée à $\Pi$ et le support cuspidal de 
cette Steinberg généralisée :

\begin{theo}
(\cite{BushnellHenniartArticle}, 6.3, démonstration du théorème 2) 
La représentation $\Pi$ est une Steinberg généralisée de support cuspidal 
$\pi_f^{\otimes r}$ où $\pi_f$ est une représentation cuspidale de niveau $0$ 
de ${\rm GL}_f(\mathbb{K})$. 
Soit $\zeta$ le caractère quadratique non ramifié de $\mathbb{K}_f^{\times}$ 
(i.e. $\zeta$ est trivial sur 
$\mathcal{O}_{\mathbb{K}_{f/2}}^{\times}$ et d'ordre $2$, en particulier  
$\zeta (\varpi_{\mathbb{K}}) = -1$). 
Alors la paire admissible modérée associée à $\pi_f$ est 
$(\zeta \chi_f, \mathbb{K}_f)$.
\end{theo}

\begin{rmq}
Nous avons choisi de travailler avec les conventions de Silberger et Zink, c'est-à-dire que 
la paire admissible modérée associée à une représentation $\rho$,  
membre de la série discrète de niveau~$0$, est préservée par Jacquet-Langlands, 
mais a besoin d'être corrigée pour connaître la paire admissible modérée 
associée au type étendu maximal de 
niveau $0$ de $\rho$ (cf. remarque \ref{LienEntreLesParametrisations}). Dans leur article 
\cite{BushnellHenniartArticle}, Bushnell et Henniart ont choisi de paramétrer les représentations 
par les paires admissibles modérées associées à leurs types étendus maximaux. 
Ainsi, d'après \cite{BushnellHenniartArticle}, 6.3, démonstration du théorème 2, 
$\Pi$ est une Steinberg généralisée de support cuspidal 
$\pi_f^{\otimes r}$ où la paire admissible modérée associé au type étendu maximal de 
$\pi_f$ est $(\zeta \chi_f, \mathbb{K}_f)$. Or, d'après \cite{SilbergerZink2}, exemples 0.10 page 184, 
puisque $\pi_f$ est cuspidale, la paire admissible modérée associée à $\pi_f$ est aussi 
$(\zeta \chi_f, \mathbb{K}_f)$ (il n'y a pas de correction à faire).
\end{rmq}

\begin{nota}
On notera $\eta$ le caractère quadratique de $\mathbb{F}^{\times}$ associé à 
$\mathbb{K} / \mathbb{F}$ (i.e. $\eta$ est l'unique caractère non 
trivial de $\mathbb{F}^{\times}$ qui soit trivial sur 
${\rm N}_{\mathbb{K} / \mathbb{F}} (\mathbb{K}^{\times})$). 
\end{nota}

\begin{rmq}
Soit $\psi$ un générateur du groupe de Galois 
${\rm Gal} (\mathbb{K} / \mathbb{F})$, alors 
$\psi (\varpi_{\mathbb{K}}) = - \varpi_{\mathbb{K}}$ et 
${\rm N}_{\mathbb{K} / \mathbb{F}} (\varpi_{\mathbb{K}}) 
= \varpi_{\mathbb{K}} \psi (\varpi_{\mathbb{K}}) 
= - \varpi_{\mathbb{K}}^2 = - \varpi_{\mathbb{F}}$.
On en déduit que 
$\eta (\varpi_{\mathbb{F}}) = \eta (-\varpi_{\mathbb{F}}) \eta (-1) = \eta (-1)$.
\end{rmq}

Nous aurons besoin du résultat suivant, du à N. Matringe, 
qui donne un critère de distinction d'une représentation membre de 
la série discrète de ${\rm GL}_n (\mathbb{K})$ en fonction de son support cuspidal :

\begin{theo}
 (\cite{Matringe}, Corollaire 4.2)
La représentation $\Pi$ est ${\rm GL}_n (\mathbb{F})$-distinguée si et seulement si 
$\pi_f$ est $\eta$-distinguée.
\end{theo}

On vérifie facilement la propriété suivante :

\begin{lem}
On fixe $\alpha$ une racine carrée de $\eta (-1)$. 
Pour $x$ dans $\mathcal{O}_{\mathbb{K}}^{\times}$, on note $\overline{x}$ sa réduction dans 
$k^{\times} = \frac{\mathcal{O}_{\mathbb{K}}^{\times}}{1+\mathcal{P}_{\mathbb{K}}} 
\simeq \frac{\mathcal{O}_{\mathbb{F}}^{\times}}{1+\mathcal{P}_{\mathbb{F}}}$. On définit 
$\widehat{\eta} : \mathbb{K}^{\times} \rightarrow \mathbb{C}^{\times}$ tel que pour tout 
$i$ dans $\mathbb{Z}$ et tout $x$ dans $\mathcal{O}_{\mathbb{K}}^{\times}$ :
$$
\widehat{\eta} (\varpi_{\mathbb{K}}^i x) = \alpha^i \eta (x) 
= \alpha^i \overline{\eta} (\overline{x})
$$
Alors $\widehat{\eta}$ est un caractère de $\mathbb{K}^{\times}$ qui prolonge $\eta$. 
\end{lem}

\begin{nota}
Nous noterons $\widehat{\pi}_f = \widehat{\eta} \otimes \pi_f$. 
\end{nota}

\begin{lem}
On a les isomorphismes de $\mathbb{C}$-espaces vectoriels suivants :
$$
{\rm Hom}_{{\rm GL}_f (\mathbb{F})} (\pi_f, \eta \circ {\rm det}) 
= {\rm Hom}_{{\rm GL}_f (\mathbb{F})} (\pi_f, \widehat{\eta} \circ {\rm det}) 
\simeq {\rm Hom}_{{\rm GL}_f (\mathbb{F})} (\widehat{\eta} \otimes \pi_f, \mathds{1})
$$
Ainsi $\Pi$ est ${\rm GL}_n (\mathbb{F})$-distinguée si et seulement 
si $\widehat{\pi}_f$ est ${\rm GL}_f (\mathbb{F})$-distinguée.

\end{lem}

\begin{rmq}
 On vérifie aussi que $\widehat{\pi}_f$ est une représentation cuspidale de niveau $0$ de 
${\rm GL}_f (\mathbb{K})$ de paire admissible modérée associée 
$(\zeta \times \chi_f \times \widehat{\eta} \circ {\rm N}_{\mathbb{K}_f / \mathbb{K}}, \mathbb{K}_f)$. 
\end{rmq}

Notre objectif ici est de montrer que si on se trouve dans le cas du théorème \ref{ConclusionCNDistinction} 
où la série discrète paramétrée par $(\chi_f, \mathbb{K}_f)$ n'est pas $\mathcal{D}^{\times}$-distinguée, alors son 
image par la correspondance de Jacquet-Langlands n'est pas non plus ${\rm GL}_n (\mathbb{F})$-distinguée.

\subsubsection{Premier cas : cas où $f$ est pair.}

On suppose dans tout ce paragraphe que $f$ est pair. 
Nous rappelons un résultat du théorème 3.2.31 dans \cite{Coniglio} :

\begin{theo}\label{RappelsConditionsDistinctionCuspidalesTotRamfPair}
Soit $\Lambda$ une représentation cuspidale de niveau $0$ de 
${\rm GL}_f (\mathbb{K})$ (où $f$ est pair) 
de paire admissible modérée $(\theta, \mathbb{K}_f)$ 
(ici $\mathbb{K} / \mathbb{F}$ est totalement ramifiée, modérément 
ramifiée et les uniformisantes de $\mathbb{F}$ et $\mathbb{K}$ vérifient l'égalité 
$\varpi_{\mathbb{K}}^2 = \varpi_{\mathbb{F}}$). Alors 
$\Lambda$ est ${\rm GL}_f (\mathbb{F})$-distinguée si et seulement si $\theta$ est trivial sur $\mathbb{F}^{\times}$, 
$\overline{\theta}$ est trivial sur $k_{f/2}^{\times}$ (où $k_{f/2}$ est une extension de 
degré $f/2$ de $k$) et, pour $k_f = k_{f/2} [\delta]$, on a 
$\theta (\varpi_{\mathbb{K}}) \overline{\theta} (\delta) = -1$.
\end{theo}

\begin{propo}\label{JacquetLanglandsConditionsDistCasfpairSeriesDiscretes}
Supposons que $\pi_f$ est $\eta$-distinguée. 
Soit $\delta$ dans $k_f^{\times}$ tel que $\delta \notin k_{f/2}$ et $\delta^2 \in k_{f/2}^{\times}$. 
Alors $\zeta \chi_f$ est trivial sur $\mathbb{F}^{\times}$, 
$\overline{\chi}_f$ et $\overline{\zeta \chi_f}$ sont 
triviaux sur $k_{f/2}^{\times}$, $\overline{\chi}_f (\delta) = -1$, et 
$\chi_f (\varpi_{\mathbb{K}}) = -1$.
\end{propo}

\begin{démo}
Il suffit d'appliquer le théorème \ref{RappelsConditionsDistinctionCuspidalesTotRamfPair} au caractère 
$\theta 
= \zeta \times \chi_f \times \widehat{\eta} \circ {\rm N}_{\mathbb{K}_f / \mathbb{K}}$. 

\end{démo}

D'après ce qui précède, si la représentation $\pi_f$ est $\eta$-distinguée, alors 
la paire admissible modérée associée à $\pi_f$ vérifie les conditions du théorème 
\ref{RappelsConditionsDistinctionCuspidalesTotRamfPair}. Le résultat suivant, du à 
A. C. Kable, va nous permettre de conclure :

\begin{theo}\label{HakimMurnDistinctionEtaDistinction}
(\cite{Kable}) 
Une représentation cuspidale de ${\rm GL}_f (\mathbb{K})$ ne peut pas être à la fois 
${\rm GL}_f (\mathbb{F})$-distinguée et $\eta$-distinguée.
\end{theo}

On en déduit directement le théorème suivant :

\begin{theo}\label{JacquetLanglandsNonDistinctionfPairSeriesDiscretes}
 La représentation $\Pi = JL (\rho)$ n'est pas ${\rm GL}_n (\mathbb{F})$-distinguée. 
Ainsi $\rho$ n'est pas $\mathcal{D}^{\times}$-distinguée et son image par la correspondance de Jacquet-Langlands 
n'est pas ${\rm GL}_n (\mathbb{F})$-distinguée. 
\end{theo}

\subsubsection{Deuxième cas : cas où $f$ est impair et $\overline{\chi}_0$ est trivial sur $k^{\times}$.}

On suppose dans tout ce paragraphe que $f$ est impair. 
On rappelle un résultat du théorème 3.1.21 dans \cite{Coniglio} :

\begin{theo}\label{RappelsConditionsDistinctionCuspidalesTotRamfImpair}
Soit $\Lambda$ une représentation cuspidale de niveau $0$ de ${\rm GL}_f (\mathbb{K})$. 
Puisque $f$ est impair et l'extension $\mathbb{K} / \mathbb{F}$ est totalement ramifiée, modérément 
ramifiée, la représentation 
$\Lambda$ n'est pas  
${\rm GL}_f (\mathbb{F})$-distinguée.
\end{theo}

On en déduit directement le résultat suivant :

\begin{theo}\label{JacquetLanglandsNonDistinctionfImpairSeriesDiscretes}
 La représentation $\Pi = JL (\rho)$ n'est pas ${\rm GL}_n (\mathbb{F})$-distinguée. 
Ainsi $\rho$ n'est pas $\mathcal{D}^{\times}$-distinguée et son image par la correspondance de Jacquet-Langlands 
n'est pas ${\rm GL}_n (\mathbb{F})$-distinguée. 
 \end{theo}

\section{Cas particulier de la représentation de Steinberg de ${\rm GL}_2 (\Delta)$.}

\subsection{Paramétrisation de Silberger et Zink de la représentation de Steinberg de $G$.}

Commençons par donner une définition de la représentation de Steinberg 
de $G = {\rm GL}_2 (\Delta)$ :

\begin{déf}
Notons $P_0$ le sous-groupe parabolique minimal de 
$G$ :
$$
P_0 = \left(\begin{array}{cc} 
             \Delta & \Delta\\
             0 & \Delta^{\times}
            \end{array} \right)
$$
On note $V_0 = {\rm Ind}_{P_0}^{G} \mathds{1}_{P_0}$. Pour tout parabolique $P$ de 
$G$ tel que $P_0 \subseteq P$, on note $V_P$ l'image canonique de 
${\rm Ind}_{P}^{G} \mathds{1}_{P}$ dans $V_0$. 
Alors :
$$
{\rm St}_G = V_0 / (\sum_{P_0 \varsubsetneq P} V_P)
$$
Ici, on a donc :
${\rm St}_G = V_0 / \mathds{1}_G$.
\end{déf}

Rappelons quel est le paramètre de Silberger et Zink de la représentation de Steinberg de 
$G$ :

\begin{propo}
(\cite{SilbergerZink2}, exemple 0.10).\\
Notons $(\pi, V)$ la représentation de Steinberg de $G$. 
On utilise les mêmes notations que dans \ref{ParamSilbergerZinkNonCuspidales}. 
Soit $\chi$, caractère modéré de $\mathbb{K}_{d}^{\times} = \mathbb{K}_{2 \delta}^{\times}$, qui paramétrise en partie 
$\pi$ (alors $\chi_{\vert \mathbb{K}^{\times}}$ est le caractère central de $\pi$) et $f$ la longueur de 
la ${\rm Gal} (k_{\mathbb{K}, d} / k_{\mathbb{K}})$-orbite de $\overline{\chi}$. Alors  
$f = 1$, $\chi_f$ est le caractère trivial de $\mathbb{K}^{\times}$, $\overline{\chi}_0$ est le caractère trivial de 
$k_{\Delta}^{\times}$ et donc pour tout sommet $s$ :
$$
V_s = {\rm St}_{\overline{G}_s}
$$
\end{propo}

\subsection{Cas où l'extension $\mathbb{K} / \mathbb{F}$ est totalement ramifiée.}

D'après \ref{SyntheseCNDistinctionTotRam}, on a: 

\begin{theo}\label{TheoremeNonDistinctionSteinbergTotRam}
On suppose que l'extension $\mathbb{K} / \mathbb{F}$ est totalement ramifiée, modérément ramifiée.
Dans ce cas, la représentation de Steinberg $(\pi, V)$ de $G$ ne peut pas être 
$\mathcal{D}^{\times}$-distinguée.
\end{theo}

\subsection{Cas où l'extension $\mathbb{K} / \mathbb{F}$ est non ramifiée.}

\begin{nota}
Rappelons que $A_{\mathbb{K}}$ désigne l'appartement standard de 
$X_{\mathbb{K}}$ tel que :
$$
\mathcal{S}_0 
= \{ s_k = [ \mathcal{O}_{\mathbb{K}} \oplus \mathcal{P}_{\mathbb{K}}^k] 
: k \in \mathbb{Z} \}
$$
est l'ensemble des sommets de $A_{\mathbb{K}}$ et $j(X_{\mathbb{F}}) = m_0$ 
est le milieu de l'arête $\{ s_0, s_1 \}$.
\end{nota}

\begin{lem}
Soit $\mathcal{A}_2$ l'ordre héréditaire principal minimal de ${\rm M}_2 (\Delta)$ :
$$
\mathcal{A}_2 = \left( \begin{array}{cc}
                        \mathcal{O}_{\Delta} & \mathcal{O}_{\Delta}\\
                        \mathcal{P}_{\Delta} & \mathcal{O}_{\Delta}
                       \end{array} \right)
$$
et $\mathcal{R}_2$ le normalisateur de $\mathcal{A}_2^{\times}$ dans $G$ 
$\mathcal{R}_2 = \langle t_2 \rangle \mathcal{A}_2^{\times}$ 
où $t_2 = \left( \begin{array}{cc}
                        0 & 1\\
                        \varpi_{\Delta} & 0
                       \end{array} \right)$.
Le type étendu de niveau $0$ maximal pour $\pi = {\rm St}_G$ est la représentation 
$(\widetilde{\Sigma}, \mathbb{C})$ de $\mathcal{R}_2$ telle que :
$$ 
\widetilde{\Sigma} (t_2^i x) : \mathbb{C} \rightarrow \mathbb{C}, 
v \mapsto (-1)^i v
$$
pour $i$ dans $\mathbb{Z}$ et $x$ dans $\mathcal{A}_2^{\times}$. 
\end{lem}

\begin{lem}
Soit $a_0$ l'arête reliant les sommets $s_0$ et $s_1$.
Pour tout $v$ dans $V_{a_0}$, on a :
$$
\pi (\varpi_{\mathcal{D}}).v = -v
$$
\end{lem}

\begin{démo}
On remarque que 
$\varpi_{\mathcal{D}}^d = \varpi_{\mathbb{F}} = \varpi_{\mathbb{K}} = 
\varpi_{\Delta}^{\delta} = (t_2^2)^{\delta} = t_2^d$. 
Ainsi $\varpi_{\mathcal{D}} = t_2 u$ pour un élément $u$ dans ${\rm GL}_2 (\mathcal{O}_{\Delta})$. 
Donc, pour tout $v$ dans $V_{a_0}$ :
$\pi (\varpi_{\mathcal{D}}).v = \pi (t_2).v = \widetilde{\Sigma} (t_2).v = -v$. 

\end{démo}

\begin{lem}
On a $\varphi_{s_0} = 0$. 
\end{lem}

\begin{démo}
Rappelons que, d'après \ref{OrbitesCasNonRam}, les 
$\mathcal{D}^{\times}$-orbites des sommets de $X_{\mathbb{K}}$ sont exactement les sphères de centre 
$m_0$, où $m_0 = j ([\mathcal{O}_{\mathcal{D}} ])$ est le milieu de l'arête $\{ s_0, s_1 \}$.
Ainsi, le seul sommet voisin de $s_0$ qui est aussi dans sa $\mathcal{D}^{\times}$-orbite est $s_1$. 
Notons $a_0$ l'arête $\{ s_0, s_1 \}$. On a 
$\varpi_{\mathcal{D}} = t_2 u$
où $u \in {\rm GL}_2 (\mathcal{O}_{\Delta})$. 
Puisque $t_2$ échange les sommets $s_0$ et $s_1$, et $u$ fixe $s_0$, 
$\varpi_{\mathcal{D}}$ échange les sommets $s_0$ et $s_1$. 
Sur le module de Jacquet $V_{a_0} = (V_{s_0})^{\mathcal{U}_{a_0}^1}$, on a :
$$
\varphi_{s_0} = \varphi_{s_1} = \varphi_{\varpi_{\mathcal{D}}.s_1} \circ \pi (\varpi_{\mathcal{D}}) 
= \varphi_{s_0} \circ \pi (\varpi_{\mathcal{D}})
$$
D'après le lemme précédent, pour tout $v$ dans $V_{a_0}$, on a  
$\pi (\varpi_{\mathcal{D}}).v = -v$. 
Ainsi, pour tout $v$ dans $V_{a_0}$, 
$\varphi_{s_0} (v - \pi (\varpi_{\mathcal{D}}).v) = 0 = 2 \varphi_{s_0} (v)$, donc 
$\varphi_{s_0}$ est triviale sur $V_{a_0}$.\\
Notons $\delta_0$ la droite de $\mathbb{P}_{s_0}^1 (k_{\Delta})$ correspondant à l'arête 
$a_0$ et $\delta_1, \cdots, \delta_{Q^2}$ les droites correspondant aux autres voisins de $s_0$. 
D'après ce qui précède, $\varphi_{s_0} = 0$ sur $V_{a_0} = V_{\delta_0}$. 
Comme précédemment, on fixe $\widetilde{\varphi}_{s_0}$ dans 
$\{ \mathbb{P}_{s_0}^1 (k_{\Delta}) \rightarrow \mathbb{C} \} / \{ \text{fonctions constantes} \}$ 
telle que, pour tout $h$ dans $V_{s_0}$ :
$$
\varphi_{s_0} (h) = 
\sum_{d \in \mathbb{P}_{s_0}^1 (k_{\Delta})} \widetilde{\varphi}_{s_0} (d) h(d)
$$
Fixons $h$ dans $V_{s_0}$ telle que $h (\delta_0) \neq 0$. 
On définit $h_{\delta_0}$ dans $V_{\delta_0}$ de sorte que 
$h_{\delta_0} (d) = - \frac{h (\delta_0)}{Q^2}$ si $d \neq \delta_0$ et 
$h_{\delta_0} (\delta_0) = h (\delta_0)$. Alors :
$$
\varphi_{s_0} (h_{\delta_0}) = 0 
= \widetilde{\varphi}_{s_0} (\delta_0) h (\delta_0)
- \frac{h(\delta_0)}{Q^2} \sum_{\delta \neq \delta_0} \widetilde{\varphi}_{s_0} (\delta)
$$
D'où 
$\widetilde{\varphi}_{s_0} (\delta_0) 
= \frac{1}{Q^2} \sum_{i=1}^{Q^2} \widetilde{\varphi}_{s_0} (\delta_i)$. 
Les $Q^2$ autres sommets voisins de $s_0$ sont tous à une distance $\frac{1}{2} + 1$ de $m_0$, 
donc sont dans la même $\mathcal{D}^{\times}$-orbite. 
Notons $t_1, \cdots, t_{Q^2}$ ces voisins ($t_i$ correspond à la droite $\delta_i$). 
Soit $d_i$ dans $\mathcal{D}^{\times}$ tel que 
$d_i.t_i = t_1$. 
On vérifie facilement que l'on 
peut supposer que $d_i \in \mathcal{O}_{\mathcal{D}}^{\times}$. 
Ainsi, pour tout $i$ dans 
$\{ 1, \cdots, Q^2 \}$, il existe $g_i$ dans $\mathcal{O}_{\mathcal{D}}^{\times}$ tel que 
$g_i.t_i = t_1$ (et alors $g_i.\delta_i = \delta_1$ et $g_i.s_0 = s_0$). 
Alors,  
$\varphi_{s_0} = \varphi_{g_i.s_0} \circ \pi (g_i) = \varphi_{s_0} \circ \pi (g_i)$. 
On en déduit alors que pour $i$ dans $\{ 1, \cdots, Q^2 \}$, 
$\widetilde{\varphi}_{s_0} (\delta_i) = \widetilde{\varphi}_{s_0} (\delta_1)$. 
Par conséquent :
$$
\widetilde{\varphi}_{s_0} (\delta_0) 
= \frac{1}{Q^2} \sum_{i=1}^{Q^2} \widetilde{\varphi}_{s_0} (\delta_i) 
= \widetilde{\varphi}_{s_0} (\delta_0) 
= \frac{Q^2}{Q^2} \widetilde{\varphi}_{s_0} (\delta_1) = \widetilde{\varphi}_{s_0} (\delta_1)
$$
Par suite, $\widetilde{\varphi}_{s_0}$ est constante et 
$\varphi_{s_0} = 0$. 
Puis, en utilisant la $\mathcal{D}^{\times}$-équivariance de $\varphi$, on montre que 
si $s$ est un sommet de $X_{\mathbb{K}}$ tel que $d(s, m_0) = 1/2$, 
alors $\varphi_s = 0$.

\end{démo}

Un raisonnement par récurrence nous permet de montrer le lemme suivant :

\begin{lem}
On a $\varphi_{s_k} = 0$ pour tout entier naturel $k$. 
\end{lem}

Une conséquence immédiate de ce lemme est le théorème suivant :

\begin{theo}\label{TheoremeNonDistinctionSteinbergNonRam}
On suppose que l'extension $\mathbb{K} / \mathbb{F}$ est non ramifiée.
Dans ce cas, la représentation de Steinberg $(\pi, V)$ de $G$ n'est pas  
$\mathcal{D}^{\times}$-distinguée.
\end{theo}

\subsection{Représentation de Steinberg et correspondance de Jacquet-Langlands.}

Dans cette partie, on suppose seulement que l'extension $\mathbb{K} / \mathbb{F}$ 
est quadratique modérément ramifiée (elle peut être non ramifiée ou totalement ramifiée). 

\begin{rmq}
Comme en \ref{DistinctionEtJLNonCusp}, on note $JL$ la correspondance de Jacquet-Langlands :
$$
JL : \mathcal{R}_0^2 ({\rm GL}_2 (\Delta)) 
\rightarrow \mathcal{R}_0^2 ({\rm GL}_n (\mathbb{K}))  
$$
Notons $\Pi$ la représentation de Steinberg de ${\rm GL}_2 (\Delta)$. Alors, 
d'après \cite{SilbergerZink2}, $JL(\Pi)$ est la représentation de Steinberg 
de ${\rm GL}_n (\mathbb{K})$ où $n = 2 \delta$ (et $\delta$ est l'indice de $\Delta$). 
\end{rmq}

Une conséquence immédiate du corollaire 4.2 de \cite{Matringe} est le fait suivant :

\begin{theo}\label{NonDistinctionSteinbergEtJacquetLanglands}
La représentation de Steinberg $JL(\Pi)$ de ${\rm GL}_n (\mathbb{K})$ n'est pas 
${\rm GL}_n (\mathbb{F})$-distinguée.\\
Ainsi $\Pi$ n'est pas $\mathcal{D}^{\times}$-distinguée et son image par 
la correspondance de Jacquet-Langlands n'est pas ${\rm GL}_n (\mathbb{F})$-distinguée.
\end{theo}

\begin{démo}
D'après le corollaire 4.2 de \cite{Matringe}, la représentation de Steinberg 
$JL(\Pi)$ de ${\rm GL}_n (\mathbb{K})$ est 
${\rm GL}_n (\mathbb{F})$-distinguée si et seulement si le caractère trivial est 
$\eta^{n-1}$-distingué où $\eta$ est le caractère quadratique de $\mathbb{F}^{\times}$ 
associé à $\mathbb{K} / \mathbb{F}$. Or, ici $n$ est pair et le caractère trivial n'est pas $\eta$-distingué. 
On en déduit que $JL(\Pi)$ n'est pas 
${\rm GL}_n (\mathbb{F})$-distinguée.

\end{démo}

\appendix
\section{Un résultat de multiplicité $1$ en caractéristique nulle}

\subsection{Rappels sur les paires de Gelfand-Kazhdan.}

Dans cette partie, on note $G$ un groupe topologique localement compact, totalement discontinue, 
unimodulaire. Soit $\mathbb{F}$ un corps local non archimédien. 

\begin{déf}
Soit $H$ un sous-groupe fermé de $G$. 
On suppose que $H \backslash G$ est muni d'une mesure invariante à droite, qu'il existe une 
anti-involution continue $i$ de $G$ telle que $i(H) = H$ et que toute distribution bi-$H$-invariante de 
$G$ est fixée par $i$. On dit alors que $(G, H)$ (ou plus précisément $(G, H, i)$) 
est une paire de Gelfand-Kazhdan.
\end{déf}

\begin{déf}
 Soit $H$ un sous-groupe fermé de $G$. 
On dit que le couple $(G, H)$ est une paire de Gelfand si pour toute représentation lisse 
irréductible $\pi$ de $G$ on a :
$$
{\rm dim} ({\rm Hom}_H (\pi, \mathds{1})) \times 
{\rm dim} ({\rm Hom}_H (\widetilde{\pi}, \mathds{1})) \leq 1
$$
(où $\widetilde{\pi}$ désigne la contragrédiente (lisse) de $\pi$).
\end{déf}

Le résultat suivant est démontré dans un cas particulier dans \cite{GelfandKazhdan2} et 
plus généralement dans \cite{Gross} (proposition 4.2) :

\begin{theo}
Une paire de Gelfand-Kazhdan est en particulier une paire de Gelfand.
\end{theo}

Dans la suite, notre objectif sera de montrer qu'un certain couple $(G, H)$ est une paire de Gelfand. 
Pour cela, nous utiliserons un résultat plus général, 
donné dans \cite{GelfandKazhdan2} : 

\begin{déf}
 Soit $\mathbb{X}$ une variété algébrique définie sur $\mathbb{F}$. 
Soit $\mathbb{G}$ un groupe algébrique défini sur $\mathbb{F}$ qui agit sur $\mathbb{X}$ de manière algébrique. 
On note $X = \mathbb{X} (\mathbb{F})$ et $G = \mathbb{G} (\mathbb{F})$. 
On dit que la paire $(\mathbb{X}, \mathbb{G})$ est régulière si pour tout $x$ dans $X$ 
l'application :
$$
\mathbb{G} \rightarrow \mathbb{G}.x, g \mapsto g.x
$$
est une submersion.
\end{déf}

\begin{rmq}
Pour tout $x$ dans $X$, il existe au plus (à multiplication par un réel strictement positif près) 
une mesure positive $\mu_x$ sur $G.x$ qui soit $G$-invariante. 
\end{rmq}

\begin{theo}\label{TheoremeConditionsPairesGelfandKazhdan}
Soit $(\mathbb{X}, \mathbb{G})$ une paire régulière. 
On suppose qu'il existe $\sigma$ un automorphisme algébrique de $\mathbb{X}$, 
défini sur $\mathbb{F}$, tel que :
\begin{itemize}
 \item[1.] Pour tout $x$ dans $X$, on a $\sigma (G.x) = G.x$.
 \item[2.] Pour tout $x$ dans $\mathbb{X}$, on a $\sigma (\mathbb{G}.x) = \mathbb{G}.x$. 
 \item[3.] L'automorphisme $\sigma$ fixe les mesures $\mu_x$.
\end{itemize}
Alors $\sigma$ fixe toute distribution $G$-invariante sur $X$.
\end{theo}

\begin{rmq}
 Si l'automorphisme $\sigma$ est d'ordre fini, alors la troisième condition est satisfaite.
\end{rmq}

\subsection{Théorème de Hilbert $90$ généralisé.}

On introduit les notations suivantes pour toute la suite :

\begin{nota}
 On considère $\mathbb{K} / \mathbb{F}$ une extension quadratique séparable de corps locaux non 
archimédiens. On suppose que $\mathbb{F}$ est de caractéristique nulle. 
On note $\Gamma = {\rm Gal} (\mathbb{K} / \mathbb{F}) = \{ Id, \sigma_0 \}$. 
Soit $A$ une $\mathbb{F}$-algèbre centrale simple d'indice $n$. 
Soit $\mathcal{D}$ une $\mathbb{F}$-algèbre à division centrale d'indice $d$ telle que:
$$
A \simeq {\rm M}_m (\mathcal{D})
$$
On note $l = n^{2} = m^{2} \times d^{2}$. 
On fixe $(b_1, \cdots, b_l)$ une base de $A$ sur $\mathbb{F}$. 
On pose $A_{\mathbb{K}} = A \otimes_{\mathbb{F}} \mathbb{K}$. 
On a l'injection canonique de $\mathbb{K}$ dans $A_{\mathbb{K}}$:
$$
\mathbb{K} \hookrightarrow A_{\mathbb{K}}, \, k \mapsto Id_{A} \otimes k
$$
$A_{\mathbb{K}}$ est une $\mathbb{K}$-algèbre centrale simple.
On considère $\mathcal{B} = (b_1 \otimes 1, \cdots, b_l \otimes 1)$ une $\mathbb{K}$-base de 
$A_{\mathbb{K}}$. 
Soit $G = A_{\mathbb{K}}^{\times}$.

\end{nota}

\begin{rmq}
On remarque que $\Gamma$ agit naturellement sur $A_{\mathbb{K}}$:
$$
\sigma.(a \otimes k) = a \otimes \sigma (k) ; a \in A, k \in \mathbb{K}
$$
Donc pour tout $\mu_1, \cdots, \mu_l \in \mathbb{K}$, on a:
$$
\sigma.(\mu_1 (b_1 \otimes 1) + \cdots + \mu_l (b_l \otimes 1)) 
= \mu_1^{\sigma} (b_1 \otimes 1) + \cdots +  \mu_l^{\sigma} (b_l \otimes 1)
$$
Par conséquent, si $c \in A_{\mathbb{K}}$ est donné par ses coordonnées dans la base~$\mathcal{B}$, 
$c = (\mu_1, \cdots, \mu_l)$, alors pour tout $\sigma \in \Gamma$, on a:
$$
\sigma (c) = (\mu_1^{\sigma}, \cdots, \mu_l^{\sigma})
$$
On a donc une action naturelle de $\Gamma$ sur $G$. 
\end{rmq}

Pour démontrer le théorème de Hilbert $90$ dans notre cas, nous commençons par rappeler 
un résultat d'indépendance algébrique des automorphismes 
(cf. \cite{Bourbaki1}, chapitre 5 (\textit{Corps commutatifs}), paragraphe 10 (\textit{Extensions galoisiennes}), 
partie 7 (\textit{Indépendance algébrique des automorphismes}), Théorème 4):

\begin{theo}
Soit $\mathbb{E}$ un corps infini, $\mathbb{L}$ une extension galoisienne de $\mathbb{E}$ de degré fini~$m$.
Soit $\lbrace \sigma_1, \cdots, \sigma_m \rbrace = {\rm Gal} (\mathbb{L} / \mathbb{E})$. 
Soit $\Omega$ une extension quelconque de $\mathbb{L}$. 
Soit $P \in \Omega [X_1, \cdots, X_m]$ tel que pour tout $x \in \mathbb{L}$, 
$P (\sigma_1 (x), \cdots, \sigma_m (x) )~=~0$. 
Alors $P = 0$.
\end{theo}

On montre le corollaire suivant par récurrence :

\begin{cor}
Soit $\mathbb{E}$ un corps infini, $\mathbb{L}$ une extension galoisienne de $\mathbb{E}$ de degré fini~$m$.
Soit $\lbrace \sigma_1, \cdots, \sigma_m \rbrace = {\rm Gal} (\mathbb{L} / \mathbb{E})$. 
Soient $N \geq 1$ et $N.m$ variables indépendantes $Z_1, \cdots, Z_{N.m}$. 
Soit $P$ dans $\mathbb{L} [Z_1, \cdots, Z_{N.m}]$ tel que, pour tout 
$x_1, \cdots, x_N \in \mathbb{L}$, on a :
$$
P (x_1^{\sigma_1}, \cdots, x_1^{\sigma_m}, x_2^{\sigma_1}, \cdots, x_2^{\sigma_m}, \cdots, x_N^{\sigma_m}) = 0
$$
Alors $P = 0$. 
\end{cor}

\begin{theo}
(Théorème de Hilbert $90$ généralisé)\\
Tout $1$-cocycle de $\Gamma$ dans $G$ est un $1$-cobord, i.e. :
$$
H^1 (\Gamma, G) = \{ \overline{\mathds{1}} \}
$$
\end{theo}

\begin{démo}
Considérons 
$\varphi : \Gamma \rightarrow G, \sigma \mapsto \varphi (\sigma)  = g_{\sigma}$
un $1$-cocycle. Alors, pour tout $\sigma, \tau \in \Gamma$ :
$$
\varphi (\sigma \tau) = \varphi (\sigma) . \sigma (\varphi (\tau))
$$
Montrons que $\varphi$ est un $1$-cobord, i.e montrons qu'il existe $b$ dans $G$ 
tel que, pour tout $\sigma$ dans $\Gamma$ :
$$
\varphi (\sigma) = b^{-1} \sigma (b)
$$ 
Ainsi, pour tout $b$ dans $G$, on a $\varphi (Id) = b^{-1} Id (b)$. 
Puisque $\varphi (Id) = 1$ et $\Gamma = \{ Id, \sigma_0 \}$, il suffit de trouver $b$ dans $G$ tel que 
$\varphi (\sigma_0) = b^{-1} \sigma_0 (b)$. 
On notera $\varphi (\sigma_0) = g_{\sigma_0} = (\lambda_1, \cdots, \lambda_l)$ (sa décomposition dans $\mathcal{B}$). 
Pour tout $c \in A_{\mathbb{K}}$, on forme la série de Poincaré:
$$
b = \sum_{\sigma \in \Gamma} \varphi (\sigma) \sigma (c) = c + g_{\sigma_0} . \sigma_0 (c)
$$
Soit $\tau_0 \in \Gamma$, alors:
\begin{eqnarray*}
\tau_0 (b)
& = & \tau_0 (\sum_{\sigma \in \Gamma} \varphi (\sigma) \sigma (c)) 
    =  \sum_{\sigma \in \Gamma} \tau_0 (\varphi(\sigma)) \tau_0 ( \sigma (c)) 
    = \sum_{\sigma \in \Gamma} \varphi (\tau_0)^{-1} \varphi (\tau_0 \sigma) \tau_0 \sigma (c)\\
& = & \varphi (\tau_0)^{-1} (\sum_{\sigma \in \Gamma} \varphi (\tau_0 \sigma) \tau_0 \sigma (c)) 
    = \varphi (\tau_0)^{-1} b
\end{eqnarray*}
Pour montrer que $\varphi$ est un $1$-cobord, il suffit de montrer que l'on peut choisir $c$ 
tel que $b$ soit inversible, i.e tel que ${\rm Nrd}_{A_{\mathbb{K}}} (b) \neq 0$ 
(où ${\rm Nrd}_{A_{\mathbb{K}}}$ désigne la 
norme réduite de $A_{\mathbb{K}}$).\\
Raisonnons par l'absurde.
Supposons que pour tout $c \in A_{\mathbb{K}}$ on ait ${\rm Nrd}_{A_{\mathbb{K}}} (b) = 0$. 
Soient $\mu_1, \cdots, \mu_l$ dans~$\mathbb{K}$, soit $c = (\mu_1, \cdots, \mu_l)$, alors :
$$
{\rm Nrd}_{A_{\mathbb{K}}} (c + g_{\sigma_0} . \sigma_0 (c)) = 0
$$
avec $\sigma_0 (c) = (\mu_1^{\sigma_0}, \cdots, \mu_l^{\sigma_0})$ et $g_{\sigma_0} = (\lambda_1, \cdots, \lambda_l)$. 
Et donc:
$$
g_{\sigma_0} . \sigma_0 (c) = g_{\sigma_0} . (\sum_{i} \mu_i^{\sigma_0} (b_i \otimes 1)) 
= \sum_{i} \mu_i^{\sigma_0} 
(g_{\sigma_0} (b_i \otimes 1)) 
$$
La norme réduite étant polynomiale en les coordonnées, il existe $P$ 
dans $\mathbb{K} [Z_1, \cdots, Z_{2 \times l}]$ tel que :
$$
{\rm Nrd}_{A_{\mathbb{K}}} (c + g_{\sigma_0} . \sigma_0 (c)) 
= P (\mu_{1}, \mu_{2}, \cdots, \mu_{l}, \mu_{1}^{\sigma_0},\cdots, \mu_{l}^{\sigma_0})
$$
On en déduit que, pour tout $\mu_1, \cdots, \mu_l \in \mathbb{K}$ :
$$
P (\mu_{1}, \mu_{2}, \cdots, \mu_{l}, \mu_{1}^{\sigma_0}, \cdots, \mu_{l}^{\sigma_0}) = 0
$$
D'après le corollaire précédent, $P = 0$ dans $\mathbb{K} [Z_1, \cdots, Z_{2 \times l}]$. 
Ainsi, pour tout $x_1, \cdots, x_{2 \times l} \in \mathbb{K}$, on a:
$$
P (x_1, \cdots, x_{2 \times l}) = 0
$$
On choisit $x_{1}, x_{2}, \cdots, x_{l} \in \mathbb{K}$ tels que $c_0 = (x_1, \cdots, x_l) \in G $. 
On pose $y_{1} = y_{2} = \cdots = y_{l} = 0$. 
Alors $P (x_{1}, x_{2}, \cdots, x_{l}, y_{1}, \cdots, y_{l}) = 0 $. 
De plus, $P (x_{1}, x_{2}, \cdots, x_{l}, y_{1}, \cdots, y_{l})$ correspond à:
$$
{\rm Nrd}_{A_{\mathbb{K}}} (c_0 + 0) = {\rm Nrd}_{A_{\mathbb{K}}} (c_0) 
$$
Or ${\rm Nrd}_{A_{\mathbb{K}}} (c_0) \neq 0$ puisque $c_0 \in G$, ce qui est absurde puisque $P = 0$. 
On en déduit qu'il existe $c \in A_{\mathbb{K}}$ tel que $b = c + g_{\sigma_0} . \sigma_0 (c) \in G$. 
Par suite $\varphi$ est un $1$-cobord et on a bien le résultat.

\end{démo}

\subsection{Résultat de multiplicité $1$.}

\begin{theo}
La paire $(G, H)$ est une paire de Gelfand-Kazhdan.
\end{theo}

\begin{démo}
Pour montrer ce résultat, bien connu dans le cas déployé, 
nous allons utiliser le théorème \ref{TheoremeConditionsPairesGelfandKazhdan} 
en adaptant la preuve de \cite{Serre} à notre contexte. 
\begin{itemize}
 \item[1.] Par définition on a :
$$
G = \mathbb{G} (\mathbb{F}) = \mathbb{H} (\mathbb{K} \otimes_{\mathbb{F}} \mathbb{F})
 \simeq \mathbb{H} (\mathbb{K}) = ({\rm M} _m (\mathcal{D}) \otimes_{\mathbb{F}} \mathbb{K} )^{\times}
 \simeq ({\rm M} _m (\mathcal{D}\otimes_{\mathbb{F}} \mathbb{K}))^{\times} 
$$
De plus, on sait qu'il existe une $\mathbb{K}$-algèbre à division centrale $\Delta$ telle que :
$$
{\rm M} _m (\mathcal{D}\otimes_{\mathbb{F}} \mathbb{K})
 \simeq {\rm M} _{\mu} (\Delta)
$$
Donc $G \simeq {\rm GL} _{\mu} (\Delta)$. 
On a également:
$$
H = \mathbb{H} (\mathbb{F}) = (A \otimes_{\mathbb{F}} \mathbb{F} )^{\times} \simeq A^{\times}
$$
d'où $H \simeq {\rm GL} _m (\mathcal{D})$. 
On fixe $\sigma$ un générateur du groupe de Galois ${\rm Gal} (\mathbb{K} / \mathbb{F})$. 
Alors $\sigma$ agit naturellement sur les coefficients de 
${\rm M} _m (\mathcal{D}\otimes_{\mathbb{F}} \mathbb{K})$ :
$$
\sigma . (d \otimes k) = d \otimes \sigma (k)
$$
On notera $\theta$ l'action de $\sigma$ sur ${\rm M} _m (\mathcal{D}\otimes_{\mathbb{F}} \mathbb{K})$ 
et donc sur $G$. 
On remarque en particulier que $\theta$ est une involution et que $H$ est l'ensemble 
des points fixes de $G$ sous l'action de $\theta$. 
On rappelle que l'on a un ismorphisme de $\mathbb{F}$-algèbres (cf. \cite{Bourbaki2}, 
chapitre 5 (\textit{Corps commutatifs}), paragraphe 10 (\textit{Extensions galoisiennes}), 
partie 4 (\textit{Descente galoisienne}), Proposition 8 et remarque 1) :
$$
\alpha : \mathbb{K} \otimes_{\mathbb{F}} \overline{\mathbb{F}} 
\rightarrow \overline{\mathbb{F}} \times \overline{\mathbb{F}} , \, 
x \otimes y \mapsto (\sigma (x) y, x y)
$$
On a 
$H_1 = \mathbb{H} (\overline{\mathbb{F}}) 
= ({\rm M} _m (\mathcal{D}) \otimes_{\mathbb{F}} \overline{\mathbb{F}})^{\times}
= ({\rm M} _m (\mathcal{D} \otimes_{\mathbb{F}} \overline{\mathbb{F}}))^{\times}$. 
Il existe une extension $\mathbb{F} \subseteq \mathbb{L} \subseteq \overline{\mathbb{F}}$ 
telle~que:
$$
\mathcal{D} \otimes_{\mathbb{F}} \mathbb{L} \simeq {\rm M}_d (\mathbb{L})
$$
Par conséquent 
$\mathcal{D} \otimes_{\mathbb{F}} \overline{\mathbb{F}} \simeq {\rm M}_d (\overline{\mathbb{F}})$
d'où 
$H_1 \simeq ({\rm M} _m ({\rm M}_d (\overline{\mathbb{F}})))^{\times}$. 
Par suite, $H_1 \simeq {\rm GL}_n (\overline{\mathbb{F}})$. 
En utilisant l'isomorphisme $\alpha$, on remarque de plus que:
\begin{eqnarray*}
G_1 = \mathbb{G} (\overline{\mathbb{F}})
& = & \mathbb{H} (\mathbb{K} \otimes_{\mathbb{F}} \overline{\mathbb{F}})
   =  (A \otimes_{\mathbb{F}} (\mathbb{K} \otimes_{\mathbb{F}} \overline{\mathbb{F}}) )^{\times}\\
& \simeq & (A \otimes_{\mathbb{F}} (\overline{\mathbb{F}} \oplus \overline{\mathbb{F}}) )^{\times} 
  \simeq  ((A \otimes_{\mathbb{F}} \overline{\mathbb{F}}) 
\oplus (A \otimes_{\mathbb{F}} \overline{\mathbb{F}}))^{\times} 
  \simeq  ({\rm M}_n (\overline{\mathbb{F}}) \times {\rm M}_n (\overline{\mathbb{F}}))^{\times}\\
& \simeq & {\rm GL}_n (\overline{\mathbb{F}}) \times {\rm GL}_n (\overline{\mathbb{F}})
\end{eqnarray*}
On remarque également que:
$$
G_1 \simeq ({\rm M} _m (\mathcal{D} \otimes_{\mathbb{F}} \overline{\mathbb{F}}) \otimes_{\mathbb{F}} \mathbb{K})^{\times} 
\simeq {\rm GL}_n (\mathbb{K} \otimes_{\mathbb{F}} \overline{\mathbb{F}})
$$
On a donc $G_1 \simeq H_1 \times H_1$ avec:
$$
\widetilde{\alpha} : G_1 \rightarrow H_1 \times H_1 , \, 
[x_{i,j} \otimes y_{i,j} ]_{i,j} \mapsto ( [\sigma (x_{i,j}) y_{i,j} ]_{i,j}, [x_{i,j} y_{i,j} ]_{i,j})
$$
où $G_1$ est identifié à ${\rm GL}_n (\mathbb{K} \otimes_{\mathbb{F}} \overline{\mathbb{F}})$, 
$x_{i,j} \in \mathbb{K}$, $y_{i,j} \in \overline{\mathbb{F}}$. 
Notons encore $\theta$ l'action naturelle de $\sigma$ sur les coefficients de $G_1$. 
Via ces identifications, l'action de $\sigma$ sur $H_1 \times H_1$ (encore notée $\theta$), 
correspond à la permutation des coordonnées. 
On va désormais considérer deux actions.\\ 
L'une de $H_1 \times H_1$ sur $G_1$:
$$
h_1, h_2 \in H_1, g \in G_1, (h_1, h_2) . g = h_1 g h_2^{-1}
$$
L'autre de $H \times H$ sur $G$:
$$
h_1, h_2 \in H, g \in G, (h_1, h_2) . g = h_1 g h_2^{-1}
$$
On notera également:
$$
i : G_1 \rightarrow G_1, g \mapsto \theta (g)^{-1}
$$
(et on notera de même $i : G \rightarrow G, g \mapsto \theta (g)^{-1}$). 

\item[2.] D'après ce qui précède, $i$ est une anti-involution continue et $i(H) = H$.

\item[3.] D'après \cite{Borel}, chapitre II, paragraphe 6, proposition 6.7, puisque 
$\mathbb{F}$ est de caractéristique nulle, pour tout $x$ dans ${\rm GL}_n (\mathbb{F})$, 
l'application :
$$
\chi_x : H_1 \times H_1 \rightarrow H_1 \times H_1.x
$$
est submersive. 
On en déduit que $(\mathbb{G}, \mathbb{H} \times \mathbb{H})$ est une paire régulière. 

\item[4.] Puisque $i$ est une anti-involution, $i$ est d'ordre $2$. 
Ainsi, pour tout $x$ dans $G$, $i$ fixe toute mesure positive $\mu_x$ sur 
$H \times H.x$ qui soit $H \times H$-invariante.

\item[5.] Montrons que $i$ fixe les $H_1 \times H_1$-orbites de $G_1$. \\
On peut tout d'abord remarquer que $H_1$ est l'ensemble des points fixes de $G_1$ 
sous l'action de $\theta$, donc pour tout $h \in H_1$, on a $i (h) = h^{-1}$. 
De plus, on a l'injection de $H_1$ dans $G_1$ par injection diagonale:
$$
H_1 \rightarrow G_1 \simeq H_1 \times H_1, x \mapsto (x, x)
$$
On remarque que pour tout $g \in G_1$ et tous $h_1, h_2 \in H_1$, on a:
$$
i ( (h_1, h_2).g) = i (h_1 g h_2^{-1}) = i (h_2^{-1}) i(g) i(h_1) = h_2 i(g) h_1^{-1}
$$
Il nous suffit donc de montrer que $g$ et $i(g)$ sont dans la même $H_1 \times H_1$-orbite. 
Grâce à l'application $\widetilde{\alpha}$ définie auparavant, on a une bijection entre $G_1$ et 
$H_1 \times H_1$. 
On peut alors voir l'action de $H_1 \times H_1$ sur $G_1$ comme une action de $H_1 \times H_1$ sur lui même :
$$
(h_1, h_2).(x_1, x_2) = (h_1 x_1 h_2^{-1}, h_1 x_2 h_2^{-1})
$$
On notera $( H_1 \times H_1 / \sim)$ ses orbites. 
D'après le calcul précédent, l'application $i$ permet de définir une application :
$$
j : (H_1 \times H_1 / \sim) \rightarrow (H_1 \times H_1 / \sim), 
\overline{g} = \overline{(x_1, x_2)} \mapsto \overline{i (g)} = \overline{(x_2^{-1}, x_1^{-1})}
$$
Montrons que $j$ est l'application identité. 
L'application $H_1 \times H_1 \rightarrow H_1, \, (a,b) \mapsto a b^{-1}$ 
permet d'identifier les $H_1 \times H_1$-orbites de $G_1$ aux 
classes de similitudes de $H_1$. On peut ainsi définir l'application :
$$
\eta :  (H_1 \times H_1 / \sim) \rightarrow  Ad (H_1) \backslash H_1 ,  
\overline{(x_1, x_2)} \mapsto Ad (H_1). (x_1 x_2^{-1})
$$
où $Ad (H_1) \backslash H_1$ désigne les classes de similitudes de $H_1$.
On vérifie rapidement que $\eta$ est une bijection. 
Soit $\zeta : Ad (H_1) \backslash H_1 \rightarrow Ad (H_1) \backslash H_1$ telle que 
$\zeta = \eta \circ j \circ \eta^{-1}$.  
Pour tout $h \in H_1$, il existe $x_1, x_2$ dans $H_1$ tels que $Ad (H_1). h = Ad (H_1). (x_1 x_2^{-1})$ et :
$$
\zeta (Ad (H_1). (x_1 x_2^{-1})) 
= \zeta \circ \eta (\overline{(x_1, x_2)}) 
= \eta \circ j (\overline{(x_1, x_2)}) = Ad (H_1) . (x_2^{-1} x_1)
$$
Comme $x_1 x_2^{-1} = x_2 (x_2^{-1} x_1) x_2^{-1}$, $x_1 x_2^{-1}$ et $x_2^{-1} x_1$ sont $H_1$-conjugués. Par 
conséquent $\zeta$ est l'application identité.
On en déduit que $i$ fixe bien les $H_1 \times H_1$-orbites de $G_1$. 

\item[6.] Montrons que $i$ fixe les $H \times H$-orbites de $G$. \\
On pose:
$$
S = \{ g \in G : i (g) = g \} = \{ g \in G : \theta (g) = g^{-1} \}
$$
On vérifie facilement que pour tout $g$ dans $G$, 
$g i(g) \in S$ et pour tout $h$ dans $H$, $ghi(gh) = gi(g)$. 
Considérons alors l'application :
$$
\psi : G / H \rightarrow S, gH \mapsto g i(g)
$$
Montrons que $\psi$ est surjective.  
Soit $s \in S$. On définit l'application :
$$
\varphi_s : \Gamma \rightarrow G, \, \sigma \mapsto s , \, Id \mapsto 1
$$
On vérifie facilement que $\varphi_s$ est un $1$-cocycle. 
D'après le théorème de Hilbert $90$ généralisé, $\varphi_s$ est un $1$-cobord. 
Il existe donc $g \in G$ tel que, pour tout $\tau$ dans $\Gamma$, 
$\varphi (\tau) = g \tau (g)^{-1}$.
En particulier, pour $\tau = \sigma$, on a:
$$
\varphi (\sigma) = g \sigma (g)^{-1} \Leftrightarrow s = g i(g) = \psi (g)
$$
Donc $\psi$ est bien surjective. 
On remarque de plus que pour tout $g \in G$ et tout $h \in H$, on a :
$$
\psi (hg) = h g i(g) i(h) = h (\psi (g))h^{-1} = Ad (h).\psi (g)
$$
L'application suivante est donc bien définie:
$$
\xi : H \textbackslash G / H \rightarrow Ad(H) \textbackslash S, \, 
HgH \mapsto Ad (H) . (g i(g))
$$
Par surjectivité de $\psi$, $\xi$ est clairement surjective. 
Montrons qu'elle est injective. 
Soient $g_1, g_2 \in G$ tels que $\xi (H g_1 H) = \xi (H g_2 H)$, alors 
$Ad (H) . (g_1 i(g_1)) = Ad (H) . (g_2 i(g_2))$. 
Il existe donc $h \in H$ tel que :
$$
g_1 i(g_1) = h g_2 i(g_2) h^{-1} = (h g_2 h^{-1}) (i (h g_2 h^{-1}))
$$
d'où
$(h g_2 h^{-1})^{-1} g_1 = i (g_1^{-1} (h g_2 h^{-1}))$. 
Posons $x = g_1^{-1} (h g_2 h^{-1})$, alors $i (x) = \theta (x)^{-1} = x^{-1}$ donc 
$x = \theta (x)$ et $x \in H$. 
Par suite $g_1^{-1} (h g_2 h^{-1}) \in H$ et $h g_2 h^{-1} \in g_1 H$ donc 
$H g_1 H = H g_2 H$. 
L'application $\xi$ est donc bien injective. 
On remarque que:
$$
i (h_1 g h_2) = h_2^{-1} i(g) h_1^{-1}
$$
pour tous $h_1, h_2$ dans $H$ et tout $g$ dans $G$.  
L'application suivante est donc bien définie:
$$
j : H \textbackslash G / H \rightarrow H \textbackslash G / H, H g H \mapsto i (H g H) = H i(g) H
$$
Montrons que l'application $j$ est l'identité. 
L'application $\xi$ étant bijective, on peut définir:
$$
\Lambda = \xi \circ j \circ \xi^{-1} : Ad(H) \textbackslash S \rightarrow Ad(H) \textbackslash S
$$
Pour montrer que $j$ est l'identité, il suffit de montrer que $\Lambda$ est l'identité. 
On remarque que pour tout $s \in S$ il existe $g \in G$ tel que $s = g i(g)$ et:
$$
\Lambda (Ad (H) . s) = \Lambda ( \xi (H g H)) = \xi \circ j (H g H) 
= \xi (H i(g) H) = Ad (H) . (i(g) g)
$$
On a donc:
$$
\Lambda : Ad(H) \textbackslash S \rightarrow Ad(H) \textbackslash S, \, 
Ad (H) . (g i(g)) \mapsto Ad (H) . (i(g) g)
$$
Montrons que pour tout $g \in G$, $g i(g)$ et $i(g) g$ sont conjugués dans $H$. 
Or, $i(g) g = g^{-1} (g i(g)) g$ donc ces deux éléments sont conjugués dans $G$. 
Pour montrer notre résultat, il suffit de vérifier que deux éléments de $S$ conjugués dans $G$ 
sont conjugués dans $H$. 
Soient $s_1$ et $s_2 \in S$ tels qu'il existe $g \in G$ tel que 
$s_1 = g^{-1} s_2 g$, alors $g s_1 = s_2 g$. 
Comme $g s_1 = s_2 g$, on a 
$(\lambda g) s_1 = s_2 (\lambda g)$ pour tout $\lambda \in \mathbb{K}^{\times}$ 
(car $\mathbb{K}^{\times}$ est contenu dans le centre de $A_{\mathbb{K}}$). 
Soit $\lambda$ dans $\mathbb{K}^{\times}$. On applique $\theta$ à l'égalité précédente :
$$
\theta (\lambda g) s_1^{-1} = s_2^{-1} \theta (\lambda g) \, \, \text{donc} \, \, 
\theta (\lambda g) s_1 = s_2 \theta (\lambda g)
$$
On en déduit que pour tout $\lambda \in \mathbb{K}^{\times}$, 
$(\lambda g + \theta (\lambda g)) s_1 = s_2 (\lambda g + \theta (\lambda g))$. 
Comme $\lambda g + \theta (\lambda g)$ est invariant sous l'action de $\theta$,  
$\lambda g + \theta (\lambda g) \in {\rm M}_m (\mathcal{D})$. 
Puisque 
$\lambda g + \theta (\lambda g) 
= \lambda^{\sigma} \theta (g) ( g^{-1} \theta (g) + \frac{\lambda}{\lambda^{\sigma}} Id)$, 
il suffit de trouver $\lambda \in \mathbb{K}^{\times}$ 
tel que $g^{-1} \theta (g) + \frac{\lambda}{\lambda^{\sigma}} Id$ soit inversible. 
On fixe $\mathbb{L} / \mathbb{K}$ une extension telle que $A_{\mathbb{K}}$ se 
déploie sur $\mathbb{L}$ :
$$
A_{\mathbb{K}} \otimes_{\mathbb{K}} \mathbb{L} \simeq {\rm M}_p (\mathbb{L})
$$
Alors $A_{\mathbb{K}} \hookrightarrow {\rm M}_p (\mathbb{L})$ et pour tout $x \in A_{\mathbb{K}}$, 
${\rm Nrd}_{A_{\mathbb{K}}} (x) = {\rm det}_{{\rm M}_p (\mathbb{L})} (x)$ et 
$x \in A_{\mathbb{K}}^{\times}$ si et seulement si ${\rm det}_{{\rm M}_p (\mathbb{L})} (x) \neq 0$. 
Soit $\lambda \in \mathbb{K}^{\times}$, alors :
$$
{\rm Nrd}_{A_{\mathbb{K}}} (g^{-1} \theta (g) + \frac{\lambda}{\lambda^{\sigma}} Id) 
= {\rm det}_{{\rm M}_p (\mathbb{L})} (g^{-1} \theta (g) + \frac{\lambda}{\lambda^{\sigma}} Id) 
= \chi_{g^{-1} \theta (g)} (- \frac{\lambda}{\lambda^{\sigma}})
$$
Donc $g^{-1} \theta (g) + \frac{\lambda}{\lambda^{\sigma}} Id$ est inversible 
si et seulement si $- \frac{\lambda}{\lambda^{\sigma}}$ n'est pas dans le spectre de $g^{-1} \theta (g)$. 
Par théorème de Hilbert $90$, on a un isomorphisme :
$$
\mathbb{K}^{\times} / \mathbb{F}^{\times} \rightarrow {\rm ker} ({\rm N}_{\mathbb{K} / \mathbb{F}}), 
\, \lambda \mapsto \frac{\lambda}{\lambda^{\sigma}}
$$
On vérifie facilement que $\mathbb{K}^{\times} / \mathbb{F}^{\times}$ est infini 
en l'identifiant à $\mathbb{P}^{1} (\mathbb{F})$, l'ensemble des droites de $\mathbb{F}^{2}$. 
On en déduit le résultat, c'est-à-dire que deux éléments de $S$ conjugués dans $G$ sont conjugués dans~$H$.

\end{itemize}
D'après le théorème \ref{TheoremeConditionsPairesGelfandKazhdan}, $i$ fixe toute 
distribution bi-$H$-invariante sur $G$. Par conséquent, $(G, H)$ est bien une paire de 
Gelfand-Kazhdan, et donc une paire de Gelfand.

\end{démo}


\vspace{4mm}
\noindent 
Charlène Coniglio-Guilloton, 
Charlene.Coniglio@math.univ-poitiers.fr\\
Département de Mathématiques, 
Téléport 2 - BP 30179, 
Boulevard Marie et Pierre Curie \\
86962 Futuroscope Chasseneuil Cedex, 
France

\end{document}